\documentclass[10pt]{amsart}
\usepackage{amsfonts}
\usepackage{graphicx}
\usepackage{amscd}
\usepackage{amsmath,amsfonts,amssymb,amsthm,mathrsfs,xypic,cite,color}
\xyoption{curve}
\usepackage{fancybox}
\newcommand{\m}{\Lambda}
\newcommand{\Hom}{\operatorname{Hom}}
\newcommand{\End}{\operatorname{End}}
\newcommand{\Ker}{\operatorname{Ker}}
\newcommand{\K}{\operatorname{K}}
\newcommand{\co}{\operatorname{Ck}}
\newcommand{\cok}{\operatorname{Cok}}
\newcommand{\Cok}{\operatorname{Coker}}
\newcommand{\Ima}{\operatorname{Im}}

\newcommand{\ha}{\operatorname{\mathcal{A}}}
\newcommand{\x}{\operatorname{\mathscr{X}}}
\newcommand{\D}{{\rm D}}

\newcommand {\hp}{\mathcal P}
\newcommand{\s}{\hfill \blacksquare}
\newtheorem{thm}{Theorem}[section]

\newtheorem{cor}[thm]{Corollary}

\newtheorem{lem}[thm]{Lemma}

\newtheorem{exm}[thm]{Example}

\newtheorem{facts}[thm]{Facts}

\newtheorem{prop}[thm]{Proposition}
\newtheorem{rem}[thm]{Remark}
\newtheorem{defn}[thm]{Definition}
\newcommand{\lra}{\longrightarrow}
\newcommand{\lraf}[1]{\stackrel{#1}{\lra}}
\begin{document}

\title [Separated monic representations II]
{Separated monic representations II: \\ Frobenius subcategories and RSS
equivalences}
\author [Pu Zhang, Bao-Lin Xiong] {Pu Zhang, Bao-Lin Xiong}
\thanks{Supported by the NSFC 11271251 and 11431010.}
\thanks{pzhang$\symbol{64}$sjtu.edu.cn \ \ \ xiongbaolin@gmail.com}

\begin{abstract} \ This paper aims at looking for Frobenius subcategories, via the separated
monomorphism category ${\rm smon}(Q, I, \x)$;  and on the other hand, to establish
an {\rm RSS} equivalence from ${\rm smon}(Q, I, \x)$ to its dual ${\rm sepi}(Q, I, \x)$. For a bound quiver $(Q, I)$ and an algebra $A$, where $Q$ is acyclic and $I$ is generated by
monomial relations, let $\Lambda=A\otimes_k kQ/I$. For any
additive subcategory $\x$ of $A$-mod, we construct ${\rm smon}(Q, I, \x)$ combinatorially. This construction describe Gorenstein-projective $\m$-modules as $\mathcal {GP}(\m)
= {\rm smon}(Q, I, \mathcal {GP}(A))$. It admits
a homological interpretation, and
enjoys a reciprocity ${\rm smon}(Q, I, \  ^\bot T)= \ ^\bot (T\otimes
kQ/I)$ for a cotilting
$A$-module $T$. As an application, ${\rm smon}(Q, I, \x)$ has
Auslander-Reiten sequences if $\x$ is resolving and contravariantly
finite with $\widehat{\x}=A$-mod. In particular, ${\rm smon}(Q, I,
A)$ has Auslander-Reiten
sequences. It also admits
a filtration interpretation as ${\rm smon}(Q, I, \mathscr{X})={\rm Fil}(\mathscr{X}\otimes \mathcal P(kQ/I))$, provided that $\x$ is extension-closed.  As an application, ${\rm smon}(Q, I, \x)$ is an extension-closed Frobenius subcategory if and only if so
is $\x$. This gives ``new" Frobenius subcategories of $\m$-mod in the sense that
they are not $\mathcal{GP}(\m)$. Ringel-Schmidmeier-Simson equivalence ${\rm smon}(Q, I,
\x)\cong{\rm sepi}(Q, I, \x)$ is introduced and the existence is proved for arbitrary extension-closed subcategories $\x$. In particular, the
Nakayama functor $\mathcal N_\m$ gives an {\rm RSS} equivalence ${\rm smon}(Q, I,
A)\cong{\rm sepi}(Q, I, A)$ if
and only if $A$ is Frobenius.
For a chain $Q$ with
arbitrary $I$, an explicit formula of an {\rm RSS} equivalence is
found for arbitrary additive subcategories $\x$.

\vskip5pt

\end{abstract}
\maketitle \vskip-30pt\section {\bf Introduction and preliminaries}
\subsection {} Throughout, $A$ is a
finite-dimensional algebra over a field $k$, $Q$ a finite acyclic
quiver, $I$ an admissible ideal of the path algebra $kQ$ generated
by monomial relations, and $\m\colon= A\otimes \ kQ/I$. Let $A$-mod (resp. mod$A$) be the
category of finitely generated left (resp. right) $A$-modules,  and $\mathcal {P}(A)$ (resp. $\mathcal
{GP}(A)$) the subcategory of $A$-mod of projective (resp. Gorenstein-projective) modules.
Tensor
product $\otimes$ is over $k$ if not specified, and subcategories
are full subcategories closed under isomorphisms.
If no otherwise stated, $\x$ is an arbitrary additive subcategory of $A$-mod.

\vskip5pt

\subsection {} This paper arises from looking for the extension-closed Frobenius subcategories $\mathcal F$ of $\m$-mod, with the canonical exact structure ([Q], [K1]).
Such an $\mathcal F$ induces an algebraic triangulated category ([H], [MSS]). There are several approaches to deal with this question (see e.g. [Buch], [H], [K2], [C2], [KLM2], [IT], [KIWY]).
Since $\mathcal {GP}(\m)$ is the largest extension-closed Frobenius subcategories $\mathcal F$ with $\mathcal P(\mathcal F) \subseteq \mathcal P(A)$ (see Proposition \ref{frobenius1}),
and $\mathcal {GP}(\m)$ can be described as ${\rm smon}(Q, I, \mathcal {GP}(A))$, we deal with this question via
the separated monomorphism category ${\rm smon}(Q, I, \x)$ and the correspondence $\x\mapsto {\rm smon}(Q, I, \x)$.

\vskip5pt

To know when ${\rm smon}(Q, I, \x)$ is Frobenius, the key is to know when it have enough injective objects and enough projective objects.
The latter is not hard with a direct argument; however,  it is difficult to prove that it has enough injective objects if so has $\x$.
This also leads the study of an {\rm RSS} equivalence. This difficulty is overcame after a filtration interpretation of ${\rm smon}(Q, I, \x)$
as ${\rm Fil}(\x\otimes \mathcal P(kQ/I))$ (Theorem \ref{filt}).  Thus, ${\rm smon}(Q, I, \x)$ is an extension-closed Frobenius subcategory if and only if so
is $\x$ (Corollary \ref{frobenius2}). This gives ``new" Frobenius subcategories of $\m$-mod, in the sense that they are not $\mathcal {GP}(\m)$, as indicated in Example \ref{exmfrobenius}.

\vskip5pt

\subsection {} In studying the representations of the tensor product $A\otimes B$ of $k$-algebras $A$ and $B$, the category $A\mbox{-}{\rm mod}\otimes B\mbox{-}{\rm mod}$ of
$L\otimes U$ with $L\in A\mbox{-}{\rm mod}$ and $U\in B\mbox{-}{\rm
mod}$ usually is properly contained in $(A\otimes B)\mbox{-}{\rm mod}$, and we have the Cartan-Eilenberg isomorphism ([CE, Thm. 3.1, p.209, p.205]):
$${\rm Ext}^{m}_{A\otimes B}(L\otimes
U, M\otimes V)\cong\bigoplus\limits_{p+q = m}({\rm Ext}^p_A(L,
M)\otimes {\rm Ext}^q_B(U, V)), \ \forall \ m\ge 0.$$
An advantage of taking $B= kQ/I$ is that the representations of bound quiver $(Q, I)$ ([R], [ARS],
[ASS]) can be applied to study $\m$-mod ([RS1-RS3], [S1-S3], [KLM1, KLM2], [RZ]).
This choice of $\m = A\otimes kQ/I$ is not restricted in the sense that, principally speaking,
any algebra is of this form.

\vskip5pt

{\it A
representation $X$ of $(Q, I)$ over $A$ } is a datum $X=(X_i,
X_{\alpha}, i\in Q_0, \alpha\in Q_1)$, where $X_i\in A$-mod,
and $X_\alpha: X_{s(\alpha)}\rightarrow X_{e(\alpha)}$ is an
$A$-map, such that $X_\gamma: =X_{\alpha_l }\cdots X_{\alpha_1} = 0$
for each $\gamma =\alpha_l\cdots\alpha_1\in \rho$, where $\rho$ is a minimal
set of generators of $I$. We call $X_i$
{\it the $i$-th branch} of $X$. A morphism $f$ from $X$ to $Y$ is
$(f_i, i\in Q_0)$ with each $f_i: X_i \rightarrow Y_i$ an $A$-map, such that for each arrow $\alpha: j\rightarrow i$ there
holds $$Y_\alpha f_j = f_iX_\alpha. \eqno(1.1)$$ Denote by ${\rm
rep}(Q, I, A)$ the category of representations of
$(Q,  I)$ over $A$. As in the case of $A = k$ ([ARS, p.57], [R,
p.44]),  we have $\m$-mod $\cong {\rm rep}(Q, I, A)$. {\bf
Throughout, we will identify left $\m$-modules with representations of
$(Q, I)$ over $A$}. Also, mod$\m \cong {\rm rep}(Q^{op}, I^{op}, A^{op})$, where $(Q^{op}, I^{op})$ is the opposite bound quiver of $(Q, I)$, and $A^{op}$ is the opposite algebra of $A$.

\vskip5pt

This realization of $\m$-modules leads up to the separated monomorphism category ${\rm smon}(Q, I, A)$.
Since G.Birkhoff [Bir], there is a long history of
studying  ${\rm smon}(Q, 0, A)$, where $Q$ is chain: it is called {\it the submodule category} for $n =2$ ([RS1-RS3]), and {\it the filtered chain category} ([S1-S3, SW]).
C. M. Ringel and M. Schmidmeier [RS2] establish the Auslander-Reiten
theory of the submodule category (see also [XZZ]); and D. Simson  ([S2], [S3]) studies the representation type. By  D. Kussin, H. Lenzing and H. Meltzer ([KLM1, KLM2]; see also [C1]), it is related to the singularity theory.

\vskip5pt

\subsection{} In this paper we introduce separated monic representations of $(Q, I)$
over an arbitrary additive subcategory $\x$ of $A$-mod (see Defintion \ref{maindef1}), and denote by ${\rm smon}(Q, I, \x)$ the category of
separated monic representations of $(Q, I)$ over $\x$. This generalizes the notion ${\rm smon}(Q, I, A)$ in [LZ].
This wider setting is useful in particular for finding Frobenius subcategories of $\m$-mod. Another motivation is that $\mathcal
{GP}(\m)$ can be described as ${\rm smon}(Q, I, \mathcal
{GP}(A))$ ([LZ, Thm. 4.1]) and also $\mathcal {P}(\m) = {\rm smon}(Q, I, \mathcal
{P}(A))$ (see [LZ, Lemma 4.9]). In particular, $\mathcal {P}(kQ/I) = {\rm smon}(Q, I, k)$.

\vskip5pt For an $A$-module $T$, let ${\rm add}(T)$ be the
subcategory of $A$-mod of direct summands of finite direct sums of
copies of $T$, and $^\perp T$ the subcategory given by \ $^\perp T:
= \{M\in A\mbox{-}{\rm mod} \ | \ {\rm Ext}^m_A(X, T) = 0, \ \forall
\ m\ge 1\}.$ Denote by $\D: = \Hom_k(-, k)$.

\vskip5pt

Although ${\rm smon}(Q, I, \x)$ is defined
combinatorially, it admits a homological description as
$${\rm smon}(Q, I, A) = \{X\in \m\mbox{-}{\rm mod} \ | \ {\rm Tor}^\m_1(A_A\otimes \D S, X) = 0\}
= \ ^\bot(\D(A_A)\otimes kQ/I)$$ and
$${\rm smon}(Q, I, \x) = \{X\in {\rm smon}(Q, I, A) \ | \ (A\otimes \D S(i))\otimes_\m X\in\x, \ \forall \ i\in Q_0\};$$
and a striking property is that it relates to the tilting theory
via a reciprocity ${\rm smon}(Q, I, \ ^\bot T)= \ ^\bot(T\otimes kQ/I)$ for any cotilting module $_AT$. See Theorem \ref{homodescriptionofsmon}.
As an application, ${\rm smon}(Q, I, \x)$
has Auslander-Reiten sequences, provided that $\x$ is resolving and
contravariantly finite with $\widehat{\x}=A$-mod (Theorem \ref{ar}).
In particular, ${\rm smon}(Q, I, A)$ always has Auslander-Reiten sequences.

\vskip5pt

\subsection{} The dual version ${\rm sepi}(Q, I, \x)$ can be interpreted as $\D{\rm smon}(Q^{op}, I^{op}, \D\x)$ (Proposition \ref{smonsepi}).
Ringel-Schmidmeier-Simson equivalence $F\colon{\rm smon}(Q, I,
\x)\cong{\rm sepi}(Q, I, \x)$ is introduced (Definition \ref{rss}). Such an equivalence was first observed in [RS2] and [S1] for a chain
$Q$ with $I = 0$. It implies a strong symmetry, in
particular, that the separated monic representations are as many as
the separated epic representations. We prove the existence of an {\rm RSS} equivalence $\D\Hom_\m(-, \D(A_A)\otimes kQ/I):
{\rm smon}(Q, I, \x)\cong {\rm sepi}(Q, I, \x)$ for extension-closed subcategories $\x$ (Theorem \ref{RSS0}).
As a consequence, if $\x \supseteq \mathcal P(A)$, then the Nakayama functor $\mathcal N_\m$ gives an {\rm RSS} equivalence if and only if $A$ is a Frobenius algebra.

\vskip5pt

For a chain $Q$
with arbitrary monomial ideal $I$, an  {\rm RSS} equivalence
is found (Theorem \ref{RSS1}), by using ``jumping of ${\rm Cok}$", which
is a generalization of Ringel-Schmidmeier and Simson's functor ${\rm Cok}$.
Quite surprising, this {\rm RSS} equivalence is for arbitrary additive subcategories $\x$, not necessarily extension-closed, and it is combinatorial and operable.

\subsection{} Let $Q_0$ (resp. $Q_1$, $\mathcal P$)
denote the set of vertices (resp. arrows, paths)
of $Q$, $s(p)$ (resp. $e(p)$, $l(p)$) the starting vertex (resp.
the ending vertex, the length) of $p\in \mathcal P$.
We write the conjunction of paths from the right to the left. A
vertex $i$ is viewed as a path of length $0$, denoted by $e_i$.
Label $Q_0$ as $1, \cdots, n$, such that $j
> i$ if $\alpha:
j\rightarrow i$ is in $Q_1$. Thus $n$ is a source and $1$ is a sink.
Let  $P(i)$ (resp. $I(i)$, $S(i)$) be the indecomposable projective
(resp. injective, simple) $(kQ/I)$-module at $i\in Q_0$. For $i\in Q_0$, put
$\ha(\to i): = \{\alpha\in Q_1 \ | \ e(\alpha) =i\}$,  and for $\alpha\in Q_1$, put
$$K_\alpha: = \{p\in \mathcal P \ | \ l(p) \ge 1, \ e(p) = s(\alpha), \ p\notin I, \ \alpha p \in
I\}. \eqno (1.2)$$
For  $X=(X_i, X_{\alpha}, i\in Q_0, \alpha\in Q_1)\in {\rm
rep}(Q, I, A)$ and $i\in Q_0$, put
$$\delta_i(X):
=(X_\alpha)_{\alpha\in \ha(\to i)}:  \bigoplus\limits_{\alpha\in
\ha(\to i) } X_{s(\alpha)}\longrightarrow X_i\eqno(1.3)$$ and for
$\alpha\in Q_1$, put
$$X^\alpha: = (X_p)_{p\in K_\alpha}: \bigoplus\limits_{p\in
K_\alpha} X_{s(p)} \longrightarrow X_{s(\alpha)}.\eqno(1.4)$$

\vskip5pt \noindent {\bf Conventions:} If $i$ is a source then $\delta_i(X):=0$; and if
$K_\alpha = \emptyset$ then $X^\alpha: = 0$.

\vskip5pt

Define $\cok_i: \m$-mod $\longrightarrow A$-mod to be the functor given
by $$\cok_i(X):= {\rm Coker} \ \delta_i(X)= X_i/\sum\limits_{\begin
{smallmatrix} \alpha\in \ha(\to i)
\end{smallmatrix}} \Ima X_{\alpha}\eqno (1.5)$$ (if $i$ is a source then
$\cok_i(X):=X_i$).

\vskip5pt

For example, if $T\in A$-mod and $M\in (kQ/I)$-mod with $M=(M_i,
M_{\alpha}, i\in Q_0, \alpha\in Q_1)\in {\rm rep}(Q, I, k)$, then
$T\otimes M\in \m$-mod with $$T\otimes M=(T\otimes M_i, {\rm
Id}_T\otimes M_{\alpha}, i\in Q_0, \alpha\in Q_1)\in {\rm rep}(Q, I,
A),  \ \ \mbox{and} \ \ \cok_i(T\otimes M) = T\otimes \cok_i(M),$$
where $\cok_i(M)\in k$-mod. In
particular, for $v\in Q_0$, the indecomposable projective left
$(kQ/I)$-module $P(v)\colon = (kQ/I)e_v$, viewed as a representation
in ${\rm rep}(Q, I, k)$, has the form
$$P(v) = (e_i(kQ/I)e_v, \ P(v)_\alpha, \ i\in Q_0, \ \alpha\in
Q_1),\eqno (1.6)$$ where $P(v)_\alpha: \ e_{s(\alpha)}(kQ/I)e_v
\longrightarrow e_{e(\alpha)}(kQ/I)e_v$ is the $k$-map sending path
$w$ to $\alpha w$, and
$$\cok_i(T\otimes P(v)) = \begin{cases} T, & \mbox{if} \ i=v,
\\ 0, & \mbox{if}\ i\neq v.\end{cases}  \eqno (1.7)$$

\begin{lem}\label{adjointpair} \ For each $i\in
Q_0$,  we have

\vskip5pt

$(1)$ \  $(\cok_i, -\otimes  S(i))$ is an adjoint pair.

\vskip5pt

$(2)$ \ The pairwise non-isomorphic indecomposable projective $\m$-modules are exactly $P\otimes P(i)$, where $P$ runs over pairwise non-isomorphic indecomposable projective $A$-module, and $i$ runs over $Q_0$. Thus, the branches of
projective $\m$-modules are projective $A$-modules.
\end{lem}

\section {\bf Separated monomorphism categories}

\subsection{}  A main notion is:

\begin{defn}\label{maindef1}  \ Let $\x$ be an additive subcategory of $A$-mod. A representation $X = (X_i,
X_{\alpha}, i\in Q_0, \alpha\in Q_1)\in {\rm rep}(Q, I, A)$ is separated
monic over $\x$, or a separated monic
left $\Lambda$-module over $\x$, provided that $X$ satisfies the
following conditions$:$

\vskip5pt

${\rm (m1)}$  \ For each $i\in Q_0$,  $\sum\limits_{\begin
{smallmatrix} \alpha\in \ha(\to i)
\end{smallmatrix}}\Ima
X_\alpha$ is the direct sum $\bigoplus\limits_{\begin {smallmatrix}
\alpha\in \ha(\to i) \end{smallmatrix}}\Ima X_\alpha;$

\vskip5pt

${\rm (m2)}$     \ For each $\alpha\in Q_1$,  \  $\Ker
X_\alpha =\sum\limits_{\begin {smallmatrix} p\in
K_\alpha\end{smallmatrix}} \Ima X_{p};$

\vskip5pt

${\rm (m3)}$  \ For each $i\in Q_0$, $\cok_i(X)\in\x$.
\end{defn}

We keep conventions: if $i$ is a source, then $\sum\limits_{\begin
{smallmatrix} \alpha\in \ha(\to i)\end{smallmatrix}} \Ima X_\alpha = 0 = \bigoplus\limits_{\begin {smallmatrix}
\alpha\in \ha(\to i) \end{smallmatrix}}\Ima X_\alpha$; if $K_\alpha = \emptyset$ (e.g. if
$I = 0$ or $s(\alpha)$ is a source),  then $\sum\limits_{\begin
{smallmatrix} p\in K_\alpha\end{smallmatrix}} \Ima X_p = 0$.

\vskip5pt

The conditions ${\rm (m1)}-{\rm (m3)}$ are independent.  Denote by
${\rm smon}(Q, I, \x)$ the subcategory of ${\rm rep}(Q, I, A)$ consisting of
separated monic representations over $\x$. We write
${\rm smon}(Q, I, A)$ for ${\rm smon}(Q, I, \ A\mbox{-}{\rm mod})$.
If $\x$ is closed under direct
summands, then so is ${\rm smon}(Q, I, \x)$, and hence ${\rm
smon}(Q, I, \x)$ is a Krull-Schmidt category ([R, p.52]); in this case, indecomposable objects of ${\rm smon}(Q, I, \x)$ are also indecomposable as $\m$-modules.

\begin{rem}\label{morita} Note that ${\rm smon}(Q, I, \x)$ is invariant under {\rm Morita} equivalences, in the sense that if $A\mbox{-}{\rm mod} \stackrel{G}\cong A'\mbox{-}{\rm mod}$, then
${\rm smon}(Q, I, \x)\cong {\rm smon}(Q, I, G\x),$ sending $(X_i,
X_{\alpha}, i\in Q_0, \alpha\in Q_1)\in {\rm smon}(Q, I, \x)$ to $(G(X_i),
G(X_{\alpha}), i\in Q_0, \alpha\in Q_1)$.
\end{rem}

\begin{exm} \label{proj} Let $T\in \x$, and  $M=(M_i,
M_{\alpha}, i\in Q_0, \alpha\in Q_1)\in {\rm rep}(Q, I, k)$. By definition one can see that \
$T\otimes M = (T\otimes M_i, {\rm Id}\otimes M_{\alpha},
i\in Q_0, \alpha\in Q_1)\in {\rm smon}(Q, I, \x)$ if and only if $M\in {\rm
smon}(Q, I, k)= \mathcal P(kQ/I)$. Thus, for each $i\in Q_0$
we have a functor $-\otimes P(i): \x\longrightarrow {\rm smon}(Q, I, \x),$  and
${\rm smon}(Q, I, A)\cap (A\mbox{-}{\rm mod}\otimes
(kQ/I)\mbox{-}{\rm mod}) = A\mbox{-}{\rm mod}\otimes \mathcal
P(kQ/I).$  If $\x$ is closed under direct summands, then
${\rm smon}(Q, I, \x)\cap (A\mbox{-}{\rm mod}\otimes
(kQ/I)\mbox{-}{\rm mod}) = \x\otimes \mathcal
P(kQ/I).$

\end{exm}

We need a reformulation of ${\rm (m1)} - {\rm (m3)}$.

\begin{lem} \label{anothermono} Let $X=(X_i, X_{\alpha}, i\in Q_0, \alpha\in Q_1)\in {\rm
rep}(Q, I, A)$. Then

\vskip5pt

$(1)$ \ $X\in {\rm smon}(Q, I, A)$ if and only if for
each $i\in Q_0$, the sequence of $A$-maps $$\bigoplus\limits_{\alpha\in
\ha(\to i)}\bigoplus\limits_{\begin{smallmatrix}q\in K_\alpha
\end{smallmatrix}} X_{s(q)} \stackrel {\bigoplus\limits_{\alpha\in
\ha(\to i)} X^\alpha}\longrightarrow \bigoplus\limits_{\begin {smallmatrix}
\alpha\in \ha(\to i)
\end{smallmatrix}}X_{s(\alpha)}\stackrel {\delta_i(X)}\longrightarrow X_i \eqno(2.1)$$
is exact, where $\delta_i(X)$ and $X^\alpha$  are defined in
$(1.3)$ and $(1.4)$, respectively.

\vskip5pt

$(2)$ \ For each $i\in Q_0$,  $\cok_i(X) = (A\otimes \D S(i)) \otimes_{\m} X$.
\end{lem}

\noindent{\bf Proof.}  \ $(1)$ If $X\in {\rm smon}(Q, I, A)$, then
\begin{align*}{\rm Ker} \ \delta_i(X) & = \{ (x_\alpha)_{\alpha\in \ha(\to i)} \in \bigoplus\limits_{\alpha\in
\ha(\to i)} X_{s(\alpha)} \ | \ \sum\limits_{\alpha\in \ha(\to i)}
X_\alpha(x_\alpha) = 0 \  \}\\ & \stackrel{{\rm (m1)}} =
\bigoplus\limits_{\alpha\in \ha(\to i)}{\rm Ker}
X_\alpha\stackrel{{\rm (m2)}} =\bigoplus\limits_{\alpha\in \ha(\to
i)} \sum\limits_{\begin {smallmatrix} p\in
K_\alpha\end{smallmatrix}} \Ima X_{p}\\ &
=\bigoplus\limits_{\alpha\in \ha(\to i)} {\rm Im} X^\alpha = {\rm
Im}(\bigoplus\limits_{\alpha\in \ha(\to i)} X^\alpha). \end{align*}

Conversely, if $(2.1)$ is exact for each $i\in Q_0$, then by ${\rm
Im}(\bigoplus\limits_{\alpha\in \ha(\to i)} X^\alpha) = {\rm Ker} \
\delta_i(X)$ and the fact
\begin{align*}{\rm
Im}(\bigoplus\limits_{\alpha\in \ha(\to i)} X^\alpha)&
=\bigoplus\limits_{\alpha\in \ha(\to i)} {\rm Im} X^\alpha =
\bigoplus\limits_{\alpha\in \ha(\to i)} \sum\limits_{\begin
{smallmatrix} p\in K_\alpha\end{smallmatrix}} \Ima X_{p} \subseteq
\bigoplus\limits_{\alpha\in \ha(\to i)}{\rm Ker} X_\alpha
\\ &\subseteq \{
(x_\alpha)_{\alpha\in \ha(\to i)} \in \bigoplus\limits_{\alpha\in
\ha(\to i)} X_{s(\alpha)} \ | \ \sum\limits_{\alpha\in \ha(\to i)}
X_\alpha(x_\alpha) = 0 \  \} = {\rm Ker} \ \delta_i(X)
\end{align*}
we see that ${\rm Ker} X_\alpha = \sum\limits_{\begin {smallmatrix}
p\in K_\alpha\end{smallmatrix}} \Ima X_{p}$ for each $\alpha\in
\ha(\to i)$, and that $\sum\limits_{\begin {smallmatrix}
\alpha\in \ha(\to i)
\end{smallmatrix}}\Ima
X_\alpha$ is the direct sum $\bigoplus\limits_{\begin {smallmatrix}
\alpha\in \ha(\to i) \end{smallmatrix}}\Ima X_\alpha$ for each $i\in
Q_0$, i.e., $X\in {\rm smon}(Q, I, A)$.

\vskip5pt

$(2)$ \ We have a projective presentation of right
$(kQ/I)$-module $\D S(i):$
$$\bigoplus\limits_{\begin
{smallmatrix} \alpha\in \ha(\to i) \end{smallmatrix}}
e_{s(\alpha)}(kQ/I)\stackrel{(\alpha.)}\longrightarrow e_i
(kQ/I)\longrightarrow \D S(i)\longrightarrow 0,$$ where $\alpha.$ is the right $(kQ/I)$-map given by the left multiplication by
$\alpha$. Applying $A\otimes -$ we get a projective
presentation of right $\m$-module $A\otimes \D S(i)$:
$$\bigoplus\limits_{\begin {smallmatrix} \alpha \in \ha(\to i)
\end{smallmatrix}} (1\otimes e_{s(\alpha)})\m \lraf{((1\otimes
\alpha).)} (1\otimes e_i)\m \longrightarrow  A\otimes \D S(i))\longrightarrow
0$$ where $1$ is the identity of $A$, $(1\otimes\alpha).$ is the right $\m$-map given by the left
multiplication by $(1\otimes\alpha)$. For
$X\in \m$-mod, applying $-\otimes_\m X$ we get the
exact sequence
$$\bigoplus\limits_{\begin {smallmatrix} \alpha\in \ha(\to i)
\end{smallmatrix}} X_{s(\alpha)} \lraf{(X_{\alpha})_{\alpha\in
\ha(\to i)}} X_i  \longrightarrow (A\otimes \D S(i)) \otimes_{\m} X
\longrightarrow 0.$$ Since $\delta_i(X) = (X_{\alpha})_{\alpha\in
\ha(\to i)}, $  we get $\cok_i(X) = {\rm Coker} \ \delta_i(X) = (A\otimes \D S(i)) \otimes_{\m} X$.  $\s$

\subsection{} A subcategory $\x$ of $A$-mod is {\it
resolving} if $\x$ is
closed under extensions, kernels of epimorphisms, and direct
summands,  and $\x\supseteq \mathcal P(A)$ ([AR]).

\begin{lem} \label{cokexact} $(1)$ \ For $i\in Q_0$,
the restriction of \ $\cok_i:\ \m$-mod $\longrightarrow A$-mod to ${\rm
smon}(Q, I, A)$ is an exact functor.

\vskip5pt

$(2)$ \  ${\rm smon}(Q, I, \x)$ is closed under extensions $($resp.
kernels of epimorphisms$;$ direct summands$)$ if and only if $\x$ is
closed under extensions $($resp. kernels of epimorphisms$;$ direct
summands$)$.

\vskip5pt

Thus, ${\rm smon}(Q, I, \x)$ is a resolving subcategory of $\m$-mod
if and only if $\x$ is a resolving subcategory of $A$-mod.
\end{lem} \noindent{\bf Proof.} $(1)$ \ Let $0\rightarrow (X_i,
X_{\alpha})\lraf{f} (Y_i, Y_{\alpha}) \lraf{g} (Z_i,
Z_{\alpha})\rightarrow 0$ be an exact sequence in ${\rm smon}(Q, I,
A)$. By $(1.1)$ and $(1.4)$ we have a commutative diagram
\[\xymatrix @R=0.85cm@C=0.65cm{\bigoplus \limits_{\begin{smallmatrix} q\in K_\alpha \end{smallmatrix}}Y_{_{s(q)}}\ar@{>>}[r]^{\oplus g_{{s(q)}}}
\ar@{>>}[d]^-{Y^{\alpha}}& \bigoplus \limits_{\begin{smallmatrix} q\in K_\alpha \end{smallmatrix}}Z_{s(q)}\ar@{>>}[d]^-{Z^{\alpha}} \\
\Ima Y^{ \alpha}\ar[r] & \Ima Z^{\alpha}}\] So $\Ima
Y^\alpha\rightarrow \Ima Z^\alpha$ is epic. By $(2.1)$,
$\Ker {\delta_i}(Y)=\bigoplus
\limits_{\begin{smallmatrix}\alpha\in \ha(\to i) \end{smallmatrix}}
\Ima Y^\alpha$ and $ \Ker  {\delta_i}(Z)=\bigoplus
\limits_{\begin{smallmatrix} \alpha\in \ha(\to i)
\end{smallmatrix}}\Ima Z^\alpha$. Thus
$\Ker{\delta_i}(Y)\rightarrow \Ker{\delta_i}(Z)$ is epic, and the assertion follows from the
commutative diagram

$$\small{\CD 0 @>>> \bigoplus\limits_{\begin {smallmatrix} \alpha\in \ha(\to i)\end{smallmatrix}}  X_{s(\alpha)}@>\oplus f_{s(\alpha)}>> \bigoplus\limits_{\begin {smallmatrix} \alpha\in \ha(\to i)\end{smallmatrix}}
   Y_{s(\alpha)} @>\oplus g_{s(\alpha)}>>\bigoplus\limits_{\begin {smallmatrix} \alpha\in \ha(\to i)\end{smallmatrix}}  Z_{s(\alpha)} @>>> 0   \\
  && @V {\delta_i}(X) VV @V {\delta_i}(Y) VV @V {\delta_i}(Z) VV\\
  0 @>>> X_i @>f_i>> Y_i @>g_i>> Z_i @>>>0
\endCD}$$

\vskip5pt \noindent and the Snake Lemma.

\vskip5pt

$(2)$ \ Using $(1)$ and the fact that ${\rm smon}(Q, I, A)$ is a
resolving subcategory of $\m$-mod (see [LZ, Thm. 3.1]), the
assertion can be easily proved: the argument is as follows.

\vskip5pt For example, assume that $\x$ is extension-closed. Let
$0\rightarrow (X_i, X_{\alpha})\rightarrow (Y_i,
Y_{\alpha})\rightarrow (Z_i, Z_{\alpha})\rightarrow 0$ be an
exact sequence in ${\rm rep}(Q, I, A)$ with $X= (X_i,
X_{\alpha})\in{\rm smon}(Q, I, \x)$ and $Z=(Z_i, Z_{\alpha})\in{\rm
smon}(Q, I, \x).$ Since ${\rm smon}(Q, I, \x)\subseteq {\rm smon}(Q,
I, A)$, and ${\rm smon}(Q, I, A)$ is closed under extensions, it
follows that $Y=(Y_i, Y_{\alpha})\in {\rm smon}(Q, I, A)$, and
hence by $(1)$ we get an exact sequence $0\rightarrow
\cok_i(X)\rightarrow \cok_i(Y)\rightarrow
\cok_i(Z)\rightarrow 0$ with $\cok_i(X)\in\x$ and
$\cok_i(Z)\in\x$ for $i\in Q_0$. By the assumption
$\cok_i(Y)\in\x$, and hence $Y\in {\rm smon}(Q, I, \x)$.

\vskip5pt

Also, for example, assume that ${\rm smon}(Q, I, \x)$ is
extension-closed. Let $0\rightarrow X\rightarrow Y\rightarrow
Z\rightarrow 0$ be an exact sequence in $A$-mod with $X\in\x$ and
$Z\in \x$. Define $(X_i, X_{\alpha})\in {\rm rep}(Q, I, A)$ by
$X_1:=X$ and $X_i:=0$ for $i\neq 1$ (note that $1$ is a
sink of $Q$). Similarly for $(Y_i, Y_{\alpha})$ and $(Z_i,
Z_{\alpha})$. Since $1$ is a sink of $Q$, it follows that
$0\rightarrow (X_i, X_{\alpha})\rightarrow (Y_i,
Y_{\alpha})\rightarrow (Z_i, Z_{\alpha})\rightarrow 0$ is exact in
${\rm smon}(Q, I, A)$ with $(X_i, X_{\alpha})\in{\rm smon}(Q, I,
\x)$ and $(Z_i, Z_{\alpha})\in{\rm smon}(Q, I, \x)$, and hence $(Y_i, Y_{\alpha})\in {\rm smon}(Q, I, \x)$. So $Y =
\cok_1((Y_i, Y_{\alpha}))\in \x$.

\vskip5pt

By the same argument we see that ${\rm smon}(Q, I, \x)$ is closed under kernels of epimorphisms (resp. direct summands) if and only
if $\x$ is closed under kernels of
epimorphisms (resp. direct summands).

\vskip5pt

By the same argument as above,  $\mathcal
X\supseteq \mathcal P(A)$ if and only if ${\rm smon}(Q, I,
\x)\supseteq \mathcal P(\m)$. Thus ${\rm smon}(Q, I, \x)$ is a
resolving subcategory of $\m$-mod if and only if so is $\x$ of
$A$-mod. $\s$

\subsection{A homological interpretation and a reciprocity}
The separated monomorphism category is defined combinatorially, but it admits a homological interpretation,
which derives a reciprocity between the perpendicular operator and the monomorphism operator. It was observed for a chain $Q$ with $I = 0$ in [Z].
Put $S: =\bigoplus\limits_{\begin {smallmatrix} i\in
Q_0\end{smallmatrix}}S(i)$.

\begin{thm}\label{homodescriptionofsmon}  $(1)$ \ We have
\begin{align*}{\rm smon}(Q, I, A) = \{X\in \m\mbox{-}{\rm mod} \ | \ {\rm Tor}^\m_1(A_A\otimes \D S, X) = 0\}
= \ ^\bot(\D(A_A)\otimes kQ/I)\end{align*}
and
\begin{align*}{\rm smon}(Q, I, \x)& = \{X\in {\rm smon}(Q, I, A) \ | \ (A\otimes \D S(i))\otimes_\m X\in\x, \ \forall \ i\in Q_0\} \\ &
= \{X\in {\rm smon}(Q, I, A) \ | \ (A\otimes M)\otimes_\m X\in\x \ \mbox{for each right} \ (kQ/I)\mbox{-module} \ M\}.\end{align*}
where for the second equality we assume that $\x$ is extension-closed.

\vskip5pt

$(2)$ \ Let $T$ be an $A$-module. Then \ ${\rm smon}(Q, I, \ ^\bot T)= {\rm smon}(Q, I,
A)\cap \ ^\bot(T\otimes kQ/I)$.

\vskip5pt

Moreover, if there is an exact sequence $0\rightarrow T_m\rightarrow
\cdots\rightarrow T_0\rightarrow \D(A_A)\rightarrow 0$ with each
$T_j\in {\rm add} (T)$, then ${\rm smon}(Q, I, \ ^\bot T)= \ ^\bot
(T\otimes kQ/I )$.
\end{thm}

\vskip5pt

To prove Theorem \ref{homodescriptionofsmon} we need the following fact.

\begin{lem}\label{Tor} If $X\in
\m$-mod and ${\rm Tor}^{\m}_1(A_A\otimes \D S, X)=0$, then
${\rm Tor}^{\m}_m(A_A\otimes \D S, X)=0, \ \forall \ m\geq 1.$
\end{lem}
\noindent{\bf Proof.} \ For $i\in Q_0$, define
$l_i:=\mbox{max}\{ l(p)\ | \ p\in \hp \ \mbox{with}\ e(p)=i\}$.
Since $Q$ is acyclic, we can use induction on $l_i$ to prove the assertion: for $i\in Q_0$ there holds ${\rm
Tor}^{\m}_m(A_A\otimes \D S(i), X)=0, \ \forall \ m\geq 1$.

If $l_i=0$, then $i$ is a source, $S(i)$ is an injective
$(kQ/I)$-module, and $A_A\otimes \D S(i)$ is a projective
right $\m$-module. So ${\rm Tor}^{\m}_m(A_A\otimes
\D S(i), X)=0, \ \forall \ m\geq 1.$  Assume that the assertion holds for $i\in Q_0$ with $l_i\leq l-1.$
For $i\in Q_0$ with $l_i=l$, let $M$ be the right $(kQ/I)$-module such that
 $0\rightarrow M\rightarrow e_i (kQ/I) \stackrel \pi \rightarrow \D S(i)\rightarrow 0$ is an exact sequence of right $(kQ/I)$-modules, where $\pi$
 is a projective cover.  Since $Q$ is acyclic, any composition factor of $M$ is of the form $\D S(j)$ with
$j\in Q_0$ satisfying $l_j\leq l-1$: in fact, if $\D S(j)$ is a
composition factor of $M$, then $\Hom_{kQ/I}( e_j(kQ/I),
M)\ne 0$ and $\Hom_{kQ/I}( e_j(kQ/I), e_i (kQ/I))\ne
0$, which means that there is a path from $j$ to $i$; also since the multiplicity of $\D S(i)$ in  $e_i (kQ/I)$ is $1$, it follows that $j\ne i$, and hence
$l_j< l$. By induction,
${\rm Tor}^{\m}_m(A_A\otimes \D S(j), X)=0, \ \forall \
m\geq 1$, and hence ${\rm Tor}^{\m}_m(A_A\otimes M, X)=0, \ \forall
\ m\geq 1$. Apply $-\otimes_\m X$ to the exact sequence
$0\rightarrow A_A\otimes M\rightarrow A_A\otimes e_i(kQ/I)
\rightarrow A_A \otimes \D S(i)\rightarrow 0$. Since $A_A\otimes
e_i(kQ/I)$ is a projective right $\m$-module and ${\rm Tor}^{\m}_1(A_A\otimes \D S(i), X)=0$, we see
that ${\rm Tor}^{\m}_m(A_A\otimes \D S(i), X)=0, \ \forall \ m\ge 1$.  $\s$

\vskip5pt

\noindent{\bf Proof of Theorem \ref{homodescriptionofsmon}.} \ $(1)$ \ For each $i\in Q_0$, we have a projective resolution of right
$(kQ/I)$-module $\D S(i):$
$$\cdots\rightarrow \bigoplus\limits_{\begin {smallmatrix} \alpha\in \ha(\to i) \end{smallmatrix}} \bigoplus\limits_{\begin {smallmatrix} q\in K_\alpha \end{smallmatrix}}
e_{s(q)}(kQ/I)\stackrel{\bigoplus\limits_{\begin {smallmatrix}
\alpha\in \ha(\to i)
\end{smallmatrix}}(q.)}\longrightarrow\bigoplus\limits_{\begin
{smallmatrix} \alpha\in \ha(\to i) \end{smallmatrix}}
e_{s(\alpha)}(kQ/I)\stackrel{(\alpha.)}\rightarrow  e_i
(kQ/I)\rightarrow \D S(i)\rightarrow  0,$$ where $\alpha.$ (resp.
$q.$) means the right $(kQ/I)$-map given by left multiplication by
$\alpha$ (resp. $q$). Applying $A\otimes -$ we get a projective
resolution of right $\m$-module $A\otimes \D S(i)$:
$$\cdots\rightarrow   \bigoplus\limits_{\begin {smallmatrix} \alpha\in \ha(\to i) \end{smallmatrix}}
\bigoplus\limits_{\begin {smallmatrix} q\in K_\alpha
\end{smallmatrix}} (1\otimes e_{s(q)})\m
\stackrel{\bigoplus\limits_\alpha (1\otimes q).} \longrightarrow
\bigoplus\limits_{\begin {smallmatrix} \alpha \in \ha(\to i)
\end{smallmatrix}} (1\otimes e_{s(\alpha)})\m \lraf{((1\otimes
\alpha).)} (1\otimes e_i)\m \rightarrow  A\otimes \D S(i)\rightarrow
0$$ where $1$ is the identity of $A$, $(1\otimes\alpha).$ (resp.
$(1\otimes q).$) is the right $\m$-map given by left
multiplication by $(1\otimes\alpha)$ (resp. $(1\otimes q)$). For
$X\in \m$-mod, applying $-\otimes_\m X$ we get the
sequence (cf. $(2.1)$)
$$\cdots\rightarrow  \bigoplus\limits_{\begin {smallmatrix} \alpha\in \ha(\to i) \end{smallmatrix}}  \bigoplus\limits_{\begin {smallmatrix}q\in K_\alpha
\end{smallmatrix}} X_{s(q)} \lraf{\bigoplus\limits_\alpha(X_q)_{q\in K_\alpha}}
\bigoplus\limits_{\begin {smallmatrix} \alpha\in \ha(\to i)
\end{smallmatrix}} X_{s(\alpha)} \lraf{(X_{\alpha})_{\alpha\in
\ha(\to i)}} X_i \rightarrow (A\otimes \D S(i)) \otimes_{\m} X
\rightarrow  0$$ where  $(X_i, X_\alpha)$  is the representation of
$(Q, I)$ over $A$ corresponding to $_\m X$ (thus $X_i=(1\otimes
e_i)X$).

\vskip5pt

By Lemma \ref{anothermono}$(1)$,   $X\in {\rm smon}(Q, I, A)$ if and only if
the sequence
$$\bigoplus\limits_{\begin {smallmatrix} \alpha\in \ha(\to i) \end{smallmatrix}}  \bigoplus\limits_{\begin {smallmatrix}q\in K_\alpha \end{smallmatrix}}
X_{s(q)} \lraf{\bigoplus\limits_\alpha (X_{q})_{q\in K_\alpha}}
\bigoplus\limits_{\begin {smallmatrix} \alpha\in \ha(\to i)
\end{smallmatrix}} X_{s(\alpha)} \lraf{(X_{\alpha})_{\alpha\in
\ha(\to i)}} X_i \rightarrow (A\otimes \D S(i)) \otimes_{\m} X
\rightarrow  0$$ is an exact sequence of left $A$-modules for all
$i\in Q_0$, and if and only if ${\rm Tor}^{\m}_1(A\otimes \D S(i),
X)=0$, $\forall \ i\in Q_0$, or equivalently, if and only if ${\rm
Tor}^{\m}_1(A\otimes \D S, X)=0$. So by Lemma \ref{Tor}, $X\in {\rm
smon}(Q, I, A)$ if and only if ${\rm Tor}^{\m}_m(A\otimes \D S,
X)=0$, \ $\forall \ m\geq 1$.

\vskip5pt

Since
${\rm proj.dim}\ _{kQ/I}S < \infty,$ it is clear that ${\rm Tor}^{\m}_m(A\otimes \D S,
X)=0$ for all \ $m\geq 1$ if and only if
${\rm Tor}^{\m}_m(A\otimes \D(kQ/I),
X)=0$ for $m\geq 1$.

\vskip5pt

By the well-known isomorphism
$$\D{\rm Tor}^{\m}_m( A\otimes \D(kQ/I), X)\cong {\rm Ext}_{\m}^m(X,
\D(A)\otimes kQ/I)$$ (To see this isomorphism, take a
projective resolution of $_\m X$ and consider the Tensor-Hom adjoint
pair $((A\otimes \D(kQ/I))\otimes_\m -, \Hom_k(A\otimes \D(kQ/I), -)$ between  $\m\mbox{-}{\rm
mod}$ and $k\mbox{-}{\rm mod}$.) we
see that $X\in {\rm smon}(Q, I, A)$ if and only if ${\rm Ext}_{\m}^m(X,
\D(A)\otimes kQ/I)=0$ for $m\geq 1$, i.e., if and only if
$X\in \ ^\bot(\D(A_A)\otimes kQ/I)$.

\vskip5pt

By Lemma \ref{anothermono}$(2)$ we have ${\rm smon}(Q, I, \x) = \{X\in {\rm smon}(Q, I, A) \ | \ (A_A\otimes \D S(i))\otimes_\m X\in\x, \ \forall \ i\in Q_0\}$.
Since any right $(kQ/I)$-module $M$ has composition series and $\x$ is extension-closed,
it follows that $X\in {\rm smon}(Q, I, \x)$ if and only if $X\in {\rm smon}(Q, I, A)$ and $(A\otimes M)\otimes_\m X\in\x \ \mbox{for each right} \ (kQ/I)\mbox{-module} \ M$.

\vskip5pt

$(2)$ \  Let $X\in {\rm smon}(Q, I, A)$, and $i\in Q_0$. We
claim that  $\cok_i(X)\in \ ^\bot T$ if and only if $X\in \
^\bot (T\otimes S(i)).$

In fact, we take a $\m$-projective resolution $\cdots \rightarrow
P^1 \stackrel{d_1}\rightarrow P^0 \stackrel{d_0}\rightarrow
X\rightarrow 0$. By Lemma \ref{cokexact}$(2)$, ${\rm
Ker}d^i\in {\rm smon}(Q, I, A)$ for each $i\ge 0$. By Lemma
\ref{cokexact}$(1)$ we get an exact sequence $\cdots \rightarrow
\cok_i(P^1)\rightarrow \cok_i(P^0)\rightarrow
\cok_i(X)\rightarrow 0,$  and by $(1.7)$,
$\cok_i(A\otimes kQ/I)=A,$ so it is a projective resolution of
$\cok_i(X)$. By Lemma \ref{adjointpair}$(1)$ we get a commutative
diagram
$$\xymatrix@R=0.65cm@C=0.4cm{0\ar[r] & \Hom_A(\cok_i(X), T)\ar[r]\ar[d]^-{\wr} &\Hom_A(\cok_i(P^0), T)\ar[r]\ar[d]^-{\wr} & \Hom_A(\cok_i(P^1), T)
\ar[r]\ar[d]^-{\wr} & \cdots\\
 0 \ar[r] &\Hom_\m(X, T\otimes S(i))\ar[r] &\Hom_\m(P^0, T\otimes S(i))\ar[r] &\Hom_\m(P^1, T\otimes S(i))\ar[r] &\cdots }$$ Note that $\cok_i(X) \in \ ^\bot T$ if and only the upper row is
exact, if and only if the lower row is exact, and if and only if
$X\in \ ^\bot (T\otimes S(i))$.

\vskip5pt

By definition $X\in {\rm smon}(Q, I, \ ^\bot T)$ if and only if
$X\in {\rm smon}(Q, I, A)$ and $\cok_i(X)\in \ ^\bot T, \  \forall \
i\in Q_0$. So by the above claim, $X\in {\rm smon}(Q, I, \ ^\bot T)$ if
and only if $X\in {\rm smon}(Q, I, A)\cap \ ^\bot (T\otimes S)$.
Since ${\rm proj.dim} \ _{kQ/I}S < \infty$,   it is clear that $^\bot (T\otimes S) = \ ^\bot (T\otimes kQ/I)$. Thus $X\in {\rm smon}(Q, I, \ ^\bot T)$ if and
only if $ X\in {\rm smon}(Q, I, A)\cap \ ^\bot (T\otimes kQ/I)$.
This proves the first assertion.

\vskip5pt

Assume that there is an exact sequence $0\rightarrow T_m\rightarrow
\cdots\rightarrow T_0\rightarrow \D(A_A)\rightarrow 0$ with each
$T_j\in {\rm add} (T)$. To show ${\rm smon}(Q, I, \ ^\bot T)= \
^\bot (T\otimes kQ/I)$, by $(1)$ it suffices to show
$^\bot (T\otimes kQ/I)\subseteq \ ^\bot (\D(A_A)\otimes kQ/I)$. This
follows from the exact sequence $0\rightarrow T_m\otimes kQ/I
\rightarrow \cdots \rightarrow T_0\otimes kQ/I \rightarrow
\D(A_A)\otimes kQ/I \rightarrow 0$ with each $ T_j \otimes kQ/I\in
\mbox{add}(T\otimes kQ/I)$. $\s$

\begin{cor} \label{X_i} \ Let $\x$ be an extension-closed subcategory of $A$-mod. If $X\in {\rm smon}(Q, I, \x)$, then $X_i\in \x$ for each $i\in Q_0$.
\end{cor}
\noindent{\bf Proof.} \ For each $i\in Q_0$ and $X\in {\rm smon}(Q, I, \x)$, by Theorem \ref{homodescriptionofsmon}$(1)$ we have
$X_i = (1\otimes e_i)X = (1\otimes e_i)\m\otimes_\m X = (A\otimes e_i(kQ/I))\otimes_\m X\in \x.$
$\s$

\vskip5pt

We stress that  $X_i\in \x$ for each $i\in Q_0$ can {\bf not} replace ${\rm(m3)}$. For example, let $Q$ be the quiver
$2\longrightarrow 1$, $A= kQ$ and $\x = \mathcal P(A)$. Then
$(\binom{1}{0}\hookrightarrow \binom{1}{1})$ satisfies ${\rm(m1)}$
and ${\rm(m2)}$ with $\binom{1}{0}\in \x$ and $\binom{1}{1}\in \x$,
but  $(\binom{1}{0}\hookrightarrow \binom{1}{1})\notin {\rm smon}(Q,
0, \x),$ since the quotient $\binom{1}{1}/\binom{1}{0} \cong
\binom{0}{1}\notin \x$.

\subsection{} The following properties are nontrivial. Since they are not used later, we omit the proof.
\begin{facts} \label{simplerG} \ Let $\x$ be an extension-closed subcategory of
$A$-mod, and $X = (X_i, \ X_\alpha, \ i\in Q_0, \ \alpha\in Q_1)\in
{\rm smon}(Q, I, \x)$. Then

\vskip5pt

$(1)$ \ For each path $p$, \ ${\rm Im}X_p\in\x$.

\vskip5pt

$(2)$ \ For each $\alpha\in Q_1$, \ ${\rm Coker}X_\alpha\in \x$.

\vskip5pt

$(3)$ \ If $\x$  is also closed under kernels of
epimorphisms, then ${\rm Ker}X_p\in \x$ for each path $p$. \end{facts}

\section {\bf Connections with tilting theory and Auslander-Reiten theory}

With the reciprocity in Theorem \ref{homodescriptionofsmon}$(2)$, we will in particular prove that ${\rm smon}(Q, I, A)$ has Auslander-Reiten sequences.

\subsection{Cotilting modules} An $A$-module $T$ is {\it an $r$-cotilting module} (\cite{HR},
\cite{Rin}, \cite{H}, \cite{AR},  \cite{M}), provided that  inj.dim$T\leq r$, ${\rm Ext}^s_A(T, T)=0$ for $s\geq 1$,  and there is an exact sequence $0\rightarrow
T_m\longrightarrow \cdots\longrightarrow T_0\longrightarrow \D(A_A)\rightarrow 0$
with each $T_j\in \mbox{add}(T)$.

\vskip5pt

Following
\cite{AR}, let $\widehat{\x}$ denote the subcategory of $A$-mod
of those $A$-modules $X$ such that there is an exact
sequence $0\rightarrow X_m \rightarrow \cdots\rightarrow
X_0\rightarrow X\rightarrow 0$ with each $X_i\in \x$. The following
lemma is a key in proving Proposition \ref{tensorcotilting}, and
has independent interest. It has been proved by M. Auslander and R.
O. Buchweitz [AB, Prop. 3.5] under different conditions: one of
these conditions is that $\x$ is closed under kernels of
epimorphisms (see [AB, p. 23, line 16]); while this condition can
{\bf not} be satisfied in our application.

\begin{lem} \label{repro} \ Let $\x$ be a self-orthogonal $($i.e.,  ${\rm Ext}^s_A(M, \ N)=0, \ \forall \
M\in \x, \ \forall \ N\in \x, \ \forall \ s\geq 1)$ additive
subcategory of $A$-mod. Then  $\widehat{\x}$ is closed under cokernels of
monomorphisms.
\end{lem}
\noindent{\bf Proof.} \ Let $0 \rightarrow X \stackrel
{f}\longrightarrow Y \stackrel {g}\longrightarrow  Z \rightarrow 0$
be an exact sequence with $X\in\widehat{\mathcal X}$ and $Y\in
\widehat{\mathcal X}$. So there are exact sequences $0\rightarrow
X_n\rightarrow \cdots \rightarrow X_0\xrightarrow {c_0} X\rightarrow
0$ and $0\rightarrow Y_n\xrightarrow {d_n} \cdots \xrightarrow {d_1}
Y_0\xrightarrow {d_0}Y\rightarrow 0$ with each $X_i\in \mathcal X$
and $Y_i\in \mathcal X$. Since $\mathcal X$ is self-orthogonal,
${\rm Ext}^1_A(\x, {\rm Ker}d^i) = 0$ for each $i\ge 0$, and hence
$f: X\rightarrow Y$ induces a chain map $f^\bullet:
X^\bullet\rightarrow Y^\bullet$, where $X^\bullet$ is the complex
$0\rightarrow X_n\rightarrow X_{n-1}\rightarrow \cdots \rightarrow
X_0\rightarrow 0$, and similarly for $Y^\bullet$. Then we get a
morphism of distinguished triangles in the bounded derived category
$\mathcal D^b(A):$
\[\xymatrix@R=0.65cm@C=0.55cm{X^\bullet \ar[r]^{f^\bullet} \ar[d]^{c_0} &Y^\bullet
\ar[r]^-{\binom{0}{{\rm Id}_{Y^\bullet}}} \ar[d]^{d_0}
&{\rm Cone}(f^\bullet) \ar @{.>}[d] \ar[rr]^-{(0, {\rm Id}_{X^\bullet[1]})}&& X^\bullet [1] \ar[d]^-{c_0[1]} \\
X \ar[r]^-f &Y \ar[r]^-g & Z \ar[rr] && X[1]}
\]
\noindent where distinguished triangle at the lower row is induced
by the exact sequence $0 \rightarrow X \stackrel {f}
\rightarrow Y \rightarrow Z \rightarrow 0$, and ${\rm
Cone}(f^\bullet)$ is the mapping cone $0\rightarrow
X_n\rightarrow X_{n-1}\oplus Y_n\rightarrow
\cdots\rightarrow X_0\oplus Y_1\xrightarrow {\partial}
Y_0\longrightarrow 0$. Since $c_0$ and $d_0$ are isomorphisms in
$\mathcal D^b(A)$, we have $Z\cong {\rm Cone}(f^\bullet)$ in
$\mathcal D^b(A)$. It follows that the cohomology group ${\rm H}^i{\rm Cone}(f^\bullet)$ is isomorphic to the $i$-th cohomology group
of the stalk complex $Z$ for each $i\in \Bbb Z$. So ${\rm
Cone}(f^\bullet)$ is exact except at the $0$-th position, and
$Y_0/\Ima\partial \cong Z$ induced by $gd_0: Y_0\longrightarrow Z$.
Thus
$$0\longrightarrow X_n\longrightarrow X_{n-1}\oplus Y_n\longrightarrow
\cdots\longrightarrow X_0\oplus
Y_1\stackrel{\partial}\longrightarrow
 Y_0 \stackrel{gd_0}\longrightarrow Z\longrightarrow 0$$ is exact. This proves $Z\in \widehat{\mathcal X}$. \ $\s$

\begin{prop}\label{tensorcotilting} \ Let $A$ and $B$ be
finite-dimensional $k$-algebras, $T$ an $A$-module, and $L$ a
$B$-module. Then $T\otimes L$ is a cotilting $(A\otimes B)$-module
if and only if \ $_AT$ and \ $_BL$ are cotilting
modules. In this case ${\rm End}_{A\otimes B}(T\otimes L)\cong
{\rm End}_{A}(T)\otimes {\rm End}_{B}(L)$ as algebras.
\end{prop}

\noindent{\bf Proof.}  Assume that $_AT$ and $_BL$ are cotilting modules. Let ${\rm inj.dim} T =r$  with an injective resolution $0\rightarrow T\rightarrow I_0 \rightarrow \cdots
\rightarrow I_r \rightarrow0$,  and ${\rm inj.dim}
L =s$.
Then we get an exact sequence $0\rightarrow T\otimes L\rightarrow I_0
\otimes L\rightarrow\cdots \rightarrow I_r\otimes L\rightarrow0$.
Since ${\rm inj.dim} L = s < \infty$, each ${\rm
inj.dim} (I_i\otimes L) \le s$, and hence ${\rm
Ext}^{r+s+1}_{A\otimes B}(X, T\otimes L) = {\rm Ext}^{s+1}_{A\otimes
B}(X, I_r\otimes L) = 0$ for each $(A\otimes B)$-module $X$. Thus
${\rm inj.dim} (T\otimes L)\leq r+s$.  By the Cartan - Eilenberg isomorphism ([CE, Thm. 3.1, p.209, p.205]) we
have

$${\rm Ext}^{s}_{A\otimes B}(T\otimes
L, T\otimes L)\cong\bigoplus\limits_{p+q = s}({\rm Ext}^p_A(T,
T)\otimes {\rm Ext}^q_B(L, L)), \ \forall \ s\ge 0. \eqno(3.1)
$$ Thus ${\rm Ext}^s_{\m}(T\otimes L, T\otimes L)=0, \
\forall \ s\geq 1$, and  ${\rm End}_{A\otimes B}(T\otimes L)\cong
{\rm End}_{A}(T)\otimes {\rm End}_{B}(L)$.
Put $\x:={\rm add}(T\otimes L)$. To see that $T\otimes L$ is a
cotilting $(A\otimes B)$-module, it remains to prove $\D(A\otimes B)\cong \D(A_A)\otimes \D(B_B)\in
\widehat{\x}$. In fact, since $_BL$ is cotilting, we get
an exact sequence $0\rightarrow L_m \rightarrow\cdots
\rightarrow L_0\rightarrow \D(B_B) \rightarrow0$ with
each $L_i\in {\rm add}(L)$. So we get an exact sequence of $(A\otimes B)$-modules:
$$0\longrightarrow \D(A_A)\otimes L_m \longrightarrow\cdots
\longrightarrow \D(A_A)\otimes L_0\longrightarrow \D(A_A)\otimes
\D(B_B) \longrightarrow 0.\eqno(3.2)$$ Since
$_AT$ is cotilting, we have an exact sequence $0\rightarrow
T_t\rightarrow\cdots \rightarrow T_0 \rightarrow
\D(A_A)\rightarrow0$ with each $T_i \in {\rm add}(T)$. So
$0\rightarrow T_t\otimes L_j\rightarrow\cdots
\rightarrow T_0\otimes L_j \rightarrow \D(A_A)\otimes
L_j\rightarrow 0$ is exact  for $0\leq j\leq m$, with each $ T_i\otimes
L_j\in {\rm add}(T\otimes L)$. Thus each $\D(A_A)\otimes L_j\in
\widehat{\x}$. By Lemma \ref{repro} and $(3.2)$ we have
$\D(A_A)\otimes \D(B_B)\in \widehat{\x}.$

\vskip5pt

Conversely, assume that $T\otimes L$ is cotilting.  Since the tensor product is over field $k$, an
injective resolution of the $(A\otimes B)$-module $T\otimes L$ gives an
injective resolution of the $A$-module $T^{\oplus {\rm dim}_kL}\cong
T\otimes L$, and hence ${\rm inj.dim}(T^{\oplus {\rm dim}_kL}) <
\infty$. Thus ${\rm inj.dim}T < \infty$. By a similar argument, $_AT$ satisfies the third condition of a
cotilting module. By $(3.1)$, ${\rm Ext}^s_A(T, T)=0$ for
$s\geq 1$. Thus $_AT$ is cotilting. Similarly, $_BL$ is cotilting. $\s$

\begin{cor}\label{main} Let $_AT$ be a cotilting module.
Then $T\otimes kQ/I$ is the unique cotilting $\m$-module with ${\rm
End}_{\m}(T\otimes kQ/I)\cong {\rm End}_{A}(T)\otimes (kQ/I)^{op}$,
up to multiplicities of indecomposable direct summands, such that
${\rm smon}(Q, I, \ ^\bot T)= \ ^\bot (T\otimes kQ/I)$.
\end{cor}
\noindent{\bf Proof.}  Since $Q$ is acyclic, $kQ/I$ is a cotilting
$(kQ/I)$-module. By Proposition \ref{tensorcotilting}, $T\otimes
kQ/I$ is a cotilting $\m$-module, and by
Theorem  \ref{homodescriptionofsmon}$(2)$, ${\rm smon}(Q, I, \ ^\bot T)= \ ^\bot (
T\otimes kQ/I)$. If $_\m L$ is another cotilting module such that
$^\bot L={\rm smon}(Q, I, \ ^\bot T)= \ ^\bot (T\otimes kQ/I )$,
then $(T\otimes kQ/I)\oplus L$ is also a cotilting $\m$-module.
Since the number of pairwise non-isomorphic direct summands of a
cotilting module is equal to the number of pairwise non-isomorphic
simple modules (see E. Cline, B. Parshall and L. Scott [CPS, Corol.
2.5]; see also D. Happel [H, p.101], where this result is stated for
algebras of finite global dimension), we get $L\in {\rm
add}(T\otimes kQ/I).$ $\s$

\subsection{Auslander-Reiten sequences}  M. Auslander and I. Reiten [AR, Thm. 5.5(a)] claim that
$\x$ is resolving and contravariantly finite in $A$-mod with
$\widehat{\x}=A$-mod if and only if $\x= \ ^\bot T$ for some
cotilting $A$-module $T$. As an application we get

\begin{thm}\label{ar} \ Let $\x$ be an additive subcategory of $A$-mod. Then ${\rm smon}(Q, I, \x)$ is resolving and
contravariantly finite in $\m$-mod with $\widehat{{\rm smon}(Q, I,
\x)}=\m$-mod if and only if $\x$ is resolving and contravariantly
finite in $A$-mod with $\widehat{\x}=A$-mod. In this case, ${\rm
smon}(Q, I, \x)$ is functorially finite in $\m$-mod, and ${\rm
smon}(Q, I, \x)$ has Auslander-Reiten sequences.

\vskip5pt

In particular, ${\rm smon}(Q, I, A)$ is functorially finite in
$\m$-mod, and ${\rm smon}(Q, I, A)$ has Auslander-Reiten sequences.
\end{thm}
\noindent{\bf Proof.} \ If $\x$ is resolving and contravariantly
finite with $\widehat{\x}=A$-mod, then $\x= \ ^\bot T$ for some
cotilting module $_AT$ ([AR, Thm. 5.5(a)]), and hence by
Corollary \ref{main}, $T\otimes kQ/I$ is a cotilting $\m$-module
with ${\rm smon}(Q, I, \x)={\rm smon}(Q, I, \ ^\bot T)= \ ^\bot
(T\otimes kQ/I)$. Thus ${\rm smon}(Q, I, \x)$ is resolving and
contravariantly finite in $\m$-mod with $\widehat{{\rm smon}(Q, I,
\x)}=\m$-mod, again by [AR, Thm. 5.5(a)].

\vskip5pt

Conversely, assume that ${\rm smon}(Q, I, \x)$ is resolving and
contravariantly finite in $\m$-mod with $\widehat{{\rm smon}(Q, I,
\x)}$ $=\m$-mod. By Lemma \ref{cokexact}$(2)$, $\x$ is a
resolving subcategory of $A$-mod. To see that $\x$ is
contravariantly finite in $A$-mod, we take sink vertex $1$ of $Q_0$
and consider the functor $-\otimes P(1)\colon \ A$-mod $\rightarrow {\rm
smon}(Q, I, A)$ (cf. Example \ref{proj}). For $0\ne L\in
A$-mod, $L\otimes P(1)$ has only one non-zero branch and the first
branch of $L\otimes P(1)$ is $L$. Let $f\colon X \rightarrow
L\otimes P(1)$ be a right ${\rm smon}(Q, I, \x)$-approximation. Then
$f_1\colon X_1 \rightarrow L$ is a right $\x$-approximation (note
that $X_1\in \x$ by Corollary \ref{X_i}$(1)$). Since $L\otimes
P(1)\in\widehat{{\rm smon}(Q, I, \x)}$, we have an exact sequence
$0\rightarrow X^m\rightarrow \cdots \rightarrow X^0\rightarrow
L\otimes P(1) \rightarrow 0$ with each $X^i\in {\rm smon}(Q, I,
\x)$. Taking the first branch we get an exact sequence $0\rightarrow
X^m_1\rightarrow \cdots \rightarrow X^0_1\rightarrow L \rightarrow
0$ with each $X^i_1\in \x$ (cf. Corollary \ref{X_i}$(1)$), which
means $L\in\widehat{\x}$. This proves $\widehat{\x}=A$-mod.

\vskip5pt

Assume that $\x$ is resolving and contravariantly finite in $A$-mod
with $\widehat{\x}=A$-mod. Then ${\rm smon}(Q, I, \x)$ is resolving
and contravariantly finite in $\m$-mod, as we have proven. By H.
Krause and O. Solberg [KS, Corol. 0.3], a resolving contravariantly
finite subcategory of $\m$-mod is functorially finite, and by M.
Auslander and S. O. Smal${\o}$ [AS, Thm. 2.4], an extension-closed
functorially finite subcategory of $\m$-mod has Auslander-Reiten
sequences, so ${\rm smon}(Q, I, \x)$ is functorially finite in
$\m$-mod and ${\rm smon}(Q, I, \x)$ has Auslander-Reiten sequences.
$\s$

\section{\bf Filtration interpretation}

We also have a filtration interpretation of ${\rm smon}(Q, I, \mathscr{X})$. It has important applications later.

\subsection{} Let $\mathcal A$ be an abelian category, and $\mathcal B$ a subcategory  of $\mathcal A$. Deonte by ${\rm Fil}(\mathcal B)$ the subcategory of $\mathcal A$ consisting of objects
which have a (finite) filtration with factors in $\mathcal B$.

\begin{thm}\label{filt} Let $\mathscr{X}$ be  an extension-closed subcategory of $A$-{\rm mod}. Then $${\rm smon}(Q, I, \mathscr{X})={\rm Fil}(\mathscr{X}\otimes \mathcal P(kQ/I)).$$
\end{thm}
\noindent{\bf Proof.} \ Since ${\rm smon}(Q, I, \mathscr{X})$ is extension-closed (cf. Lemma \ref{cokexact}$(2)$), ${\rm Fil}(\mathscr{X}\otimes \mathcal P(kQ/I)) \subseteq {\rm smon}(Q, I, \mathscr{X})$
(cf. Example \ref{proj}).
Let $X\in{\rm smon}(Q, I, \mathscr{X})$. We need to prove that $X$ has a filtration with the factors in $\mathscr{X}\otimes \mathcal P(kQ/I)$. We prove this by induction on $|Q_0|$.

\vskip5pt

Taking a source of $Q$, say the vertex $n$. We claim that there is an exact sequence of left $\m$-modules:
$$0\longrightarrow X_n\otimes P(n) \stackrel {\sigma_n} \longrightarrow X \longrightarrow Y \longrightarrow 0. \eqno(4.1)$$

\vskip5pt

In fact, since $n$ is a source,  $P(n)\otimes e_n(kQ/I) = (kQ/I)e_n \otimes e_n(kQ/I)\cong (kQ/I)e_n(kQ/I)$ as $(kQ/I)$-bimodules.
So we get an exact sequence of $(kQ/I)$-bimodules
$$0\longrightarrow P(n)\otimes e_n(kQ/I)\longrightarrow kQ/I \longrightarrow B'\longrightarrow 0.$$
Thus we get an
exact sequence of $\m$-bimodules $$0\longrightarrow A\otimes P(n)\otimes e_n(kQ/I)\longrightarrow \m \longrightarrow A\otimes B'\longrightarrow 0.$$
Applying  $-\otimes_{\Lambda}X$ we get an exact sequence of left $\m$-modules
$$0\longrightarrow {\rm Tor}_1^{\Lambda}(A\otimes B', X)\longrightarrow (A\otimes P(n)\otimes e_n(kQ/I))\otimes_{\Lambda}X\longrightarrow X \longrightarrow
(A\otimes B')\otimes_{\Lambda}X\longrightarrow 0.$$
Since $X\in{\rm smon}(Q, I, \mathscr{X})$, ${\rm Tor}_1^{\Lambda}(A\otimes B', X)=0$ by Theorem \ref{homodescriptionofsmon}$(1)$ (note that $B'$ is filtrated by $S(i)$ for $i\in Q_0$). So we get
an exact sequence of left $\m$-modules
$$0\longrightarrow (A\otimes P(n)\otimes e_n(kQ/I))\otimes_{\Lambda}X\longrightarrow X \longrightarrow
(A\otimes B')\otimes_{\Lambda}X\longrightarrow 0.$$
By the universal property of tensor products, we observe that $(A\otimes P(n)\otimes e_n(kQ/I))\otimes_{\Lambda}X
\cong (A\otimes e_n(kQ/I))\otimes_{\Lambda} X\otimes P(n)$ as left $\m$-modules, via $(a\otimes p \otimes e_nq) \otimes_{\Lambda} x \mapsto (a\otimes e_nq)\otimes_{\Lambda} x\otimes p.$
Since $e_n$ is a source, $e_n(kQ/I)$ is just the right simple $(kQ/I)$-module at $n$, and hence $e_n(kQ/I)\cong \D S(n)$. Thus we have isomorphisms of left $\m$-modules
\begin{align*}&(A\otimes P(n)\otimes e_n(kQ/I))\otimes_{\Lambda}X
\cong (A\otimes e_n(kQ/I))\otimes_{\Lambda} X\otimes P(n)\\& \cong (A\otimes \D S(n))\otimes_{\Lambda} X\otimes P(n)
= {\rm Cok}_n(X)\otimes P(n) = X_n\otimes P(n),\end{align*}
where the first identity
from Lemma \ref{anothermono}$(2)$. This proves the claim.

\vskip5pt

We further claim that $Y\in {\rm smon}(Q, I, \mathscr{X})$.

\vskip5pt

In fact, since $X\in{\rm smon}(Q, I, A)$, by Theorem \ref{homodescriptionofsmon}$(1)$,
${\rm Tor}^\m_1(A_A\otimes D(S), X)=0$.  Applying $(A_A\otimes \D S)\otimes_\m -$ to $(4.1)$ we get an exact sequence
$$0\longrightarrow {\rm Tor}^\m_1(A_A\otimes \D S, Y) \longrightarrow (A_A\otimes \D S)\otimes_\m (X_n\otimes P(n)) \stackrel {{\rm Id}\otimes_\m \sigma_n}
\longrightarrow (A_A\otimes \D S)\otimes_\m X.$$ By Lemma \ref{anothermono}$(2)$ and $(1.7)$ we have
$$(A_A\otimes \D S)\otimes_\m (X_n\otimes P(n)) = \bigoplus\limits_{i\in Q_0} {\rm Cok}_i(X_n\otimes P(n)) = {\rm Cok}_n(X_n\otimes P(n)) = X_n,$$
$$(A_A\otimes \D S)\otimes_\m X = \bigoplus\limits_{i\in Q_0} {\rm Cok}_i(X) = X_n\oplus (\bigoplus\limits_{i\in Q_0, i\ne n} {\rm Cok}_i(X)).$$ Since ${\rm Id}\otimes_\m -$ preserves direct sum,
it follows that ${\rm Id}\otimes_\m \sigma_n$ is monic, and hence ${\rm Tor}^\m_1(A_A\otimes \D S, Y) = 0$. By Theorem \ref{homodescriptionofsmon}$(1)$, $Y\in{\rm smon}(Q, I, A)$.
By Lemma \ref{cokexact}$(1)$ we have an exact sequence
$$0\longrightarrow {\rm Cok}_i(X_n\otimes P(n)) \longrightarrow {\rm Cok}_i(X) \longrightarrow {\rm Cok}_i(Y) \longrightarrow 0$$
from which we see that
${\rm Cok}_n(Y) =0$ and ${\rm Cok}_i(Y) = {\rm Cok}_i(X)\in \mathscr{X}$ for $i\ne n$ (cf. $(1.7)$). This proves $Y\in {\rm smon}(Q, I, \mathscr{X})$.

\vskip5pt

Put $Q'$ to be the quiver obtained from $Q$ by deleting the vertex $n$, $\rho': = \{\rho_i\in \rho \ | \ s(\rho_i)\ne n\}$, and $I': = \langle \rho_i \ | \ \rho_i\in \rho'\rangle.$
Note that ${\rm smon}(Q', I', \mathscr{X})$ is naturally regarded as a subcategory of ${\rm smon}(Q, I, \mathscr{X})$.
Since the $n$-th branch $Y_n$ of $Y$ is $0$, $Y$ is naturally regarded as an object of  ${\rm smon}(Q', I', \mathscr{X})$. Since $|Q'_0| < |Q_0|$, by the inductive hypothesis that
$Y$ has a filtration with factors in $\mathscr{X}\otimes \mathcal P(kQ'/I') \subseteq \mathscr{X}\otimes \mathcal P(kQ/I)$. It is clear that
$X$ has a filtration with factors in $\mathscr{X}\otimes \mathcal P(kQ/I)$. This completes the proof. $\s$

\subsection{\bf Projective and injective objects}

As a consequence, we can  claim,  in particular that ${\rm smon}(Q, I, \x)$ has enough injective objects if so has $\x$.

\vskip5pt

Let $\mathcal B$ be a subcategory  of an abelian category $\mathcal A$. An object $L\in\mathcal B$ is
{\it injective}, if for any exact sequence $0 \rightarrow X
\stackrel f \rightarrow Y \rightarrow Z \rightarrow 0$ in $\mathcal B$ and any $h\colon X
\rightarrow L$, there exists $g\colon Y \rightarrow L$ such that $h
= g\circ f.$ We say that $\mathcal B$ {\it has enough injective objects},
if for each object $X\in \mathcal B$, there exists an exact sequence
$0 \rightarrow X \rightarrow L \rightarrow Z \rightarrow 0$ in
$\mathcal B$, such that $L$ is an injective
object of $\mathcal B$.  Dually, we say that $\mathcal B$ {\it has enough projective objects}.

\begin{cor} \label{enoughinjectives} Let $\x$ be an extension-closed subcategory of $A$-mod. Then

\vskip5pt

$(1)$ \ ${\rm smon}(Q,
I, \x)$ has enough projective objects if and only if so has $\x$. In  this case, the indecomposable projective objects of ${\rm smon}(Q, I, \x)$ are exactly $M\otimes P(i)$, where $M$ runs over
indecomposable projective objects of $\x$ and $i$ runs over $Q_0;$ and the projective objects of $\x$ are exactly the $i$-th branches of projective objects of
${\rm smon}(Q, I, \x)$, where $i$ is a fixed sink of $Q$.

\vskip5pt

$(2)$ \ ${\rm smon}(Q, I, \x)$ has enough injective objects if and only if so has $\x$. In this case, the indecomposable injective
objects of ${\rm smon}(Q, I, \x)$ are exactly $N\otimes P(i)$, where $N$
runs over indecomposable injective objects of $\x$ and $i$ runs over $Q_0;$ and the injective objects of $\x$ are exactly the
$i$-th branches of injective
objects of ${\rm smon}(Q, I, \x)$, where $i$ is a fixed sink of $Q$.
\end{cor}

\noindent {\bf Proof.} We prove $(2)$. The assertion $(1)$ can be dually proved.
For an injective object $N$ of $\x$ and $W\otimes P\in \x\otimes \mathcal P(kQ/I)$, by the Cartan-Eilenberg isomorphism we have
\begin{align*}&{\rm Ext}^1_\m(W\otimes P, N\otimes P(i))\\ & \cong ({\rm Hom}_A(W, N) \otimes {\rm Ext}^1_{kQ/I}(P, P(i)))\bigoplus ({\rm Ext}^1_A(W, N) \otimes {\rm Hom}_{kQ/I}(P, P(i))) = 0.\end{align*}
Thus
$N\otimes P(i)$ is an injective object of $\x\otimes \mathcal P(kQ/I)$, and hence by Theorem \ref{filt}, $N\otimes P(i)$ is an injective object of ${\rm smon}(Q, I, \x)$.

\vskip5pt

Assume that $\x$ has enough injective objects. For arbitrary $W\in\x$, there is an exact sequence $0\rightarrow W\rightarrow N\rightarrow V\rightarrow 0$ in $\x$ with $N$ an injective object of $\x$.
Then for arbitrary $P\in\mathcal{P}(kQ/I)$, we get an exact sequence $0\rightarrow W\otimes P \rightarrow N\otimes P \rightarrow V\otimes P \rightarrow 0$ in $\x\otimes\mathcal{P}(kQ/I)$.
This shows that $\x\otimes\mathcal{P}(kQ/I)$ has enough injective objects.
By Theorem \ref{filt} and using the Horseshoe Lemma, it is not hard to see that for arbitrary $X\in{\rm smon}(Q, I, \x)$, there is an exact sequence in ${\rm smon}(Q, I, \x)$:
$$0\longrightarrow X\longrightarrow \bigoplus\limits_{i\in Q_0}N_i\otimes P(i)\longrightarrow Y\longrightarrow 0$$
where each $N_i$ is an injective object of $\x$. Thus ${\rm smon}(Q, I, \x)$ has enough injective objects.

\vskip5pt

If $X$ is an indecomposable injective object of ${\rm smon}(Q, I, \x)$, then the above exact sequence splits, and hence $X\cong N\otimes P(i)$ for some indecomposable injective object $N$ of $\x$. Note that
if $V$ is  an indecomposable $A$-module, then $\End_{\m}(V\otimes P(i))\cong\End_A(V)\otimes\End_{kQ/I}(P(i))\cong\End_A(V)$ is a local algebra, and thus $V\otimes P(i)$ is indecomposable.

\vskip5pt

Conversely, assume that ${\rm smon}(Q, I, \x)$ has enough injective objects.  Let $i$ be a fixed sink vertex of $Q$ and $L = (L_j, \ L_\alpha)$ an injective object of
${\rm smon}(Q, I, \x)$. We claim that the
$i$-th branch $L_i$ of $L$ is an injective object of
$\x$. In fact,  for an exact
sequence $0\rightarrow U
\stackrel f\rightarrow V \rightarrow W \rightarrow 0$ in $\x$ and an $A$-map $h_i\colon U\rightarrow L_i$,  we get an exact sequence
$0\rightarrow U\otimes P(i) \stackrel {f\otimes {\rm
Id}}\longrightarrow V\otimes P(i) \rightarrow W\otimes P(i)
\rightarrow 0$ in ${\rm smon}(Q, I, \x)$ and a $\m$-map
$h\colon U\otimes P(i)\rightarrow L$, where the $j$-th
branch $h_j$ of $h$ is $0$ if $j\ne i$, and the $i$-th branch of $h$
is $h_i$ (since $i$ is a sink, $h$ is indeed a
$\m$-map). So there is a $\m$-map $g\colon V\otimes
P(i)\rightarrow L$ with $h = g\circ ({f\otimes {\rm
Id}})$, ann hence $h_i = g_i\circ f.$ This proves the claim. For $U\in \x$, let
$0\rightarrow U\otimes P(i) \stackrel
{f}\longrightarrow L \stackrel {g}\longrightarrow Y\rightarrow 0$ be an exact sequence in ${\rm smon}(Q, I, \x)$
such that $L$ is an injective object of ${\rm smon}(Q, I, \x)$.  This induces an exact sequence
$0\rightarrow U \stackrel {f_i}\longrightarrow L_i \stackrel
{g_i}\longrightarrow Y_i\rightarrow 0$ in $\x$ such that $L_i$ is an injective
object of $\x$, by the claim above. This shows that $\x$ has enough injective
objects. If $N$ is an injective object of $\x$, then $N\otimes P(i)$ is an
injective object of ${\rm smon}(Q, I, \x)$, thus $N$ is the $i$-th branch of
an injective object of ${\rm smon}(Q, I, \x)$.  $\s$

\section{\bf Frobenius subcategories}

\subsection{}  For {\it an exact category} we refer to [Q] and [K1]. An extension-closed subcategory $\mathcal F$ of an abelian
category $\mathcal A$ is an exact category in the canonical way. In this paper we use this exact structure.
An exact category $\mathcal F$ is {\it a Frobenius subcategory} of
$\mathcal A$, if $\mathcal F$ has enough projective
objects and enough injective objects, and $X$ is a projective object
of $\mathcal F$ if and only if it is an injective object of $\mathcal F$.

\vskip5pt

Let $\mathcal A$ be an abelian category with enough projective
objects. {\it A complete $\mathcal A$-projective
resolution} is an exact sequence $P^\bullet: \  \cdots \rightarrow
P^{-1}\rightarrow P^{0} \stackrel{d^0}{\longrightarrow}
P^{1}\rightarrow \cdots$ of projective objects of $\mathcal A$, such
that ${\rm Hom}_\mathcal A(P^\bullet, P)$ is again exact
for each projective object $P$ of $\mathcal A$. An object $M$ is
{\it Gorenstein-projective}, if there is a complete $\mathcal
A$-projective resolution $P^\bullet$ such that $M\cong
\operatorname{Ker}d^0$ ([ABr], [EJ]). Then $\mathcal {GP}(\mathcal A)$  is
resolving and a Frobenius subcategory ([Bel, Prop. 3.8],
[AR, Prop. 5.1], [Hol, Thm. 2.5]). Also, any subcategory of $\mathcal {P}(\mathcal A)$
is a Frobenius subcategory, where $\mathcal {P}(\mathcal A)$ is the subcategory of projective objects
of $\mathcal A$. The following result in particular implies that $\mathcal {GP}(\mathcal A)$ is the largest resolving subcategory which is also a Frobenius subcategory.

\begin{prop}\label{frobenius1} \ Let $\mathcal A$ be an abelian category
with enough projective objects, and $\mathcal F$ an extension-closed subcategory. If $\mathcal F$ is  Frobenius with $\mathcal P(\mathcal F) \subseteq \mathcal P(\mathcal A)$, then $\mathcal F \subseteq \mathcal
{GP}(\mathcal A).$
\end{prop}
\noindent {\bf Proof.}  \ Since $\mathcal F$ has enough
projective objects, for each $X\in\mathcal F$ there
is an exact sequence $0\rightarrow K \rightarrow L
\rightarrow X \rightarrow 0$ in $\mathcal F$ with $L\in \mathcal P(\mathcal F) \subseteq \mathcal P(\mathcal A)$. For each
projective object $P$ of $\mathcal A$, applying $\Hom_A(-, P)$ we
see that ${\rm Ext}^{m+1}_\mathcal A(X, P) \cong {\rm
Ext}^m_\mathcal A(K, P), \ \forall \ m\ge 1.$ Since $K\in
\mathcal F,$ we have
${\rm Ext}^2_\mathcal A(X, P) \cong {\rm Ext}^1_\mathcal A(K, P) = 0,   \
\forall \ X\in\mathcal F.$ Since $K\in \mathcal F,$ we have
${\rm Ext}^3_\mathcal A(X, P) \cong {\rm Ext}^2_\mathcal A(K, P) = 0,   \
\forall \ X\in\mathcal F.$ Continuing this process we get ${\rm
Ext}^m_\mathcal A(X, P) = 0,  \ \forall \ X\in\mathcal F,  \
\forall \ m\ge 1.$

\vskip5pt

Since $\mathcal F$ has enough injective
objects and the injective objects of $\mathcal F$ are also projective objects of $\mathcal A$, we
get an exact sequence
$0 \rightarrow X \rightarrow P^0 \rightarrow P^1 \rightarrow \cdots
\rightarrow P^{i}\stackrel {d^i}\rightarrow \cdots $ with
each $P^i\in\mathcal P(\mathcal A)$ and
${\rm Im} d^i\in \mathcal F$ for all $i\ge 0$. Connecting it with
a projective resolution of $X$ we get an exact sequence
$$\cdots
\longrightarrow P_{1}\longrightarrow P_{0} \longrightarrow P^0
\longrightarrow P^1 \longrightarrow \cdots \longrightarrow
P^{i}\stackrel {d^i}\longrightarrow \cdots $$ Since  ${\rm
Ext}^m_\mathcal A(X, P)= 0$ and ${\rm Ext}^m_\mathcal A({\rm Im}
d^i, P)= 0$ for $m\ge 1, \ i\ge 0,$ and for any projective object
$P$ of $\mathcal A$, this exact sequence is a complete projective resolution. So $X\in \mathcal {GP}(\mathcal A).$  $\s$

\subsection{} The following result follows from Corollary \ref{enoughinjectives}. It gives a way of constructing Frobenius subcategories of $\m$-mod from the one of $A$-mod, via
the correspondence $\x\mapsto {\rm smon}(Q, I, \x)$.

\begin{cor}\label{frobenius2} Let $\x$ be an extension-closed subcategory of $A$-mod. Then
${\rm smon}(Q, I, \x)$ is a Frobenius category if and only if so
is $\x$.
\end{cor}

For $Q = A_2$ (and thus $I = 0$) and $\x = A$-mod,  it was proved in {\rm [C1, 2.1]}.

\vskip5pt

\subsection{} It is clear that ${\rm smon}(Q, I, \x)= {\rm smon}(Q, I, \mathscr Y)$ if and only if
$\x = \mathscr Y$ (to see this, note that if $M\in\x$, then $M\otimes P(1)\in {\rm smon}(Q, I, \x) = {\rm smon}(Q, I, \mathscr Y)$, and hence ${\rm Cok}_1(M\otimes P(1)) = M\in \mathscr Y$)

\vskip5pt

The following examples show that Corollary \ref{frobenius2} gives ``new" Frobenius subcategories, in the sense that
they are not $\mathcal{GP}(\m)$.

\begin{exm} \label{exmfrobenius} $(1)$ \ Let $A$ be the algebra given by the quiver $\xymatrix{1\ar@<1ex>[r]^{\alpha} & 2\ar[r]^{\gamma}\ar@<1ex>[l]^{\beta} & 3}$ with relations $\alpha\beta=0=\beta\alpha$.
Then the indecomposable projective $A$-modules are \ $\begin{smallmatrix} 1\\ 2 \\ 3\end{smallmatrix}$, $\begin{smallmatrix} & 2 & \\ 1 & & 3\end{smallmatrix}$ and $3$,
and the indecomposable injective $A$-modules are \ $\begin{smallmatrix} 2\\ 1\end{smallmatrix}$, $\begin{smallmatrix} 1 \\ 2 \end{smallmatrix}$ and $\begin{smallmatrix} 1 \\ 2 \\ 3\end{smallmatrix}$.
The {\rm Auslander-Reiten} quiver of $A$ is
$$
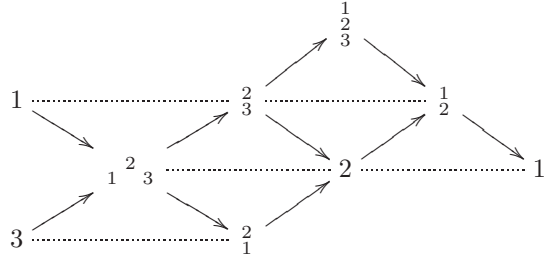
\xymatrix @R= 0.3cm @C=0.8cm {& & & {\begin{smallmatrix} 1\\ 2 \\ 3\end{smallmatrix}}\ar[dr]& & \\
1\ar[dr] & & {\begin{smallmatrix} 2\\ 3 \end{smallmatrix}}\ar[dr]\ar[ur]\ar@{.}[ll] & & {\begin{smallmatrix} 1\\ 2\end{smallmatrix}}\ar[dr]\ar@{.}[ll] & \\
& {\begin{smallmatrix} & 2 & \\ 1 & & 3\end{smallmatrix}}\ar[dr]\ar[ur] & & 2\ar[ur]\ar@{.}[ll] & & 1\ar@{.}[ll]\\
3\ar[ur] & &{\begin{smallmatrix} 2\\ 1\end{smallmatrix}}\ar[ur]\ar@{.}[ll] & & &}$$
where the two vertices $1$ represent the same simple module.
Put $\mathscr{X}: ={\rm add}(1\oplus {\begin{smallmatrix} 1\\ 2\end{smallmatrix}}\oplus 2 \oplus {\begin{smallmatrix} 2\\ 1\end{smallmatrix}})$. By the {\rm Auslander-Reiten} formula
${\rm Ext}^1_A(M, N) \cong \D\underline{{\rm Hom}}_A({\rm Tr}\D N, M)$ we see that  $\x$
is an extension-closed subcategory of $A$-{\rm mod}. By the {\rm Auslander-Reiten} quiver we know that $\x$ is a {\rm Frobenius} subcategory of $A$-{\rm mod}, with indecomposable projective-injective objects exactly being ${\begin{smallmatrix} 1\\ 2\end{smallmatrix}}$ and ${\begin{smallmatrix} 2\\ 1\end{smallmatrix}}$.

\vskip5pt

For any acyclic quiver $Q$ and an arbitrary monomial admissible ideal $I$, applying  {\rm Corollary \ref{frobenius2}} we know that
${\rm smon}(Q, I, \x)$ is a {\rm Frobenius} subcategory of $\m$-{\rm mod}, where $\m = A\otimes kQ/I$.
Since  $\x \ne \mathcal{GP}(A)$, ${\rm smon}(Q, I, \x)\ne {\rm smon}(Q, I, \mathcal{GP}(A)) = \mathcal{GP}(\m)$.
Thus ${\rm smon}(Q, I, \x)$ is a ``new" {\rm Frobenius} subcategory of $\m$-{\rm mod}.

\vskip5pt

Note that $\mathcal P(\x) \nsubseteq \mathcal P(A)$, and $\x\nsubseteq \mathcal{GP}(A)$.  This also shows that the condition
$\mathcal P(\mathcal F) \subseteq \mathcal P(\mathcal A)$ can not be dropped in {\rm Proposition \ref{frobenius1}}.

\vskip5pt

$(2)$ \ Let $A$ be the algebra given by the quiver $\xymatrix{1\ar@<1ex>[r]^{\alpha} & 2\ar@<1ex>[l]^{\beta} & 3\ar[l]_{\gamma}}$ with relations $\alpha\beta=0=\beta\alpha$.
Then the indecomposable projective $A$-modules are \ $\begin{smallmatrix} 1\\ 2 \end{smallmatrix}$, $\begin{smallmatrix} 2 \\ 1 \end{smallmatrix}$ and $\begin{smallmatrix} 3\\ 2\\1 \end{smallmatrix}$,
and the indecomposable injective $A$-modules are \ $\begin{smallmatrix} 3\\ 2\\1 \end{smallmatrix}$, $\begin{smallmatrix} 1 & & 3 \\ & 2 & \end{smallmatrix}$ and $3$.
The {\rm Auslander-Reiten} quiver of $A$ is
$$
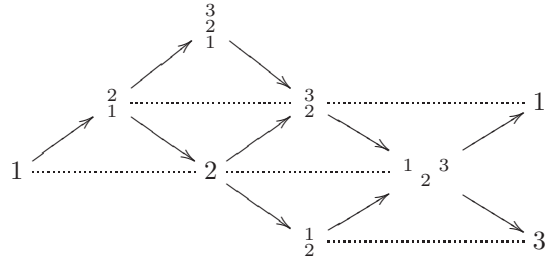
\xymatrix@R= 0.3cm @C=0.8cm {& & {\begin{smallmatrix} 3 \\ 2 \\ 1 \end{smallmatrix}}\ar[dr]& & & \\
& {\begin{smallmatrix} 2 \\1 \end{smallmatrix}}\ar[ur]\ar[dr] & & {\begin{smallmatrix} 3\\ 2 \end{smallmatrix}}\ar[dr]\ar@{.}[ll] & & 1\ar@{.}[ll] \\
1\ar[ur]& & 2\ar[ur]\ar[dr]\ar@{.}[ll] & & {\begin{smallmatrix} 1 & & 3 \\ & 2 & \end{smallmatrix}}\ar[ur]\ar[dr]\ar@{.}[ll] & \\
& & & {\begin{smallmatrix} 1\\ 2 \end{smallmatrix}}\ar[ur]& & 3\ar@{.}[ll]}$$
where the two vertices $1$ represent the same simple module. Put $\mathscr{X}: ={\rm add}(1\oplus {\begin{smallmatrix} 2\\ 1\end{smallmatrix}} \oplus 2 \oplus {\begin{smallmatrix} 1\\ 2\end{smallmatrix}}).$
Then $\x$ is an extension-closed subcategory of $A$-{\rm mod} and it is a {\rm Frobenius} subcategory of $A$-mod,  with indecomposable projective-injective objects exactly being ${\begin{smallmatrix} 1\\ 2\end{smallmatrix}}$ and ${\begin{smallmatrix} 2\\ 1\end{smallmatrix}}$.

\vskip5pt

For any acyclic quiver $Q$ and an arbitrary monomial admissible ideal $I$, by {\rm Corollary \ref{frobenius2}},
${\rm smon}(Q, I, \x)$ is a {\rm Frobenius} subcategory of $\m$-{\rm mod}.
Since  $\x \subsetneqq \mathcal{GP}(A)$, ${\rm smon}(Q, I, \x)\subsetneqq {\rm smon}(Q, I, \mathcal{GP}(A)) = \mathcal{GP}(\m)$.
Thus ${\rm smon}(Q, I, \x)$ is a ``new" {\rm Frobenius} subcategory of $\m$-{\rm mod}.
\end{exm}

\section{\bf The dual version}
We state the dual version of the main points so far, for later applications. A
representation $X = (X_i, \ X_{\alpha})\in {\rm rep}(Q, I, A)$ is {\it separated epic} over $\x$,
if $X$ satisfies the following conditions:

${\rm (e1)}$  \ For $i\in Q_0$,  \ $\Ima (X_i
\stackrel{(X_\alpha) _{\alpha\in \ha(i\to )}}\longrightarrow \bigoplus\limits_{\alpha\in \ha(i\to)}
X_{e(\alpha)})
 = \bigoplus\limits_{\begin {smallmatrix} \alpha\in \ha(
i\to)
\end{smallmatrix}}\Ima X_\alpha,$ where  $\ha(i\to ): = \{\alpha\in Q_1 \ | \
s(\alpha) =i\};$

${\rm (e2)}$  \  For $\alpha\in Q_1$, $\Ima
X_\alpha = \bigcap\limits_{\begin{smallmatrix}q\in
L_\alpha\end{smallmatrix}}\Ker X_q,$  where
$L_\alpha: \ = \{\mbox{non-zero path} \ q \ \mbox{of length}\ge 1 \
| \ s(q) = e(\alpha), \ q\alpha\in I\};$

${\rm (e3)}$  \ For $i\in Q_0$,  \ $\Ker_i(X):=\bigcap\limits_{\begin{smallmatrix}\alpha\in \ha(i\to
)\end{smallmatrix}}\Ker X_\alpha\in \x.$

\vskip5pt

Denote by
${\rm sepi}(Q, I, \x)$ the subcategory of ${\rm rep}(Q, I, A)$ of
separated epic representations over $\x$, and write
${\rm sepi}(Q, I, A)$ for ${\rm sepi}(Q, I, \ A\mbox{-}{\rm mod})$. Let
$\mathcal {I}(A)$ (resp. $\mathcal {GI}(A)$) be the category of
injective (resp. Gorenstein-injective) $A$-modules. Then $\mathcal {I}(\m) =
{\rm sepi}(Q, I, \mathcal {I}(A))$ and
$\mathcal {GI}(\m) = {\rm sepi}(Q, I, \mathcal {GI}(A)).$

\vskip5pt

\noindent {\bf Lemma 2.4'.}  {\it Let $X=(X_i, X_{\alpha}, i\in Q_0, \alpha\in Q_1)\in {\rm
rep}(Q, I, A)$. Then

\vskip5pt

$(1)$ \ $X\in {\rm sepi}(Q, I, A)$ if and only if $$\xymatrix{X_i \ar[rr]^-{(X_\alpha) _{\alpha\in \ha(i\to )}}&&
\bigoplus\limits_{\alpha\in \ha(i\to)}
X_{e(\alpha)}\ar[rr]^-{\bigoplus\limits_{\alpha\in
\ha(i\to)}(X_q)_{q\in L_\alpha }}&&\bigoplus\limits_{\alpha\in
\ha(i\to)}\bigoplus\limits_{q\in L_\alpha} X_{e(q)}}$$ is exact for each $i\in Q_0$.

\vskip5pt

$(2)$ \ For each $i\in Q_0$, we have ${\rm Ker}_i(X) = \Hom_\m(A\otimes S(i), X)$.}

\vskip5pt

\noindent {\bf Lemma 2.5'.} {\it $(1)$  For $i\in Q_0$, the
restriction of  $\Ker_i: \m$-mod $\rightarrow A$-mod to ${\rm
sepi}(Q, I, A)$ is exact.

$(2)$ \ ${\rm sepi}(Q, I, \x)$ is closed under extensions $($resp.
cokernels of monomorphisms$;$ direct summands$)$ if and only if $\x$
is closed under extensions $($resp. cokernels of monomorphisms$;$
direct summands$)$. Thus, ${\rm sepi}(Q, I, \x)$ is a coresolving subcategory of $\m$-mod
if and only if so is $\x$ of $A$-mod.}

\vskip5pt

\noindent {\bf Theorem 2.6'.}   {\it $(1)$ \ We have
\begin{align*}{\rm sepi}(Q, I, A) = \{X\in \m\mbox{-}{\rm mod} \ | \ {\rm Tor}^\m_1(\D X, \ _AA\otimes S) = 0\} = (A_A\otimes \D(kQ/I))^\bot;\end{align*}
\begin{align*}{\rm sepi}(Q, I, \x)& = \{X\in {\rm sepi}(Q, I, A) \ | \ \Hom_\m(A\otimes S(i), X)\in \x, \ \forall \ i\in Q_0\} \\ &
= \{X\in {\rm smon}(Q, I, A) \ | \ \Hom_\m(A\otimes M, X) \ \mbox{for each left} \ (kQ/I)\mbox{-module} \ M\}\end{align*}
where for the second equality we assume that $\x$ is extension-closed. Moreover, if $X\in {\rm sepi}(Q, I, \x)$, then $X_i\in \x$ for each $i\in Q_0$.

\vskip5pt

$(2)$ \ For any $A$-module $T$, \ ${\rm
sepi}(Q, I, T^\bot )= {\rm sepi}(Q, I, A)\cap (T\otimes
D(kQ/I))^\bot$.

If there is an exact sequence $0\rightarrow \ _AA
\rightarrow T_0\rightarrow \cdots\rightarrow
T_m\rightarrow 0$ with each $T_j\in {\rm add} (T)$, then ${\rm
sepi}(Q, I, T^\bot ) = (T\otimes D(kQ/I))^\bot $.}

\vskip5pt

\noindent {\bf Proposition 3.2'.} \ {\it Let $A$ and $B$ be
finite-dimensional $k$-algebras, $T$ an $A$-module, and $L$ a
$B$-module. Then $T\otimes L$ is a tilting $(A\otimes B)$-module if
and only if \ $_AT$ and \ $_BL$ are tilting
modules.}

\vskip5pt

Denote by  $\widetilde{\x}$
the subcategory of $A$-mod of those $A$-modules
$X$ with an exact sequence $0\rightarrow X
\rightarrow X_0 \rightarrow \cdots\rightarrow
X_m\rightarrow 0$ such that each $X_i\in \x$.

\vskip5pt

\noindent {\bf Theorem 3.4'.} \ {\it ${\rm sepi}(Q, I, \x)$ is coresolving
and covariantly finite in $\m$-mod with $\widetilde{{\rm sepi}(Q, I,
\x)}=\m$-mod if and only if $\x$ is coresolving and covariantly
finite in $A$-mod with $\widetilde{\x}=A$-mod. In this case, ${\rm
sepi}(Q, I, \x)$ is functorially finite in $\m$-mod, and ${\rm
sepi}(Q, I, \x)$ has {\rm Auslander-Reiten} sequences. In particular,
${\rm
sepi}(Q, I, A)$ has {\rm Auslander-Reiten} sequences}.

\vskip5pt

\noindent {\bf Theorem 4.1'.} \ {\it Let $\x$ be an extension-closed subcategory of $A$-{\rm mod}. Then
${\rm sepi}(Q, I, \x)={\rm Fil}(\x\otimes\mathcal{I}(kQ/I)).$}

\vskip5pt

\noindent {\bf Corollary 4.2'.} \ {\it Let $\x$ be an extension-closed subcategory of $A$-mod. Then

$(1)$ \ ${\rm sepi}(Q, I, \x)$ has enough projective objects if and only if so has $\x$. In this case,
the indecomposable projective objects of ${\rm sepi}(Q, I, \x)$ are exactly $M\otimes I(i)$, where $M$ runs over indecomposable projective objects of $\x$
and $i$ runs over $Q_0;$ and the projective objects of $\x$ are exactly the
$i$-th branches of projective objects of ${\rm sepi}(Q, I, \x)$, where $i$ is a fixed source of $Q$.

$(2)$ \ \ ${\rm sepi}(Q, I, \x)$ has enough injective objects if and only if so has $\x$. In this case, the indecomposable injective objects of ${\rm sepi}(Q, I, \x)$ are exactly $N\otimes I(i)$, where $N$ runs over indecomposable injective objects of $\x$
and $i$ runs over $Q_0;$ and the injective objects of $\x$ are exactly the
$i$-th branches of injective objects of ${\rm sepi}(Q, I, \x)$, where $i$ is a fixed source of $Q$.}

\vskip5pt

\noindent {\bf Proposition 5.1'.} \ {\it Let $\mathcal A$ be an abelian
category with enough injective objects, and $\mathcal F$ an extension-closed subcategory. If $\mathcal F$ is Frobenius with $\mathcal I(\mathcal{F})\subseteq\mathcal{I}(\mathcal{A})$, then $\mathcal{F}\subseteq\mathcal
{GI}(\mathcal{A})$}.

\vskip5pt

\noindent {\bf Corollary 5.2'.} \ {\it Let $\x$ be an extension-closed subcategory of $A$-mod.  Then ${\rm sepi}(Q, I, \x)$ is a Frobenius category if and only if so
is $\x.$}

\vskip5pt

We also need the following interpretation of ${\rm sepi}(Q, I, \x)$.

\begin{prop}\label{smonsepi} We have ${\rm sepi}(Q, I, \x) = \D{\rm smon}(Q^{op}, I^{op}, \D\x).$
\end{prop}
\noindent {\bf Proof.} Let $X\in \m$-mod. By Lemma 2.4'$(1)$, $X\in {\rm sepi}(Q, I,
A)$ if and only if
$$\xymatrix{\D X_i\ar@{<-}[rr]^-{(\D X_\alpha)_{\alpha\in \ha(i\to)}}&&
\bigoplus\limits_{\alpha\in \ha(i\to)}
\D X_{e(\alpha)}\ar@{<-}[rr]^-{\bigoplus\limits_{\alpha\in
\ha(i\to)}(\D X_q)_{q\in L_\alpha}}&& \bigoplus\limits_{\alpha\in
\ha(i\to)}\bigoplus\limits_{q\in L_\alpha} \D X_{e(q)}}$$
is exact for $i\in Q_0$. Rewrite it $(Q^{op}, I^{op})$. Since $\ha(i\to)$ for $(Q, I)$ is $\ha(\to i)$ for $(Q^{op}, I^{op})$ and
$L_\alpha$ for $(Q, I)$ is just $K_\alpha$ for $(Q^{op}, I^{op})$ and $\D X_i = (\D X)_i$,  $X\in {\rm sepi}(Q, I,
A)$ if and only if $$\xymatrix{\bigoplus\limits_{\alpha\in
\ha(\to i)}\bigoplus\limits_{\begin{smallmatrix}q\in K_\alpha
\end{smallmatrix}} \D X_{s(q)} \ar[rr]^-{\bigoplus\limits_{\alpha\in
\ha(\to i)} (\D X)^\alpha} && \bigoplus\limits_{\begin {smallmatrix}
\alpha\in \ha(\to i)
\end{smallmatrix}}\D X_{s(\alpha)}\ar[rr]^-{\delta_i(\D X)}&& \D X_i},$$
is exact  for $i\in Q_0$, and if and only if $\D X\in {\rm smon}(Q^{op}, I^{op}, A^{op})$.  Also, using Theorem \ref{homodescriptionofsmon}$(1)$ for $(Q^{op}, I^{op})$, for each $\in Q_0$  we have
\begin{align*}{\rm Cok}_i(\D X)\in\D\x & \Longleftrightarrow(A\otimes S(i))\otimes_{\m^{\rm op}}\D X\in\D\x\Longleftrightarrow \D X\otimes_{\m}(A\otimes S(i))\in\D\x
\\&\Longleftrightarrow \D(\D X\otimes_{\m}(A\otimes S(i))) \in \x \stackrel{\mbox{adj. iso.}}\Longleftrightarrow \Hom_\m(A\otimes S(i), X)\in\x \\ &
\stackrel{\mbox{Lem. 2.4'}}\Longleftrightarrow \Ker_i(X)\in\x. \end{align*}
Thus, $X\in {\rm sepi}(Q, I, \x)$ if and only if $\D X\in {\rm smon}(Q^{op}, I^{op}, \D\x)$. $\s$

\section{\bf Ringel-Schmidmeier-Simson equivalence}

Denote by $\mathcal N_A$ the Nakayama functor
$\D\Hom_A(-, \ A)\cong \D(A)\otimes_A-: \ A\mbox{-}{\rm mod}\longrightarrow A\mbox{-}{\rm mod}$.

\begin{defn} \label{rss} \ A {\rm Ringel-Schmidmeier-Simson} equivalence is an equivalence of categories $$F\colon
{\rm smon}(Q, I, \x)\cong {\rm sepi}(Q, I, \x)$$ such
that for $i\in Q_0$ and $M\in \x$, there is a functorial isomorphism
$F(M\otimes P(i)) \cong M\otimes I(i)$ of left $\m$-modules in both
arguments, i.e., $F|_{\x\otimes \mathcal P(kQ/I)}\cong {\rm
Id}_{\x}\otimes \mathcal N_{kQ/I}.$ \end{defn}

\vskip5pt

An {\rm RSS} equivalence implies a strong symmetry. Such an equivalence was first observed in [RS2] and [S1] for
a chain $Q$ with $I = 0$.

\subsection{The existence of an {\rm RSS} equivalence}

\begin{thm} \label{RSS0} \ If $\x$ is an extension-closed subcategory of $A$\mbox{-}{\rm mod}, then there is an {\rm RSS} equivalence \ \  $\D\Hom_\m(-, \D(A_A)\otimes kQ/I)\cong (A_A\otimes \D(kQ/I))\otimes_\m-:
{\rm smon}(Q, I, \x)\cong {\rm sepi}(Q, I, \x).$
\end{thm}
\noindent{\bf Proof.} \ Put $_\m T_\m: = \D(A)\otimes kQ/I$. Following [AR], let $\x_{_ {_{\m}} T}$ be the subcategory of $\m$-mod
consisting of $\m$-modules $_\m X$ such that there is an exact sequence
$$0\longrightarrow X \longrightarrow T_0 \stackrel {f_0}\longrightarrow T_1\longrightarrow \cdots \longrightarrow T_j \stackrel {f_j} \longrightarrow T_{j+1} \longrightarrow \cdots $$
with $T_j\in {\rm add} (_\m T)$ and ${\rm Im}f_j\in \ ^\perp (_\m T)$ for $j\ge 0$. Since $_\m T$ is a cotilting module (cf. Proposition \ref{tensorcotilting}), it follows from
M. Auslander and I. Reiten [AR, Thm. 5.4(b)] that $\x_{_ {_{\m}} T} = \ ^\perp (_{\m} T).$  By Theorem \ref{homodescriptionofsmon}$(1)$,
${\rm smon}(Q, I, A) = \ ^\perp (_{\m} T) = \x_{_ {_{\m}} T}.$

\vskip5pt

It is clear that ${\rm End}_\m(_\m T)^{op} \cong \m$ as algebras, and $T$ has the natural right ${\rm End}_\m(_\m T)^{op}$-module structure; and on the other hand, $T_\m$ is a right $\m$-module. Under the isomorphism
${\rm End}_\m(_\m T)^{op}\cong \m$, $T_{{\rm End}_\m(_\m T)^{op}}$ is exactly $T_\m.$
Note that ${\rm smon}(Q^{op}, I^{op}, A^{op})  \subseteq (A^{op}\otimes kQ^{op}/I^{op})\mbox{-}{\rm mod}
 = (A\otimes kQ/I)^{op}\mbox{-}{\rm mod} = {\rm mod}\m.$
By Theorem \ref{homodescriptionofsmon}$(1)$ we have
\begin{align*}{\rm smon}(Q^{op}, I^{op}, A^{op}) & =  \ ^\perp (\D(A^{op}_{A^{op}})\otimes kQ^{op}/I^{op}) = \ ^\perp (\D(_{A}A)\otimes (kQ/I)_{kQ/I}) = \ ^\perp (T_\m). \end{align*}
For each $X\in {\rm smon}(Q, I, A)$, then $X\in \x_{_ {_{\m}} T}$, and hence by T. Wakamatsu [W, Prop. 1] we have
$${\rm Ext}_{\m^{op}}^m(\Hom_\m(_\m X, \ _\m T), \ T_{\m}) = 0, \ \forall \ m\ge 1$$ and
the canonical left $\m$-map $\ _{\m}X \longrightarrow {\rm Hom}_{\m^{op}}(\Hom_\m(_\m X, \ _\m T), \ T_{\m})$ is an isomorphism.
So we get a contravariant functor $$\Hom_\m(-, \ _\m T): {\rm smon}(Q, I, A) \longrightarrow {\rm smon}(Q^{op}, I^{op}, A^{op}).$$

Similarly,  $T_\m $ is a cotilting module and hence ${\rm smon}(Q^{op}, I^{op}, A^{op}) = \ ^\perp (T_\m) = \x_{T_ {_{\m}}}$.
It is clear that ${\rm End}_{\m}(T_{\m}) \cong \m$ as algebras, and under this isomorphism,
the left module $_{{\rm End}_{\m}(T_{\m})} T$ is exactly $_\m T$.
For each $Y\in {\rm smon}(Q^{op}, I^{op}, A^{op}) = \x_{T_ {_{\m}}}$, again by [W, Prop. 1] we have $${\rm Ext}_{\m}^m(\Hom_{\m^{op}}(Y_\m, \ T_\m), \ _\m T) = 0, \ \forall \ m\ge 1$$ and
the canonical left $\m$-map $Y_{\m}\longrightarrow {\rm Hom}_\m(\Hom_{\m^{op}}(Y_\m, \ T_\m), \ _\m T)$ is an isomorphism. So we also get a contravariant functor
$\Hom_\m(-, T_\m): {\rm smon}(Q^{op}, I^{op}, A^{op}) \longrightarrow {\rm smon}(Q, I, A)$, and moreover $$\Hom_\m(-, \ _\m T): \ {\rm smon}(Q, I, A) \longrightarrow {\rm smon}(Q^{op}, I^{op}, A^{op})$$
is a duality. While by Proposition \ref{smonsepi}, $\D: {\rm smon}(Q^{op}, I^{op}, A^{op}) \longrightarrow {\rm sepi}(Q, I, A)$
is also a duality. Put $F: = \D\Hom_\m(-, \ _\m T) = \D\Hom_\m(-, \D(A_A)\otimes kQ/I)$. Then  $$F:
{\rm smon}(Q, I, A)\longrightarrow {\rm sepi}(Q, I, A)$$ is an equivalence, with a quasi-inverse
$$G: = \Hom_{\m^{op}}(\D(-), \ T_\m)\cong \Hom_{\m}(_AA\otimes \D(kQ/I), -): {\rm sepi}(Q, I, A)\longrightarrow {\rm smon}(Q, I, A).$$
Note that
\begin{align*}F & \cong \D\Hom_{\m^{op}}(A_A\otimes \D(kQ/I), \D(-)) \cong
\D\Hom_{\m^{op}}(A_A\otimes \D(kQ/I), \Hom_k(-, k))\\ & \cong
\D\Hom_k((A_A\otimes \D(kQ/I))\otimes_\m -, k) \cong (A_A\otimes \D(kQ/I))\otimes_\m-. \end{align*}

For $M\otimes L\in A\mbox{-}{\rm mod}\otimes (kQ/I)\mbox{-}{\rm
mod}$, we have functorial isomorphisms of left $\m$-modules:
\begin{align*}F(M\otimes L) &\cong
(A_A\otimes \D(kQ/I))\otimes_\m (M\otimes L)
\\& \cong (A_A\otimes_A M)\otimes (\D(kQ/I)\otimes_{kQ/I} L)\cong M\otimes \mathcal N_{kQ/I} (L).
\end{align*} Thus $F|_{A\mbox{-}{\rm mod}\otimes \mathcal P(kQ/I)}\cong {\rm
Id}_{A\mbox{-}{\rm mod}}\otimes \mathcal N_{kQ/I}.$  It remains to prove $F(X)\in {\rm sepi}(Q, I, \mathscr{X})$ for $X\in{\rm smon}(Q, I, \mathscr{X})$,
and $G(Y)\in {\rm smon}(Q, I, \mathscr{X})$ for $Y\in{\rm sepi}(Q, I, \mathscr{X})$.

\vskip5pt

Since $F:
{\rm smon}(Q, I, A)\longrightarrow {\rm sepi}(Q, I, A)$ is an exact functor between exact categories and ${\rm smon}(Q, I, A)$ is filtrated by
$\x\otimes \mathcal P(kQ/I)$ (cf. Theorem \ref{filt}), to show $F(X)\in {\rm sepi}(Q, I, \mathscr{X})$, it suffices to prove this
for $X = M\otimes P\in \x\otimes \mathcal P(kQ/I).$  That is, $\Hom_\m(A\otimes S(i), F(M\otimes P))\in \x$ for each $i\in Q_0$.
In fact, by the Cartan-Eilenberg isomorphism we have
\begin{align*}
&{\rm Hom}_{\Lambda}(A\otimes S(i), F(M\otimes P)) \cong  {\rm Hom}_{\Lambda}(A\otimes S(i), M\otimes \mathcal N_{kQ/I}(P)) \\
&\cong  {\rm Hom}_A(A, M) \otimes {\rm Hom}_{kQ/I}(S(i), \mathcal N_{kQ/I}(P))
\cong  M \otimes {\rm Hom}_{kQ/I}(S(i), \mathcal N_{kQ/I}(P))\in \x.
\end{align*}

Dually, since $G:
{\rm sepi}(Q, I, A)\longrightarrow {\rm smon}(Q, I, A)$ is an exact functor and ${\rm sepi}(Q, I, A)$ is filtrated by
$\x\otimes \mathcal I(kQ/I)$ (cf. Theorem 4.1'), to show $G(Y)\in {\rm smon}(Q, I, \mathscr{X})$,
it suffices to prove this for $Y = M\otimes L\in \x\otimes \mathcal I(kQ/I).$ That is,
$(A\otimes \D S(i))\otimes_\m G(M\otimes L)\in \x$ for each $i\in Q_0$. In fact,
\begin{align*}
& (A\otimes \D S(i))\otimes_{\Lambda} G(M\otimes L)
= (A\otimes \D S(i))\otimes_{\Lambda}(M\otimes \mathcal N_{kQ/I}^{-}(L))\\ &
\cong (A\otimes_A M)\otimes (\D S\otimes_{kQ/I}\mathcal N_{kQ/I}^{-}(L))\cong  M\otimes (\D S(i)\otimes_{kQ/I}\mathcal N_{kQ/I}^{-}(L))\in\mathscr{X}.
\end{align*}
This completes the proof. $\s$

\vskip5pt

We do not know whether an {\rm RSS}
equivalence $F: {\rm smon}(Q, I, \x)\cong {\rm sepi}(Q, I,
\x)$ is unique, although all the examples we have show the uniqueness.

\subsection{An {\rm RSS} equivalence and the Nakayama functor}

\begin{cor} \label{NakayamaandRSS} \ If $A$ is {\rm Frobenius},  then the restriction of the {\rm Nakayama} functor
$\mathcal N_\m$ gives an {\rm RSS} equivalence ${\rm smon}(Q, I, \x)\cong {\rm sepi}(Q, I, \x)$ for any extension-closed subcategory $\x;$
conversely, if the restriction of $\mathcal N_\m$ gives an {\rm RSS} equivalence ${\rm smon}(Q, I, \x)\cong {\rm sepi}(Q, I, \x)$, then
$\mathcal N_A|_{\x} \cong {\rm Id}_{\x}$.

\vskip5pt

Thus, if $\x \supseteq \mathcal P(A)$, then the restriction of $\mathcal N_\m$ gives an {\rm RSS} equivalence ${\rm smon}(Q, I, \x)\cong {\rm sepi}(Q, I, \x)$  if and only if $A$ is Frobenius.
In particular, the restriction of $\mathcal N_\m$ gives an {\rm RSS} equivalence ${\rm smon}(Q, I, A)\cong {\rm sepi}(Q, I, A)$  if and only if $A$ is Frobenius.
\end{cor}
\noindent{\bf Proof.} \ Assume that $A$ is a Frobenius algebra. Then  $\D(A_A) \cong \ _AA$, and hence
$$\mathcal N_\m =  \D\Hom_\m(-, \ _\m\m) = \D\Hom_\m(-, \ _AA\otimes kQ/I) \cong \D\Hom_\m(-, \ \D(A_A)\otimes kQ/I).$$ Thus the assertion follows from Theorem \ref{RSS0}.

\vskip5pt

For $M\otimes L\in A\mbox{-}{\rm mod}\otimes (kQ/I)\mbox{-}{\rm
mod}$, we have functorial isomorphisms of left $\m$-modules:
\begin{align*}&\mathcal N_\m(M\otimes L)  =
\D\Hom_\m(M\otimes L, A\otimes kQ/I)
\\& \cong \D\Hom_A(M, A)\otimes \D\Hom_{kQ/I}(L,
kQ/I)\cong \mathcal N_A (M) \otimes \mathcal N_{kQ/I} (L).
\end{align*} This shows the fact: $\mathcal N_\m|_
{A\mbox{-}{\rm mod}\otimes (kQ/I)\mbox{-}{\rm mod}} = \mathcal N_A\otimes \mathcal N_{kQ/I}.$

\vskip5pt

Conversely, assume that the restriction of $\mathcal N_\m$ gives an {\rm RSS} equivalence ${\rm smon}(Q, I, \x)\cong {\rm sepi}(Q, I, \x)$. Then for
each indecomposable left $A$-module $M\in \x$, by the definition of an {\rm RSS} equivalence we have $\mathcal N_\m(M\otimes kQ/I)\cong M\otimes \mathcal N_{kQ/I}(kQ/I)$.
By the fact above, this means that
we have a left $\m$-module
isomorphism $\mathcal N_{A}(M)\otimes \D(kQ/I)\cong M\otimes \D(kQ/I).$
It is also a left $A$-module
isomorphism, i.e.,  $\mathcal
N_{A}(M)^{\oplus {\rm dim}_k kQ/I} \cong M^{\oplus {\rm dim}_k kQ/I}$. So we get a functorial left $A$-module
isomorphism $\mathcal N_{A}(M)\cong M$, i.e., $\mathcal N_A|_{\x} \cong {\rm Id}_{\x}$.

\vskip5pt

If $\x \supseteq \mathcal P(A)$, then $\mathcal N_A|_{\x} \cong {\rm Id}_{\x}$ implies $\mathcal N_A(_AA)\cong \ _AA,$ i.e., $A$ is a Frobenius algebra. $\s$

\vskip5pt

\begin{rem} If $A$ is selfinjective which is not {\rm Frobenius} $($so $A$ is not basic$)$, we consider the basic algebra $A'$ of $A$. Then $A'$ is {\rm Frobenius} and we have
$G\colon A$-mod $\cong$ $A'$-mod. There are equivalences ${\rm smon}(Q, I, \x)\cong {\rm smon}(Q, I, G\x)$
and ${\rm sepi}(Q, I, \x)\cong {\rm sepi}(Q, I, G\x)$, given by $G$ componentwise
{\rm (cf. Remark \ref{morita}} and its dual$)$. By {\rm Corollary \ref{NakayamaandRSS}}, we get an {\rm RSS} equivalence ${\rm smon}(Q, I, \x)\cong {\rm sepi}(Q, I, \x)$ for any extension-closed subcategory $\x$, given by
$${\rm smon}(Q, I, \x)\cong {\rm smon}(Q, I, G\x) \stackrel{\mathcal N_{\m'}} \cong {\rm sepi}(Q, I, G\x) \cong {\rm sepi}(Q, I, \x).$$
\end{rem}

\subsection{The case of a chain: A combinatorial {\rm RSS} equivalence} If $Q$ is a chain with $I = 0$, then C. M.
Ringel and M. Schmidmeier [RS2, 1.2] and D. Simson [S1, Chap. 5,
Sect.2] have observed that there is an {\rm RSS} equivalence
${\rm Cok}: {\rm smon}(Q, 0, A)\cong {\rm sepi}(Q, 0, A)$
with quasi-inverse ${\rm Ker}: {\rm sepi}(Q, 0, A)\cong {\rm smon}(Q, 0,
A)$,  where ${\rm Cok}$ sends $X = (X_n \stackrel
{\varphi_{n-1}}\hookrightarrow X_{n-1} \hookrightarrow \cdots
\hookrightarrow X_2\stackrel {\varphi_1}\hookrightarrow X_1) \in
{\rm smon}(Q, 0, A)$ to
$$(X_1\twoheadrightarrow \Cok(\varphi_1\cdots\varphi_{n-1})\twoheadrightarrow \cdots \twoheadrightarrow\Cok(\varphi_1\varphi_2)\twoheadrightarrow
\Cok\varphi_1)\in{\rm sepi}(Q, 0, A)\eqno(7.1)$$ and ${\rm Ker}$ sends $Y = (Y_n \stackrel
{\psi_{n-1}}\twoheadrightarrow Y_{n-1} \hookrightarrow \cdots
\twoheadrightarrow Y_2\stackrel {\psi_1}\twoheadrightarrow Y_1) \in
{\rm sepi}(Q, 0, A)$ to
$$(\Ker\psi_{n-1} \hookrightarrow
\Ker(\psi_{n-2}\psi_{n-1})\hookrightarrow \cdots \hookrightarrow
\Ker(\psi_1\cdots\psi_{n-1}) \hookrightarrow Y_n)\in {\rm smon}(Q, 0, A).\eqno(7.2)$$
It also gives an {\rm RSS} equivalence ${\rm sepi}(Q, 0, \x)\cong {\rm smon}(Q, 0,
\x)$ for an arbitrary additive subcategory $\x$ of $A$-mod.

\vskip5pt

In general (when $I\ne 0$),  the ``jumping of the functor ${\rm Cok}$" also gives an {\rm RSS} equivalence.
The following result does not assume that $\x$ is extension-closed, so it is not a consequence of Theorem \ref{RSS0};
also,  this {\rm RSS} equivalence is given combinatorially and hence operable. We do not know whether
it coincides with the equivalence given in Theorem \ref{RSS0}.

\begin{thm} \label{RSS1} \ Let $Q$ be a chain, $I$ an arbitrary admissible ideal of $kQ$, $A$ a finite-dimensional algebra, and
$\x$ an arbitrary  additive subcategory  of $A$-mod. Then there is an {\rm RSS}
equivalence $F: {\rm smon}(Q, I, \x)\cong {\rm sepi}(Q, I,
\x)$, given by $(7.5)$ below.\end{thm}

\noindent {\bf Proof.} \ Assume that $I\ne 0$. Write the chain as
$$Q = n \stackrel{\alpha_{n-1}}\longrightarrow \cdots \longrightarrow v_1+m_1 \stackrel{\alpha_{v_1+m_1-1}}\longrightarrow
\cdots \stackrel{\alpha_{v_1}}\longrightarrow v_1
\longrightarrow\cdots \stackrel{\alpha_1}\longrightarrow 1.$$ Let
$\{\rho_1, \cdots, \rho_t\}$ be the set of minimal generators of $I$
with $t\ge 1$, and
$$\rho_1= \alpha_{v_1} \cdots \alpha_{v_1+m_1-1},  \ \ \cdots,  \ \
\rho_t = \alpha_{v_t} \cdots \alpha_{v_t+m_t-1}$$ with each $m_i\ge
2$. Different relations $\rho_i$ and $\rho_j$ may overlap, but one
can not contain another. So
$$v_1 \ge 1, \ v_{j+1}> {\rm max}\{v_{j},
v_{j}+m_{j}-m_{j+1}\} \ \mbox{for} \ 1\le j\le t-1, \ v_t+m_t\le
n.$$ Note that $v_{j+1}-v_j-m_j < 0$ if and only if $\rho_{j+1}$ and
$\rho_{j}$ overlap; and in this case, $v_j+m_j-v_{j+1}$ is the
number of the overlapped arrows of $\rho_{j+1}$ and $\rho_{j}$.
Put $u_j: = v_j+m_j-1$ for $1\le j\le t.$ We use {\bf conventions}:
$v_0: = 0, \ u_0: = 0, \ v_{t+1}: = n, \ u_{t+1}: = n,$ \
$\co:=\Cok$ and $\K:=\Ker$.

\vskip5pt

A $\m$-module $X = (X_n
\stackrel{\varphi_{n-1}}\hookrightarrow \cdots \longrightarrow X_{v_1+m_1} \stackrel{\varphi_{u_1}}\longrightarrow
\cdots \stackrel{\varphi_{v_1}}\longrightarrow
X_{v_1}\hookrightarrow \cdots \stackrel{\varphi_1}\hookrightarrow
X_1)$ is in ${\rm smon}(Q, I, A)$ if and only if $$\varphi_i \
\mbox{is monic for} \ i\notin \{v_1, \cdots, v_t\}, \ \ \mbox{and} \
\ {\rm K} (\varphi_{v_j}) = {\rm Im}(\varphi_{v_j+1}\cdots
\varphi_{u_j}) \ \mbox{for} \ 1\le j\le t; \eqno(7.3)$$
and $X\in {\rm smon}(Q, I, \x)$ if and only if $X$ satisfies $(7.3)$, $X_n\in \x$, and
${\rm Ck}(\varphi_i)\in \x$ for $1\le i\le n-1$.
A $\m$-module  $Y =
(Y_n \stackrel{\psi_{n-1}}\longrightarrow\cdots
\longrightarrow  Y_{v_1+m_1}
\stackrel{\psi_{u_1}}\longrightarrow \cdots
\stackrel{\psi_{v_1}}\twoheadrightarrow Y_{v_1}\twoheadrightarrow
\cdots \stackrel{\psi_1}\twoheadrightarrow Y_1)$ is in ${\rm
sepi}(Q, I, \chi)$ if and only if $$\psi_i \ \mbox{is epic for} \
i\notin \{u_1, \cdots, u_t\}, \ \ \mbox{and} \ \  {\rm
Im}(\psi_{u_j}) = {\rm K}(\psi_{v_j}\cdots \psi_{u_{_j}-1})\
\mbox{for} \ 1\le j\le t;  \eqno (7.4)$$
and $Y\in {\rm sepi}(Q, I, \x)$ if and only if $Y$ satisfies $(7.4)$, $Y_1\in \x$, and
${\rm K}(\psi_i)\in \x$ for $1\le i\le n-1$.

\vskip5pt

{\bf Step 1.} \ For $X = (X_i, \varphi_i)\in {\rm smon}(Q, I, A)$,
define $FX = ((FX)_i, \psi_i)$ as follows.

\vskip5pt

The $u_{_{j+1}}$-th branch $(FX)_{u_{_{j+1}}}$ of $FX$ is defined to be
$X_{v_{_{j}}+1}$ for $0\le j\le t$. Thus, the $u_{_1}$-th
branch  of $FX$ is $X_{_1}$, and the $n$-th branch of $FX$ is
$X_{v_{_t}+1}$.

\vskip5pt

If $i$ is in the interval $[u_{_j}+1, u_{_{_{j+1}}}-1]$, $0\le j\le
t$, then the $i$-th branch $(FX)_i$ of $FX$ is defined to be $\co
(\varphi_{v_{_j}+1}\cdots \varphi_{_i})$, where
$X_{_{i+1}}\stackrel{\varphi_{_i}}\longrightarrow X_{_i}
\longrightarrow \cdots \stackrel{\varphi_{v_{_j}+1}}\longrightarrow
X_{v_{_j}+1}$.  Thus
\begin{align*}FX: = (&X_{v_{_t}+1}\stackrel{\psi_{_{n-1}}}\twoheadrightarrow\co(\varphi_{v_{_t}+1}\cdots
\varphi_{_{n-1}})\twoheadrightarrow \cdots \ \ \ \ \ \ \ \ \ \
(\mbox{\small {no this line if}} \ v_{_t}+m_{_t}=n)
\\ \co(\varphi_{v_{_t}+1}\cdots
\varphi_{v_{_t}+m_{_t}})\stackrel{\psi_{u_{_t}}}\longrightarrow
&X_{v_{_{t-1}}+1} \twoheadrightarrow
\co(\varphi_{v_{_{t-1}}+1}\cdots
\varphi_{u_{_t}-1})\twoheadrightarrow \cdots
\\ \co(\varphi_{v_{_2}+1}\cdots \varphi_{v_{_2}+m_{_2}})\stackrel{\psi_{u_{_2}}}\longrightarrow &X_{v_{_1}+1}
\twoheadrightarrow \co(\varphi_{v_{_1}+1}\cdots
\varphi_{u_{_2}-1})\twoheadrightarrow \cdots
\\ \co(\varphi_{v_{_1}+1}\cdots \varphi_{v_{_1}+m_{_1}})\stackrel{\psi_{u_{_1}}}\longrightarrow
& X_{_1} \twoheadrightarrow \co(\varphi_1\cdots
\varphi_{u_{_1}-1})\twoheadrightarrow \cdots \twoheadrightarrow
\co(\varphi_{_1}\varphi_{_2}) \stackrel{\psi_{_1}}\twoheadrightarrow
\co(\varphi_{_1})) \ \ \ \ \ \ \ \ \ \ \ \ \ \ \ \ \
(7.5)\end{align*} where each $\psi_i$ is the natural epimorphism for
$i\notin \{u_1, \cdots, u_t\}$, and  for $1\le j\le t$
$$\psi_{u_j}: = \widetilde{\varphi_{v_{_{j-1}}+1}\cdots \varphi_{v_{_j}}}:
\Cok(\varphi_{v_{_j}+1}\cdots
\varphi_{v_{_j}+m_{_j}})\longrightarrow X_{v_{_{j-1}}+1}\eqno(7.6)$$
induced by $\varphi_{v_{_{j-1}}+1}\cdots \varphi_{v_j}:
X_{v_j+1}\longrightarrow X_{v_{_{j-1}}+1}$, since
$(\varphi_{v_{{_{j-1}}}+1}\cdots \varphi_{v_{_j}})(\varphi_{v_{_j}+1}\cdots
\varphi_{u_{_j}}\varphi_{v_{_j}+m_{_j}}) = 0.$

\vskip5pt

{\bf Remark:} \ If $u_{_{j+1}} -1< (v_{_j}+m_{_j})$ for some $1\le
j \le t$, then $(u_{_{j+1}}-1) - (v_{_j}+m_{_j}) = -1$, since
$v_{_{j+1}} + m_{_{j+1}}> v_{_j}+m_{_j}$. That is,  $u_{_{j+1}} =
u_{_j}+1$. So the $u_{_j}$-th branch $X_{v_{_{j-1}}+1}$ neighbors
with the $u_{_{j+1}}$-th branch $X_{v_{_j}+1}$,  and the
corresponding part of $(7.5)$ reads as
$$X_{v_{_j}+1}\stackrel{\psi_{u_{_j}}:=\varphi_{v_{_{j-1}}+1}\cdots
\varphi_{v_{_j}}}\longrightarrow X_{v_{_{j-1}}+1}.$$ In particular,
if $v_{_t} + m_{_t} = n$, then the first line of $(7.5)$ disappears,
and $\co(\varphi_{v_{_t}+1}\cdots \varphi_{v_{_t}+m_{_t}})$ in the
$t$-th relation does not make sense, and it is $X_{v_{_t}+1}$.

If $(u_{_{j+1}}-1) - (v_{_j}+m_{_j}) = 0$ for some $1\le j \le t$,
then there is only one $A$-modules between $X_{v_{_{j-1}}+1}$ and
$X_{v_{_j}+1}$. In particular,  if $v_{_t} + m_{_t} = n-1$, then the
first line of $(7.5)$ contains only one module $X_{v_{_t}+1}$.

\vskip5pt

{\bf Step 2.} \ For $X = (X_i, \varphi_i)\in {\rm smon}(Q, I, A)$, \
$FX = ((FX)_i, \psi_i)\in {\rm sepi}(Q, I, A)$. Thus, we get a functor
$F: {\rm smon}(Q, I, A)\longrightarrow {\rm sepi}(Q, I, A)$.

\vskip5pt

By construction $\psi_i$ is epic for $i\notin \{u_{_1}, \cdots,
u_{_t}\}$. We need to show ${\rm Im}(\psi_{u_{_j}}) = {\rm
K}(\psi_{v_{_j}}\cdots \psi_{u_{_j}-1})$ for $1\le j\le t$. By
$(7.6)$, ${\rm Im}(\psi_{u_{_j}}) = {\rm
Im}(\varphi_{v_{_{j-1}}+1}\cdots \varphi_{v_{_j}})$, and it is ${\rm
K}(\psi_{v_{_j}}\cdots \psi_{u_{_j}-1})$,  by taking $i = v_j$ in
{\rm Sublemma 1} below.

\vskip5pt

The following sublemma reveals a strong symmetry between ${\rm
smon}(Q, I, A)$ and ${\rm sepi}(Q, I, A)$.

\vskip5pt

{\rm \bf Sublemma 1.} \ Let $ X = (X_i, \varphi_i)\in {\rm smon}(Q,
I, A)$, and $i\in [v_{_{j-1}}+1, u_{_{j}}-1]$ for some $1\le j\le
t+1$. Then $FX = ((FX)_i, \psi_i)$ satisfies  ${\rm
K}(\psi_{_i}\cdots \psi_{u_{_{j}}-1}) = {\rm
Im}(\varphi_{v_{_{j-1}}+1}\cdots \varphi_{_i}).$

\vskip5pt

\noindent{\bf Proof of Sublemma 1.} \ If $[i,
u_{_j}-1]\cap \{u_{_1}, \cdots, u_{_t}\} = \emptyset$, then by {\bf
Step 1} all the maps $\psi_{i}, \ \cdots, \ \psi_{u_{_j}-1}$ are the natural epimorphisms induced by cokernels. By assumption $u_{_{j-1}}\notin [i,
u_{_j}-1]$, so $i\in [u_{_{j-1}}+1, u_{_{j}} -1]$, and hence
by {\bf Step 1}, $(FX)_{i} =
\co(\varphi_{v_{_{j-1}}+1}\cdots\varphi_{i}).$ Also $(FX)_{u_{_j}} =
X_{v_{_{j-1}}+1}.$ Thus
\begin{align*}{\rm K}(\psi_{i}\cdots \psi_{u_{_j}-1}) &=
\K(X_{v_{_{_{j-1}}+1}} \stackrel{\psi_{u_{_j}-1}}\twoheadrightarrow
\co(\varphi_{v_{_{j-1}}+1}\cdots\varphi_{u_{_j}-1})
\twoheadrightarrow \cdots \stackrel{\psi_{i}}\twoheadrightarrow
\co(\varphi_{v_{_{j-1}}+1}\cdots\varphi_{i}))\\&= {\rm
Im}(\varphi_{v_{_{_{j-1}}}+1}\cdots \varphi_{i}).\end{align*}

\vskip5pt

Assume $[i, u_{_j}-1]\cap \{u_1,
\cdots, u_t\} = \{u_{_{j-a}}, \cdots, u_{_{j-1}}\} \ (a\ge 1).$ Note
that $i\ge u_{j-a-1}+1$ (otherwise $u_{j-a-1}$ is also in $[i,
u_{_j}-1]$). While $i\le u_{_{j-a}}.$ Thus either $i\in
[u_{j-a-1}+1, u_{j-a}-1]$ or $i = u_{_{j-a}}.$

\vskip5pt

If $i\in [u_{j-a-1}+1, u_{j-a}-1]$, then by {\bf Step 1} $\psi_{_i}:
(FX)_{i+1}\twoheadrightarrow (FX)_{i} =
\co(\varphi_{v_{_{j-a-1}}+1}\cdots \varphi_{i}).$ Since all the maps
in $\{\psi_{_i}, \cdots, \psi_{u_{_{j-a}}}, \cdots,
\psi_{u_{_{j-1}}}, \cdots, \psi_{u_{_j}-1}\}$ are the natural epimorphisms induced by cokernels,
except those maps $\psi_{u_{j-b}}$ for $1\le b\le a$,
we have (here we use the {\bf convention}: the image $\bar x$ of an element $x$ under the natural epimorphism induced by cokernel is still simply denoted by $x$)
\begin{align*}&{\rm K}(\psi_{_i}\cdots \psi_{u_{_{j-a}}} \cdots
\psi_{u_{_{j-1}}} \cdots \psi_{u_{_j}-1}) \\& = \{x\in
X_{v_{_{j-1}+1}} = (FX)_{u_{_j}}\ | \ (\psi_{u_{_{j-a}}} \cdots
\psi_{u_{_{j-1}}})(x)\in {\rm Im}(\varphi_{v_{_{j-a-1}}+1}\cdots
\varphi_{i})\}\end{align*} where $\psi_{u_{j-a}} \cdots
\psi_{u_{j-1}}$ is the composition of only those maps
$\psi_{u_{j-b}}$ for $1\le b\le a$, namely, all the other maps,
which are the natural epimorphisms induced by cokernels, can be
taken off. By $(7.6)$ it is
\begin{align*}&\{x\in X_{v_{_{j-1}+1}} \ | \
(\varphi_{v_{_{j-a-1}}+1}\cdots \varphi_{v_{_{j-a}}})\cdots
(\varphi_{v_{_{j-2}}+1}\cdots \varphi_{v_{_{j-1}}})(x)\in {\rm
Im}(\varphi_{v_{_{j-a-1}}+1}\cdots \varphi_{i})\}. \ \ \ \
(*)\end{align*}Since all the maps before $\varphi_{v_{_{j-a}}}$ are
monic, they can be taken off. So the set $(*)$ become
\begin{align*}&\{x\in
X_{v_{_{j-1}+1}} \ | \
 \varphi_{v_{_{j-a}}}(\varphi_{v_{_{j-a}}+1}\cdots \varphi_{v_{_{j-a+1}}})\cdots
(\varphi_{v_{_{j-2}}+1}\cdots \varphi_{v_{_{j-1}}})(x)\in {\rm
Im}(\varphi_{v_{_{j-a}}}\cdots \varphi_{i})\}
\\&
= \{x\in X_{v_{_{j-1}+1}}  |
 (\varphi_{v_{_{j-a}}+1}\cdots \varphi_{v_{_{j-a+1}}})\cdots
(\varphi_{v_{_{j-2}}+1}\cdots \varphi_{v_{_{j-1}}})(x)\in {\rm
Im}(\varphi_{v_{_{j-a}}+1}\cdots \varphi_{i}) +
\K(\varphi_{v_{_{j-a}}})\}
\\& \stackrel{(7.3)}= \{x\in
X_{v_{_{j-1}+1}}  |
 \varphi_{v_{_{j-a}}+1}\cdots
\varphi_{v_{_{j-1}}}(x)\in {\rm Im}(\varphi_{v_{_{j-a}}+1}\cdots
\varphi_{i}) + {\rm Im}(\varphi_{v_{_{j-a}}+1}\cdots
\varphi_{u_{_{j-a}}})\}\end{align*} Since ${\rm
Im}(\varphi_{v_{_{j-a}}+1}\cdots \varphi_{i}\cdots
\varphi_{u_{_{j-a}}}) \subseteq {\rm
Im}(\varphi_{v_{_{j-a}}+1}\cdots \varphi_{i})$, the set $(*)$ is $$\{x\in X_{v_{_{j-1}+1}} |
 (\varphi_{v_{_{j-a}}+1}\cdots \varphi_{v_{_{j-a+1}}})\cdots
(\varphi_{v_{_{j-2}}+1}\cdots \varphi_{v_{_{j-1}}})(x)\in {\rm
Im}(\varphi_{v_{_{j-a}}+1}\cdots \varphi_{i})\}.$$ Continuing this
process  we finally get ${\rm K}(\psi_{_i}\cdots \psi_{u_{_{j-a}}}
\cdots \psi_{u_{_{j-1}}} \cdots \psi_{u_{_j}-1}) = {\rm
Im}(\varphi_{v_{_{j-1}}+1}\cdots \varphi_{_i}).$

\vskip5pt

If $i = u_{_{j-a}},$ then by {\bf Step 1} $\psi_{_i}=
\psi_{u_{_{j-a}}}: (FX)_{i+1}\longrightarrow (FX)_{i} =
X_{v_{_{j-a-1}}+1},$ and  we are computing ${\rm
K}(\psi_{u_{_{j-a}}} \cdots \psi_{u_{_{j-1}}} \cdots
\psi_{u_{_j}-1}).$ By the same arguments (in particular, using
$\K(\varphi_{v_{_{j-a}}}) = {\rm Im}(\varphi_{v_{_{j-a}}+1}\cdots
\varphi_{u_{_{j-a}}})$ in $(7.3)$) we get ${\rm K}(\psi_{u_{_{j-a}}}
\cdots \psi_{u_{_{j-1}}} \cdots \psi_{u_{_j}-1}) = {\rm
Im}(\varphi_{v_{_{j-1}}+1}\cdots \varphi_{u_{_{j-a}}}).$

\vskip5pt

This completes the proof of Sublemma 1.

\vskip5pt

{\bf Step 3.} \ If $X = (X_i, \varphi_i)\in {\rm smon}(Q, I, \x)$, \
then $FX = ((FX)_i, \psi_i)\in {\rm sepi}(Q, I, \x)$.

\vskip5pt

By assumption $X_n\in\x$ and  ${\rm Ck}(\varphi_i)\in \x$  for $1\le i\le n-1$. We need to prove
$(FX)_1\in \x$ and ${\rm K}(\psi_i)\in \x$ for $1\le i\le n-1$.

\vskip5pt

In fact, by construction $(FX)_1 = {\rm Ck}(\varphi_1)\in \x$.

\vskip5pt

By construction $\psi_{u_{_i}-1}: X_{v_{_{i-1}}+1}\twoheadrightarrow {\rm Ck}(\varphi_{v_{_{i-1}}+1}\cdots \varphi_{u_{_{i}}-1})$ for $i\in\{1, \cdots, t+1\}$.
So
${\rm K}(\psi_{u_{_i}-1})
= {\rm Im}(\varphi_{v_{_{i-1}}+1}\cdots \varphi_{u_{_{i}}-1})$. Since $\x$ is not assumed to be closed under extensions, we can not apply Remark \ref{simplerG}$(1)$ to get it is in $\x$. We claim that it is isomorphic to
${\rm Ck}(\varphi_{u_i})\in \x$ via
$$f: \varphi_{v_{_{i-1}}+1}\cdots \varphi_{v_i}\cdots\varphi_{u_{_{i}}-1}(x_{u_i})\mapsto x_{u_i}+{\rm Im}(\varphi_{u_i}).$$
To see this, it suffices to show that $f$ is well-defined (clearly $f$ is monic and epic, if it is well-defined).
Assume that $\varphi_{v_{_{i-1}}+1}\cdots \varphi_{v_i}\cdots\varphi_{u_{_{i}}-1}(x_{u_i}) = 0
$. We need to prove $x_{u_i}\in {\rm Im}(\varphi_{u_i})$.  By $(7.3)$, taking off the monomorphisms we get $\varphi_{v_i}\cdots\varphi_{u_{_{i}}-1}(x_{u_i}) = 0.$ Again by $(7.3)$,
$\varphi_{v_{_i}+1}\cdots\varphi_{u_{_{i}}-1}(x_{u_i})\in {\rm K}(\varphi_{v_i}) = {\rm Im}(\varphi_{v_{_i}+1}\cdots \varphi_{u_{_{i}}-1}\varphi_{u_i}).$
Thus $$\varphi_{v_{_i}+1}\cdots\varphi_{u_{_{i}}-1}(x_{u_i})= \varphi_{v_{_i}+1}\cdots \varphi_{u_{_{i}}-1}\varphi_{u_i}(x_{{u_{_i}}+1})$$ for some
$x_{{u_{_i}}+1}\in X_{{u_{_i}}+1}.$
If $\{v_{_i}+1, \cdots,  u_{_i}-1\}\cap \{v_1, \cdots, v_t\} = \emptyset$, then all the maps $\varphi_{v_{_i}+1}, \cdots,  \varphi_{u_{_i}-1}$ are monic,
and hence $x_{u_i}\in {\rm Im}(\varphi_{u_i})$. If
$\{v_{_i}+1, \cdots,  u_{_i}-1\}\cap \{v_1, \cdots, v_t\} = \{v_{i+1}, \cdots, v_{i+a}\} \ (a\ge 1)$. Then
$\varphi_{v_{_i}+1}\cdots \varphi_{u_{_i}-1} = \varphi_{v_{_i}+1}\cdots \varphi_{v_{_{i+1}}} \cdots \varphi_{v_{_{i+a}}} \cdots \varphi_{u_{_i}-1}$.
Taking off the monomorphisms we get
 $$\varphi_{v_{i+1}} \cdots \varphi_{v_{i+a}} \cdots \varphi_{u_{_i}-1}(x_{u_i}) = \varphi_{v_{i+1}} \cdots \varphi_{v_{i+a}} \cdots \varphi_{u_{_i}-1}\varphi_{u_{_i}}(x_{u_{_i}+1}).$$
Repeating this argument we finally get  $x_{u_i}\in {\rm Im}(\varphi_{u_i})$.

\vskip10pt

By construction $\psi_{u_i}: {\rm Ck}(\varphi_{v_{_{i}}+1}\cdots \varphi_{v_{_i}+m_{_i}})\longrightarrow X_{v_{_{i-1}}+1}$ for $i\in\{1, \cdots, t\}$.
Thus
\begin{align*}{\rm K}(\psi_{u_i}) &
\stackrel{(7.6)} = {\rm K}(\varphi_{v_{_{i-1}}+1}\cdots \varphi_{v_i})/{\rm Im}(\varphi_{v_{_{i}}+1}\cdots \varphi_{v_{_i}+m_{_i}})\\&
\stackrel{(7.3)} = {\rm K}(\varphi_{v_i})/{\rm Im}(\varphi_{v_{_{i}}+1}\cdots \varphi_{v_{_i}+m_{_i}})\\& \stackrel{(7.3)} = {\rm Im}(\varphi_{v_{_i}+1}\cdots \varphi_{u_i})/{\rm Im}(\varphi_{v_{_{i}}+1}\cdots \varphi_{u_i}\varphi_{v_{_i}+m_{_i}}).\end{align*}
We claim that it is isomorphic to $X_{v_{_i}+m_{_i}}/{\rm Im}(\varphi_{v_{_i}+m_{_i}}) = {\rm Ck}(\varphi_{v_{_i}+m_{_i}})\in \x$ via
$$g: \varphi_{v_{_i}+1}\cdots \varphi_{u_i}(x_{v_{_i}+m_{_i}}) + {\rm Im}(\varphi_{v_{_{i}}+1}\cdots \varphi_{u_i}\varphi_{v_{_i}+m_{_i}})\mapsto x_{v_{_i}+m_{_i}}+{\rm Im}(\varphi_{v_{_i}+m_{_i}}).$$
It suffices to see that $g$ is well-defined (clearly $g$ is monic and epic, if it is well-defined). This can be similarly shown as above.
We briefly sketch the process. Assume that $\varphi_{v_{_i}+1}\cdots \varphi_{u_i}(x_{v_{_i}+m_{_i}}) = \varphi_{v_{_{i}}+1}\cdots \varphi_{u_i}\varphi_{v_{_i}+m_{_i}}(x_{v_{_i}+m_{_i}+1})$ for some $x_{v_{_i}+m_{_i}+1}\in X_{{v_{_i}+m_{_i}+1}}$.
We need to show  $x_{v_{_i}+m_{_i}}\in {\rm Im}(\varphi_{v_{_i}+m_{_i}})$.
This is clearly true if $\{v_{_i}+1, \cdots,  u_{_i}\}\cap \{v_1, \cdots, v_t\} = \emptyset$. Suppose
$\{v_{_i}+1, \cdots,  u_{_i}\}\cap \{v_1, \cdots, v_t\} = \{v_{i+1}, \cdots, v_{i+a}\} \ (a\ge 1)$.
Then $$\varphi_{v_{i+1}} \cdots \varphi_{v_{i+a}} \cdots \varphi_{u_{_i}}(x_{v_{_i}+m_{_i}}) = \varphi_{v_{i+1}} \cdots \varphi_{v_{i+a}} \cdots \varphi_{u_{_i}}\varphi_{v_{_i}+m_{_i}}(x_{v_{_i}+m_{_i}+1}).$$
Thus
\begin{align*}&\varphi_{v_{_{i+1}}+1} \cdots \varphi_{v_{i+a}} \cdots \varphi_{u_{_i}}(x_{v_{_i}+m_{_i}})-\varphi_{v_{_{i+1}}+1} \cdots \varphi_{v_{i+a}} \cdots \varphi_{u_{_i}}\varphi_{v_{_i}+m_{_i}}(x_{v_{_i}+m_{_i}+1})
\\ & \subseteq {\rm K}(\varphi_{v_{i+1}}) \stackrel{(7.3)} = {\rm Im}(\varphi_{v_{_{i+1}}+1} \cdots \varphi_{u_{_{i+1}}})\subseteq {\rm Im}(\varphi_{v_{_{i+1}}+1} \cdots \varphi_{v_{i+a}} \cdots \varphi_{u_{_i}}\varphi_{v_{_i}+m_{_i}}).\end{align*}
So $$\varphi_{v_{_{i+1}}+1} \cdots \varphi_{v_{i+a}} \cdots \varphi_{u_{_i}}(x_{v_{_i}+m_{_i}})=\varphi_{v_{_{i+1}}+1} \cdots \varphi_{v_{i+a}} \cdots \varphi_{u_{_i}}\varphi_{v_{_i}+m_{_i}}(x'_{v_{_i}+m_{_i}+1})$$
some $x'_{v_{_i}+m_{_i}+1}\in X_{_{v_{_i}+m_{_i}+1}}.$ Repeating this argument we finally get  $x_{v_{_i}+m_{_i}}\in {\rm Im}(\varphi_{v_{_i}+m_{_i}})$.

\vskip10pt

If $i, i+1\in [u_j+1, u_{_{j+1}}-1]$  for some $j\in\{0, \cdots, t\}$, then by construction
$\psi_i: {\rm Ck}(\varphi_{v_{_j}+1}\cdots \varphi_{i}\varphi_{i+1})\twoheadrightarrow {\rm Ck}(\varphi_{v_{_j}+1}\cdots \varphi_{i})$.
Thus \begin{align*}{\rm K}(\psi_i) &= {\rm Im}(\varphi_{v_{_j}+1}\cdots \varphi_{i})/{\rm Im}(\varphi_{v_{_j}+1}\cdots \varphi_{i}\varphi_{i+1}).\end{align*}
By the similar argument we see that it is isomorphic to $X_{i+1}/{\rm Im}(\varphi_{i+1}) = {\rm Ck}(\varphi_{i+1})\in \x$ via
$$\varphi_{v_{_j}+1}\cdots \varphi_{i}(x_{i+1}) + {\rm Im}(\varphi_{v_{_j}+1}\cdots \varphi_{i}\varphi_{i+1})\mapsto x_{i+1}+{\rm Im}(\varphi_{i+1}).$$
This completes the proof of {\bf Step 3.}

\vskip5pt

{\bf Step 4.} \ For $Y = (Y_i, \psi_i)\in {\rm sepi}(Q, I, A)$,
define $GY = ((GY)_i, \varphi_i)$ as follows.

\vskip5pt

The $(v_{_j}+1)$-st branch $(GY)_{v_{_j}+1}$ of $GY$ is
$Y_{u_{_{j+1}}}$ for $0\le j\le t$. Thus, the first branch
of $GY$ is $Y_{u_{_1}}$, and the $(v_{_t}+1)$-st branch of $GY$ is
$Y_n$.

\vskip5pt

For $i\in [v_{_j}+2, v_{_{j+1}}]$ for $0\le j\le t$, the $i$-th
branch $(GY)_i$ of $GY$ is $\K(\psi_{i-1}\cdots
\psi_{u_{_{j+1}}-1})$, where $Y_{u_{_{j+1}}}
\stackrel{\psi_{u_{_{j+1}}-1}}\twoheadrightarrow Y_{u_{_{j+1}}-1}
\twoheadrightarrow \cdots \stackrel {\psi_{i-1}}\twoheadrightarrow
Y_{i-1}$.  Thus
\begin{align*}
GY: = (\K(\psi_{n-1})\stackrel{\varphi_{n-1}}\hookrightarrow \cdots
\hookrightarrow  \K(\psi_{v_t+1}\cdots \psi_{n-1})\hookrightarrow
&Y_n \stackrel{\varphi_{v_t}}\longrightarrow\K(\psi_{v_t-1}\cdots
\psi_{u_{_{t}}-1})
\\ \hookrightarrow \cdots\hookrightarrow\K(\psi_{v_{_{t-1}}+1}\cdots \psi_{u_{_t}-1})\hookrightarrow &Y_{u_t}
\stackrel{\varphi_{v_{
_{\small{t-1}}}}}\longrightarrow\K(\psi_{v_{_{t-1}}-1}\cdots
\psi_{u_{_{t-1}}-1})
\\ \hookrightarrow \cdots\hookrightarrow\K(\psi_{v_1+1}\cdots \psi_{u_{_2}-1})\hookrightarrow &Y_{u_2}
\stackrel{\varphi_{v_1}}\longrightarrow\K(\psi_{v_1-1}\cdots
\psi_{u_{_1}-1})
\\ \hookrightarrow \cdots \hookrightarrow\K(\psi_1\cdots
\psi_{u_{_1}-1}) \stackrel{\varphi_1}\hookrightarrow &Y_{u_1}) \ \ \
\ \ \ \ \  \ \ \ \  (\mbox{\small{no this line if}} \ v_1 = 1)  \ \
\ \ \ \ \ \  \ \ \ (7.7)\end{align*} where each $\varphi_i$ is the
natural embedding for $i\notin \{v_1, \cdots, v_t\}$, and for $1\le
j\le t$
$$\varphi_{v_j}: = \widetilde{\psi_{u_j}\cdots \psi_{u_{_{j+1}}-1}}:
Y_{u_{_{j+1}}}\longrightarrow \K(\psi_{v_j-1}\cdots
\psi_{u_{_j}-1})\eqno(7.8)$$ is induced by $\psi_{u_j}\cdots
\psi_{u_{_{j+1}}-1}: Y_{u_{_{j+1}}}\longrightarrow Y_{u_j}.$

\vskip5pt

{\bf Remark:} \ If $v_{_{j}} -1 < v_{_{j-1}}+1$ for some $1\le j
\le t$, then $v_{j}=  v_{_{j-1}}+1$. Thus the $(v_{_{j}}+1)$-th
branch $Y_{u_{_{j+1}}}$ neighbors with the $(v_{_j} =
v_{_{j-1}}+1)$-th branch $Y_{u_j}$, and the corresponding part of
$(7.7)$ is read as
$$Y_{u_{_{j+1}}}\stackrel{\varphi_{v_{_j}}:=\psi_{u_{_{j+1}}-1}\cdots
\psi_{u_{_j}}}\longrightarrow Y_{u_j}.$$ In particular, if $v_{_1} =
1$, then the last line of $(7.7)$ disappears, and
$\K(\psi_{v_{_1}-1}\cdots \psi_{u_{_1}-1})$  in the first relation
does not make sense, and it is $Y_{u_{_1}} = Y_{m_{_1}}$.

If $v_{_{j}} -1 = v_{_{j-1}}+1$ for some $1\le j \le t$, then there
is only one $A$-modules between $Y_{u_{_{j+1}}}$  and $Y_{u_j}$. In
particular, if $v_1= 2$, then the last line of $(7.7)$ contains only
one module $Y_{m_{_1}+1}$.

\vskip5pt

{\bf Step 5.} \ For $Y = (Y_i, \psi_i)\in {\rm sepi}(Q, I, A)$, $GY
= ((GY)_i, \varphi_i)\in {\rm smon}(Q, I, A)$. Thus, we get a functor
$G: {\rm sepi}(Q, I, A)\longrightarrow {\rm smon}(Q, I, A)$.

\vskip5pt

By construction $\varphi_i$ is monic for $i\notin \{v_1, \cdots,
v_t\}$.  We need to show $\K(\varphi_{v_{_j}}) = {\rm
Im}(\varphi_{v_{_j}+1}\cdots \varphi_{u_{_j}})$ for $1\le j\le t$.
By $(7.8)$, $\K(\varphi_{v_{_j}}) = \K(\psi_{u_{_j}}\cdots
\psi_{u_{_{j+1}}-1})$, and it is ${\rm Im}(\varphi_{v_{_j}+1}\cdots
\varphi_{u_{_j}})$,  by taking $i=u_j$ in {\rm Sublemma 2} below,

\vskip5pt

The following sublemma is the same with {\rm Sublemma 1}, but only
after Theorem \ref{RSS1} is proved.

\vskip5pt

{\rm \bf Sublemma 2.} \ Let $Y = (Y_i, \psi_i)\in {\rm sepi}(Q, I,
A)$, \ and $i\in [v_{_j}+1, u_{_{j+1}}-1]$ for some $0\le j\le t$.
Then $GY = ((GY)_i, \varphi_i)$ has the property ${\rm
Im}(\varphi_{v_{_j}+1}\cdots \varphi_i) =
\K(\psi_i\cdots\psi_{u_{_{j+1}}-1}).$

\vskip5pt \noindent{\bf Proof of Sublemma 2.} \ First, assume
$[v_{_j}+1, i]\cap \{v_1, \cdots, v_t\} = \emptyset.$  Then by {\bf
Step 3} all the maps $\varphi_{v_{_j}+1}, \cdots, \varphi_{i}$ are
embeddings. Since $v_{_{j+1}}\notin [v_{_j}+1, i]$,  we have $i+1\in
[v_{_j}+2, v_{_{j+1}}]$, and hence by {\bf Step 3}, $(GY)_{i+1} =
\K(\psi_{i}\cdots\psi_{u_{_{j+1}-1}})$. Also $(GY)_{v_{_j}+1} =
Y_{u_{_{j+1}}}$. Thus
\begin{align*}&{\rm
Im}(\varphi_{v_{_j}+1}\cdots \varphi_{i}) = {\rm
Im}(\K(\psi_{i}\cdots\psi_{u_{_{j+1}}-1})\stackrel{\varphi_{i}}\hookrightarrow
\cdots \stackrel{\varphi_{v_{_j}+1}}\hookrightarrow Y_{u_{_{j+1}}})
\cong \K(\psi_{i}\cdots\psi_{u_{_{j+1}}-1}).\end{align*}

It remains to consider the case that $[v_{_j}+1, i]\cap \{v_{_1},
\cdots, v_{_t}\} = \{v_{_{j+1}}, \cdots, v_{_{j+a}}\} \ (a\ge 1).$
Since all the maps in $\{\varphi_{v_{_j}+1}, \cdots
\varphi_{v_{j+1}},  \cdots \varphi_{v_{j+a}}, \cdots \varphi_i\}$
are embeddings except those maps $\varphi_{v_{j+b}}$ for $1\le b\le
a$, we have
\begin{align*}&{\rm
Im}(\varphi_{v_{_j}+1}\cdots \varphi_{i}) = {\rm
Im}(\varphi_{v_{_j}+1}\cdots \varphi_{v_{j+1}} \cdots
\varphi_{v_{j+a}}\cdots \varphi_i) = \varphi_{v_{j+1}} \cdots
\varphi_{v_{j+a}}((GY)_{i+1})
\end{align*}
where $\varphi_{v_{j+1}} \cdots \varphi_{v_{j+a}}$ is the
composition of only those maps $\varphi_{v_{j+b}}$ for $1\le b\le
a$, namely, all the kernel embeddings can be taken off. Note that
$i+1\le v_{j+a+1}$ (otherwise $v_{j+a+1}$ is also in $[v_{_j}+1,
i]$). While $i+1\ge v_{_{j+a}}+1.$ Thus either $i+1\in [v_{j+a}+2,
v_{j+a+1}]$ or $i+1 = v_{_{j+a}}+1.$

\vskip5pt

If $i+1\in [v_{j+a}+2, v_{j+a+1}]$, then $(GY)_{i+1}
=\K(\psi_{i}\cdots \psi_{u_{_{j+a+1}}-1})$, and by $(7.8)$ we have
\begin{align*}&{\rm
Im}(\varphi_{v_{_j}+1}\cdots \varphi_i) = (\psi_{u_{_{j+1}}}\cdots
\psi_{u_{_{j+2}}-1}) \cdots (\psi_{u_{_{j+a}}}\cdots
\psi_{u_{_{j+a+1}}-1})(\K(\psi_{i}\cdots \psi_{u_{_{j+a+1}}-1})).
\end{align*}
We first compute $(\psi_{u_{_{j+a}}}\cdots
\psi_{u_{_{j+a+1}}-1})(\K(\psi_{i}\cdots \psi_{u_{_{j+a+1}}-1})).$
Note that  for $A$-maps $f: U \longrightarrow V$ and $g:
V\longrightarrow W$, if $f$ is epic, then $f(\K(gf)) = \K(g)$. Since
$\psi_{u_{_{j+a}}+1}\cdots \psi_{u_{_{j+a+1}}-1}$ is epic, we have
\begin{align*}&(\psi_{u_{_{j+a}}}\cdots
\psi_{u_{_{j+a+1}}-1})(\K(\psi_{i}\cdots \psi_{u_{_{j+a+1}}-1}))
\\ & = \psi_{u_{_{j+a}}}(\psi_{u_{_{j+a}}+1}\cdots \psi_{u_{_{j+a+1}}-1})(\K((\psi_{i}\cdots\psi_{u_{_{j+a}}-1}\psi_{u_{_{j+a}}})(\psi_{u_{_{j+a}}+1}\cdots \psi_{u_{_{j+a+1}}-1}))
\\& = \psi_{u_{_{j+a}}}(\K(\psi_{i}\cdots\psi_{u_{_{j+a}}-1}\psi_{u_{_{j+a}}})).
\end{align*}
While $\psi_{u_{_{j+a}}}$ is {\bf not} epic. But thanks to $(7.4)$,
we still claim
$$\psi_{u_{_{j+a}}}(\K(\psi_{i}\cdots\psi_{u_{_{j+a}}-1}\psi_{u_{_{j+a}}})) =
\K(\psi_{i}\cdots\psi_{u_{_{j+a}}-1}).$$ In fact, it is clear that
$\psi_{u_{_{j+a}}}(\K(\psi_{i}\cdots\psi_{u_{_{j+a}}-1}\psi_{u_{_{j+a}}}))
\subseteq\K(\psi_{i}\cdots\psi_{u_{_{j+a}}-1}).$ For $y\in
\K(\psi_{i}\cdots\psi_{u_{_{j+a}}-1})$, we have
$(\psi_{v_{_{j+a}}}\cdots
\psi_{i-1})(\psi_{i}\cdots\psi_{u_{_{j+a}}-1})(y) = 0.$ So by
$(7.4)$, $y\in \K(\psi_{v_{_{j+a}}}\cdots \psi_{u_{_{j+a}}-1}) =
{\rm Im}(\psi_{u_{_{j+a}}})$. From this one easily sees that $y\in
\psi_{u_{_{j+a}}}(\K(\psi_{i}\cdots\psi_{u_{_{j+a}}-1}\psi_{u_{_{j+a}}}))$.
This proves the claim, and hence $$(\psi_{u_{_{j+a}}}\cdots
\psi_{u_{_{j+a+1}}-1})(\K(\psi_{i}\cdots \psi_{u_{_{j+a+1}}-1}))
=\K(\psi_{i}\cdots\psi_{u_{_{j+a}}-1}).$$ It follows that
\begin{align*}&{\rm
Im}(\varphi_{v_{_j}+1}\cdots \varphi_i) = (\psi_{u_{_{j+1}}}\cdots
\psi_{u_{_{j+2}}-1}) \cdots (\psi_{u_{_{j+a-1}}}\cdots
\psi_{u_{_{j+a}}-1})(\K(\psi_{i}\cdots\psi_{u_{_{j+a}}-1})).
\end{align*}
Continuing this process we get ${\rm Im}(\varphi_{v_{_j}+1}\cdots
\varphi_{v_{j+1}} \cdots \varphi_{v_{j+a}}\cdots \varphi_i) =
\K(\psi_i\cdots\psi_{u_{_{j+1}}-1}).$

\vskip5pt

If $i+1 = v_{_{j+a}}+1,$ then $(GY)_{i+1} =Y_{u_{_j}+a+1}$, and we
are computing ${\rm Im}(\varphi_{v_{_j}+1}\cdots
\varphi_{v_{j+a}})$. By the same argument, in particular, by using
${\rm Im}(\psi_{u_{_{j+a}}}) = \K(\psi_{v_{_{j+a}}}\cdots
\psi_{u_{_{j+a}}-1})$, we get ${\rm Im}(\varphi_{v_{_j}+1}\cdots
\varphi_i) = \K(\psi_i\cdots\psi_{u_{_{j+1}}-1}).$ This completes
the proof of {\rm Sublemma 2}.

\vskip5pt

{\bf Step 6.} \ If $Y = (Y_i, \psi_i)\in {\rm sepi}(Q, I, \x)$, then $GY
= ((GY)_i, \varphi_i)\in {\rm smon}(Q, I, \x)$.

\vskip5pt

We need to prove $(GY)_n\in \x$ and
${\rm Ck}(\varphi_i)\in \x$ for $1\le i\le n-1$. By assumption we already have $Y_1\in\x$ and  ${\rm K}(\psi_i)\in \x$ for $1\le i\le n-1$.

\vskip5pt

By construction $(GY)_n = {\rm K}(\psi_{n-1})\in \x$.

\vskip5pt

By construction $\varphi_{v_i}: Y_{u_{_{i+1}}}\longrightarrow {\rm K}(\psi_{v_i-1}\cdots \psi_{u_i-1})$ for $i\in\{1, \cdots, t\}$.
Thus
\begin{align*}{\rm Ck}(\varphi_{v_i}) &
= {\rm K}(\psi_{v_i-1}\cdots \psi_{u_i-1})/{\rm Im}(\varphi_{v_i}) \\ &\stackrel{(7.8)} = {\rm K}(\psi_{v_i-1}\cdots \psi_{u_i-1})/{\rm Im}(\psi_{u_i}\cdots \psi_{u_{_{i+1}}-1})\\&
\stackrel{(7.4)} = {\rm K}(\psi_{v_i-1}\cdots \psi_{u_i-1})/{\rm Im}(\psi_{u_i})\\& \stackrel{(7.4)} = {\rm K}(\psi_{v_i-1}\psi_{v_i}\cdots \psi_{u_i-1})/{\rm Ker}(\psi_{v_i}\cdots\psi_{u_{_i}-1})
\\& \cong \psi_{v_i}\cdots\psi_{u_{_i}-1}({\rm K}(\psi_{v_i-1}\psi_{v_i}\cdots \psi_{u_i-1})
.\end{align*}
It is clearly contained in ${\rm K}(\psi_{v_i-1})$. We claim that it is exactly ${\rm K}(\psi_{v_i-1})\in \x$.
In fact, if $\{v_{_i}, \cdots,  u_{_i}-1\}\cap \{u_1, \cdots, u_t\} = \emptyset$, then all the maps $\psi_{v_{_i}}, \cdots,  \psi_{u_{_i}-1}$ are surjective,
and hence the claim holds. If $\{v_{_i}, \cdots,  u_{_i}-1\}\cap \{u_1, \cdots, u_t\} = \{u_{i-1}, \cdots, u_{i-a}\} \ (a\ge 1)$, then
$\psi_{v_{_i}}\cdots \psi_{u_{_i}-1} = \psi_{v_{_i}} \cdots \psi_{u_{_{i-a}}} \cdots \psi_{u_{_{i-1}}}\cdots \psi_{u_{_i}-1}$.
For each $y_{v_{_i}}\in {\rm K}(\psi_{v_{_i}-1})$, since
all the maps $\psi_{v_{_i}}, \cdots,  \psi_{u_{_{i-a}}-1}$ are surjective, it follows that
$y_{v_{_i}} = \psi_{v_{_i}} \cdots  \psi_{u_{_{i-a}}-1}(y_{u_{_{i-a}}})$ for some $y_{u_{_{i-a}}}\in Y_{u_{_{i-a}}}.$ Thus
\begin{align*}y_{u_{_{i-a}}}\in {\rm K}(\psi_{v_{_i}-1}\psi_{v_{_i}} \cdots  \psi_{u_{_{i-a}}-1})&\subseteq {\rm K}(\psi_{v_{_{i-a}}}\cdots \psi_{v_{_i}-1}\psi_{v_{_i}} \cdots  \psi_{u_{_{i-a}}-1})
\\& \stackrel{(7.4)}= {\rm Im}(\psi_{u_{i-a}}).\end{align*}
So $y_{v_{_i}} = \psi_{v_{_i}} \cdots  \psi_{u_{_{i-a}}-1}\psi_{u_{i-a}}(y_{u_{_{i-a}}+1})$ for some $y_{u_{_{i-a}}+1}\in Y_{u_{_{i-a}}+1}.$
Repeating this argument we finally get $y_{v_{_i}} = \psi_{v_{_i}} \cdots  \psi_{u_{i}-1}(y_{u_i})$ for some $y_{u_i}\in Y_{u_i},$
and moreover $y_{u_i}\in {\rm K}(\psi_{v_{_i}-1}\psi_{v_{_i}}\cdots \psi_{u_{_i}-1}).$
This proves the claim.

\vskip5pt

By construction $\varphi_{v_{_i}+1}: {\rm K}(\psi_{v_{_i}+1}\cdots \psi_{u_{_{i+1}}-1})\hookrightarrow Y_{u_{_{i+1}}}$ for $i\in\{0, \cdots, t\}$.
Thus
\begin{align*}{\rm Ck}(\varphi_{v_{_i}+1}) &
= Y_{u_{_{i+1}}}/{\rm K}(\psi_{v_i+1}\cdots \psi_{u_{_{i+1}}-1})\cong {\rm Im}(\psi_{v_{_i}+1}\cdots \psi_{u_i}\cdots \psi_{u_{_{i+1}}-1})\\&
\stackrel{(7.4)} ={\rm Im}(\psi_{v_{_i}+1}\cdots \psi_{u_i})
\stackrel{(7.4)} = \psi_{v_{_i}+1}\cdots \psi_{u_{_i}-1}({\rm K}(\psi_{v_i}\psi_{v_{_i}+1}\cdots \psi_{u_{_i}-1}).\end{align*}
By the same argument as above we see that it is exactly ${\rm K}(\psi_{v_i})\in \x$.

\vskip5pt

If $i, i+1\in [v_{_{j}}+2, v_{_{j+1}}]$  for some $j\in\{0, \cdots, t\}$, then by construction $\varphi_i: {\rm K}(\psi_i\cdots \psi_{u_{_{j+1}}-1})\hookrightarrow {\rm K}(\psi_{i-1}\psi_i\cdots \psi_{u_{_{j+1}}-1})$.
Thus \begin{align*}{\rm Ck}(\varphi_i) &= {\rm K}(\psi_{i-1}\psi_i\cdots \psi_{u_{_{j+1}}-1})/{\rm K}(\psi_i\cdots \psi_{u_{_{j+1}}-1})\\ &\cong
\psi_i\cdots \psi_{u_{_{j+1}}-1}({\rm K}(\psi_{i-1}\psi_i\cdots \psi_{u_{_{j+1}}-1})).\end{align*}
Again by the same argument as above we see that it is exactly ${\rm K}(\psi_{i-1})\in \x$.

\vskip5pt

This completes the proof of {\bf Step 6.}

\vskip5pt

{\bf Step 7.} \ $G \circ F \cong {\rm Id}_{{\rm smon}(Q, I, \x)}.$

\vskip5pt

We first prove $(GFX)_i\cong X_i$, for $X=(X_i, \varphi_i)\in {\rm
smon}(Q, I, \x)$ and  for $$i\in \bigcup\limits_{0\le j\le
t}(\{v_{_{j}}+1\}\cup [v_{_j}+2, v_{_{j+1}}]).$$

We already have  $FX = ((FX)_i, \psi_i)\in {\rm sepi}(Q, I, \x)$. By
{\bf Step 4} and {\bf Step 1},  $(GFX)_{v_{_{j}}+1} =
(FX)_{u_{_{j+1}}} = X_{v_{_j}+1}.$  Assume that $i\in [v_{_j}+2, v_{_{j+1}}]$ for some $0\le j\le t$. By
{\bf Step 4},  $(GFX)_i = \K(\psi_{i-1}\cdots \psi_{u_{_{j+1}}-1})$.
Since $i-1\in [v_j+1, u_{_{j+1}}-1]$, we can apply {\rm Sublemma 1}
(note that in {\rm Sublemma 1} $1\le j\le t+1$) to get $(GFX)_i =
\K(\psi_{i-1}\cdots \psi_{u_{_{j+1}}-1}) = {\rm
Im}(\varphi_{v_{_j}+1}\cdots \cdots \varphi_{_{i-1}}).$
 Note that $[v_{_{j}}+1, i-1]\cap \{v_1, \cdots, v_t\} = \emptyset$
(otherwise $v_{_{j+1}}\in [v_{_{j}}+1, i-1];$ but by assumption
$i\le v_{_{j+1}}$). So $\{\varphi_{v_{_j}+1}, \cdots,
\varphi_{_{i-1}}\}\cap \{\varphi_{v_1}, \cdots, \varphi_{v_t}\}
=\emptyset$, and hence all the maps $\varphi_{v_{_j}+1}, \cdots,
\varphi_{_{i-1}}$ are embeddings. Thus $(GFX)_i={\rm
Im}(\varphi_{v_{_j}+1}\cdots \varphi_{_{i-1}}) \cong X_i$.

\vskip5pt
It is routine to verify that the diagram
\[\xymatrix @C=0.6cm {(GFX)_{i+1} \ar[r]\ar[d]^-{\cong} & (GFX)_i\ar[d]^-{\cong}\\
X_{i+1} \ar[r]^-{\varphi_i} & X_i}\]
commutes for $1\le i\le n$.

\vskip5pt

{\bf Step 8.} \ $F\circ G \cong {\rm Id}_{{\rm sepi}(Q, I, \x)}.$

\vskip5pt

We first prove $(FGY)_i\cong Y_i$, for $Y=(Y_i, \psi_i)\in {\rm sepi}(Q,
I, \x)$ and for
$$i\in \bigcup\limits_{0\le j\le t}([u_{_j}+1,
u_{_{j+1}}-1]\cup\{u_{_{j+1}}\}).$$

We already have $GY = ((GY)_i, \varphi_i)\in {\rm smon}(Q, I, \x)$.
By {\bf Step 1} and {\bf Step 4}, $(FGY)_{u_{_{j+1}}} = (GY)_{v_j+1}
= Y_{u_{_{j+1}}}.$ Assume that $i\in [u_{_j}+1, u_{_{j+1}}-1]$ for some $0\le j\le t$.
By {\bf Step 1},  $(FGY)_i = \co(\varphi_{v_{_j}+1}\cdots \varphi_i)
= Y_{u_{_{j+1}}}/{\rm Im}(\varphi_{v_{_j}+1}\cdots \varphi_i).$ By
{\rm Sublemma 2}, ${\rm Im}(\varphi_{v_{_j}+1}\cdots \varphi_i) =
\K(\psi_i\cdots\psi_{u_{_{j+1}}-1})$. It follows that
$(FGY)_i=Y_{u_{_{j+1}}}/\K(\psi_i\cdots\psi_{u_{_{j+1}}-1})\cong
{\rm Im}(\psi_i\cdots\psi_{u_{_{j+1}-1}})$. Note that $[i,
u_{_{j+1}}-1]\cap \{u_1, \cdots, u_t\} = \emptyset$ (otherwise
$u_j\in [i, u_{_{j+1}}-1];$ but by assumption $i\ge u_{_j}+1 >
u_j$). So all the maps $\psi_i, \cdots, \psi_{u_{_{j+1}-1}}$ are
epimorphisms. Thus $(FGY)_i\cong {\rm
Im}(\psi_i\cdots\psi_{u_{_{j+1}-1}}) = Y_i$.

\vskip5pt
It is routine to verify that the diagram
\[\xymatrix @C=0.6cm {(FGY)_{i+1} \ar[r]\ar[d]^-{\cong} & (FGY)_i\ar[d]^-{\cong}\\
Y_{i+1} \ar[r]^-{\varphi_i} & Y_i}\]
commutes for $1\le i\le n$.

\vskip5pt

{\bf Step 9.} \ There is a functorial left $\m$-module isomorphism
$F (M\otimes P(i)) \cong M\otimes I(i)$ for $i\in Q_0$ and $M\in
\x$.

\vskip5pt

This only needs a careful verification. Assume that $j$ is the
maximal nonnegative integer such that $i\ge u_{_j}+1$  \ ($0\le j\le
t$). Then $P(i) = (kQ/I)e_i \in {\rm rep}(Q, I, k)$ is
$$(0\rightarrow \cdots \rightarrow 0 \rightarrow e_i(kQ/I)e_i = k\stackrel {{\rm Id}}\longrightarrow k \rightarrow\cdots \stackrel {{\rm Id}}\longrightarrow
k = e_{v_{_j}+1}(kQ/I)e_i \rightarrow 0 \rightarrow \cdots
\rightarrow 0)$$and hence $M\otimes P(i)\in {\rm smon}(Q, I, \x)$ is
$$(0\rightarrow \cdots \rightarrow 0 \rightarrow (M\otimes P(i))_i = M\stackrel {{\rm Id}}\longrightarrow M \rightarrow\cdots \stackrel {{\rm Id}}\longrightarrow
M = (M\otimes P(i))_{v_{_j}+1}\rightarrow 0 \rightarrow \cdots
\rightarrow 0).$$ By $(7.3)$ we have $F (M\otimes P(i)) \cong
M\otimes I(i)$.

\vskip5pt

This completes the proof of Theorem \ref{RSS1}. $\s$

\subsection{} We include an example. Let $Q = 6
\stackrel{\alpha_5}\longrightarrow 5
\stackrel{\alpha_4}\longrightarrow 4
\stackrel{\alpha_3}\longrightarrow 3
\stackrel{\alpha_2}\longrightarrow
2\stackrel{\alpha_1}\longrightarrow 1$ with $I=\langle
\alpha_1\alpha_2\alpha_3, \ \alpha_3\alpha_4\rangle$, $A$ a
finite-dimensional algebra, and $\x$ an arbitrary additive subcategory of $A$-mod. An $\m$-module $X =
(X_6\stackrel{\varphi_5} \hookrightarrow X_5
\stackrel{\varphi_4}\hookrightarrow
X_4\stackrel{\varphi_3}\longrightarrow X_3
\stackrel{\varphi_2}\hookrightarrow
X_2\stackrel{\varphi_1}\longrightarrow X_1)\in {\rm smon}(Q, I, A)$
if and only if
$$\varphi_i \ \mbox{is monic if} \ i\notin \{1, 3\}, \ \ \Ker(\varphi_1) = {\rm Im}(\varphi_2\varphi_3), \ \
\Ker(\varphi_3) = {\rm Im}(\varphi_4);$$ and $Y =
(Y_6\stackrel{\psi_5} \twoheadrightarrow Y_5
\stackrel{\psi_4}\longrightarrow Y_4\stackrel{\psi_3}\longrightarrow
Y_3 \stackrel{\psi_2}\twoheadrightarrow
Y_2\stackrel{\psi_1}\twoheadrightarrow Y_1)\in {\rm sepi}(Q, I, \x)$
if and only if
$$\psi_i \ \mbox{is epic if} \ i\notin \{3, 4\}, \ \ {\rm Im}(\psi_3) = \Ker(\psi_1\psi_2), \ \
{\rm Im}(\psi_4) = \Ker(\psi_3).$$ Then an {\rm RSS} equivalence $F\colon
{\rm smon}(Q, I, \x) \cong {\rm sepi}(Q, I, \x)$ sends $X\in {\rm
smon}(Q, I, \x)$ to
\begin{align*} (FX: = X_4\stackrel{\psi_5}\twoheadrightarrow \Cok(\varphi_4\varphi_5)\stackrel{\psi_4:=\widetilde{\varphi_2\varphi_3}}\longrightarrow
X_2\stackrel{\psi_3:=\varphi_1}\longrightarrow X_1
\stackrel{\psi_2}\twoheadrightarrow
\Cok(\varphi_1\varphi_2)\stackrel{\psi_1}\twoheadrightarrow
\Cok(\varphi_1))\in {\rm sepi}(Q, I, \x).
\end{align*}
A quasi-inverse $G: {\rm sepi}(Q, I, \x) \longrightarrow {\rm
smon}(Q, I, \x)$ of $F$ sends $Y\in {\rm sepi}(Q, I, \x)$ to
\begin{align*} (GY:=\Ker(\psi_5)\stackrel{\varphi_5}\hookrightarrow \Ker(\psi_4\psi_5)\stackrel{\varphi_4}\hookrightarrow Y_6\stackrel{\varphi_3:=\widetilde{\psi_4\psi_5}}\longrightarrow \Ker(\psi_2\psi_3)
\stackrel{\varphi_2} \hookrightarrow
Y_4\stackrel{\varphi_1:=\psi_3}\longrightarrow Y_3)\in {\rm smon}(Q,
I, \x).
\end{align*}
In particular,  for $M\in \x$ we have
\begin{align*}& F(M\otimes P(1))\cong F(0\rightarrow 0\rightarrow 0\rightarrow
0\rightarrow 0 \rightarrow M) \cong(0\rightarrow 0\rightarrow
0\rightarrow M\stackrel{=}\rightarrow M\stackrel{=}\rightarrow
M)\cong M\otimes I(1)\\ & F(M\otimes P(2)) \cong F(0\rightarrow
0\rightarrow 0\rightarrow 0\rightarrow M\stackrel{=}\rightarrow M) \cong
(0\rightarrow 0\rightarrow M\stackrel{=}\rightarrow
M\stackrel{=}\rightarrow M\rightarrow 0) \cong M\otimes I(2)
\\ & F(M\otimes P(3))\cong F(0\rightarrow 0\rightarrow 0\rightarrow
M\stackrel{=}\rightarrow M \stackrel{=}\rightarrow M)  \cong
(0\rightarrow 0\rightarrow M\stackrel{=}\rightarrow M\rightarrow
0\rightarrow 0)\cong M\otimes I(3)\\ & F(M\otimes P(4))\cong
F(0\rightarrow 0\rightarrow M\stackrel{=}\rightarrow M\rightarrow
M\rightarrow 0) \cong (M\stackrel{=}\rightarrow M\stackrel{=}\rightarrow
M\rightarrow 0\rightarrow 0\rightarrow 0) \cong M\otimes I(4)\\ &
F(M\otimes P(5))\cong F(0\rightarrow M\stackrel{=} \rightarrow
M\rightarrow 0\rightarrow 0 \rightarrow 0)  \cong
(M\stackrel{=}\rightarrow M\rightarrow 0\rightarrow 0\rightarrow
0\rightarrow 0)\cong M\otimes I(5)\\ & F(M\otimes P(6))\cong
F(M\stackrel{=}\rightarrow M\rightarrow M\rightarrow 0\rightarrow
0\rightarrow 0) \cong (M\rightarrow 0\rightarrow 0\rightarrow
0\rightarrow 0\rightarrow 0) \cong M\otimes I(6).
\end{align*}

\subsection{} \ For finite quivers $Q$ and $Q'$ (not necessarily acyclic), let $Q\otimes Q'$
be the quiver with
$$(Q\otimes Q')_0 = Q_0\times Q_0', \ \  \mbox{and} \ \ (Q\otimes Q')_1 = (Q_1\times
Q_0') \bigcup (Q_0\times Q_1').$$ Let $A = kQ/I$ and $B = kQ'/I'$ be finite-dimensional $k$-algebras,
where $Q$ and $Q'$ are finite quivers (not necessarily acyclic), $I$
and $I'$ are admissible ideals of $kQ$ and $kQ'$, respectively. Then
$$A\otimes_k B \cong k(Q\otimes Q')/(I\Box I'),$$
where $I\Box I'$ is the ideal of $k(Q\otimes Q')$ generated by
$(I\times Q_0')\bigcup (Q_0\times I')$ and commutative relations:
$$(\alpha, t')(i, \beta')-(j, \beta')(\alpha, s')$$
where $\alpha: i\longrightarrow j$ runs over $Q_1$, and $\beta':
s'\longrightarrow t'$ runs over $Q_1'$. See e.g. [L].

\subsection{} \ We include another {\rm RSS} equivalence $F: {\rm smon}(Q, I, \x)\cong {\rm sepi}(Q, I, \x)$ for arbitrary additive subcategory $\x$ (not necessarily extension-closed).
It is also constructed combinatorially, and is easily operable. Let $A$ be the path algebra of the quiver $b\longrightarrow a$, and $\m: = A\otimes kQ/I$, where $Q$ is the
quiver $$
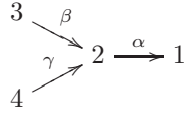
\xymatrix@R=0.1cm@C=0.2cm{3\ar[drr]^-\beta &&&& \\
& &2\ar[rr]^-\alpha & & 1\\ 4\ar[urr]^-\gamma&&&&}$$ and $I =
\langle \alpha\beta\rangle$. The
indecomposable $A$-modules are denoted by $P(a) =
\begin{smallmatrix}1\\0\end{smallmatrix} = S(a), \ \
P(b) = \begin{smallmatrix}1\\1\end{smallmatrix} = I(a)$, and  \
$I(b) =
\begin{smallmatrix}0\\1\end{smallmatrix} = S(b).$ Then $\m$
is the algebra given by the quiver
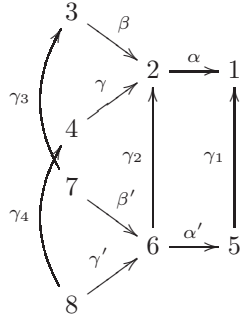
\[\xymatrix @R=0.3cm@C=0.2cm{3\ar[drr]^-\beta &&&& \\
& &2\ar[rr]^-\alpha & & 1\\ 4\ar[urr]^-\gamma&&&&\\
7\ar@/^1pc/[uuu]^-{\gamma_3}\ar[drr]^-{\beta'} &&&& \\
& &6\ar[uuu]^-{\gamma_2}\ar[rr]^-{\alpha'} & &
5\ar[uuu]^-{\gamma_1}\\
8\ar@/^1pc/[uuu]^-{\gamma_4}\ar[urr]^-{\gamma'}&&&&}\]  with
relations \ $\alpha\beta, \ \ \alpha'\beta', \ \ \alpha \gamma_2 -
\gamma_1\alpha', \ \ \  \beta \gamma_3 - \gamma_2\beta', \ \
\gamma\gamma_4 - \gamma_2\gamma'.$ A $\m$-module $X$ is
$$\xymatrix@R=0.1cm@C=0.2cm{X_3= {\begin{smallmatrix}a_3\\ b_3\end{smallmatrix}}\ar[drr]^-{X_\beta} &&&& \\
& &X_2= {\begin{smallmatrix}a_2\\ b_2\end{smallmatrix}}\ar[rr]^-{X_\alpha} & & X_1= {\begin{smallmatrix}a_1\\ b_1\end{smallmatrix}}\\
X_4= {\begin{smallmatrix}a_4\\
b_4\end{smallmatrix}}\ar[urr]^-{X_\gamma}&&&&}$$ where each $X_i$ is
an $A$-module, $X_\alpha$, $X_\beta$ and  $X_\gamma$ are $A$-maps
with $X_\alpha X_\beta = 0$. It will be
written as $\begin{smallmatrix}a_3&&&\\
a_4&a_2 &a_1\\ b_3 &b_2 &b_1 \\ b_4 & & & \end{smallmatrix}$ ({\bf not} as $\begin{smallmatrix}a_3&&&\\
b_3&a_2 &a_1\\ a_4 &b_2 &b_1 \\ b_4 & & & \end{smallmatrix}$). Then
$X\in {\rm smon}(Q, I, A)$ if and only if $X_\beta$ and $X_\gamma$ are monic, \ ${\rm Im}X_\beta + {\rm
Im}X_\gamma  = {\rm Im}X_\beta \oplus {\rm Im}X_\gamma$, and ${\rm
Ker}X_\alpha = {\rm Im}X_\beta.$ For examples,  the
indecomposable $\m$-module \
$\begin{smallmatrix}1&&&\\
1&2 &1\\ 0 &1 &1 \\ 0 & & & \end{smallmatrix}$  refers to
$$\xymatrix@R=0.1cm@C=0.2cm{{\begin{smallmatrix}1\\ 0\end{smallmatrix}}\ar[drr]^-{\binom{{\rm Id}}{0}} &&&& &\\
& &{\begin{smallmatrix}1\\ 0\end{smallmatrix}} \oplus
{\begin{smallmatrix}1\\
1\end{smallmatrix}}\ar[rr]^-{\left(\begin{smallmatrix}\sigma&{\rm
Id}\end{smallmatrix}\right)} & & {\begin{smallmatrix}1\\
1\end{smallmatrix}} & \in {\rm smon}(Q, I, A).\\
{\begin{smallmatrix}1\\
0\end{smallmatrix}}\ar[urr]^-{\binom{0}{\sigma}}&&&&&}$$ and the
indecomposable $\m$-module \
$\begin{smallmatrix}1&&&\\
1&1 &0\\ 1 &2 &1 \\ 1 & & & \end{smallmatrix}$ \  refers to
$$\xymatrix@R=0.1cm@C=0.2cm{{\begin{smallmatrix}1\\ 1\end{smallmatrix}}\ar[drr]^-{\binom{{\rm Id}}{-\pi}} &&&& &\\
& &{\begin{smallmatrix}1\\
1\end{smallmatrix}} \oplus {\begin{smallmatrix}0\\
1\end{smallmatrix}}\ar[rr]^-{\left(\begin{smallmatrix}\pi&{\rm
Id}\end{smallmatrix}\right)} & & {\begin{smallmatrix}0\\
1\end{smallmatrix}} & \notin {\rm smon}(Q, I, A).\\
{\begin{smallmatrix}1\\
1\end{smallmatrix}}\ar[urr]^-{\binom{{\rm Id}}{0}}&&&&&}$$ There are $58$ indecomposable $\m$-modules, and the
Auslander-Reiten quiver of $\m$ is given at the end of this paper, by an anticlockwise rotation,
where $\circ: = 0, \ \bullet: = 1$ and $\Box: = 2$.  For example $\begin{smallmatrix}\bullet &&&\\
\bullet&\Box&\circ\\ \circ &\bullet&\bullet \\ \circ & & &
\end{smallmatrix}$
means $\begin{smallmatrix}1&&&\\
1&2 &0\\ 0 &1&1 \\ 0 & & & \end{smallmatrix}$. From the
Auslander-Reiten quiver of $\m$ we get the Auslander-Reiten quiver
of ${\rm smon}(Q, I, A)$
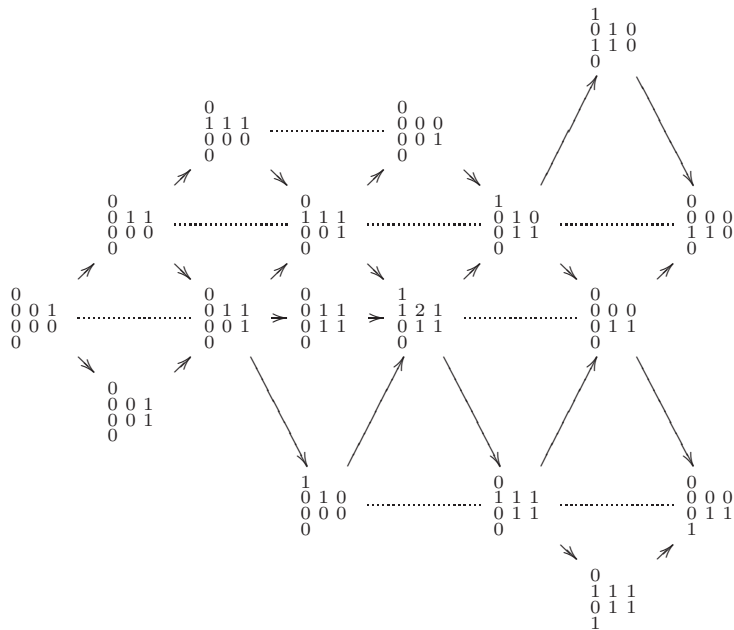
$$\xymatrix@R=0.2cm@C=0.2cm{ &&&& && {\begin{smallmatrix}1&&&\\ 0&1&0\\ 1 &1&0 \\ 0 & & & \end{smallmatrix}} \ar[ddr]
\\
&& {\begin{smallmatrix}0&&&\\ 1&1 &1\\ 0 &0&0 \\ 0 & & &
\end{smallmatrix}}\ar[dr] && {\begin{smallmatrix}0&&&\\ 0&0 &0\\ 0
&0&1 \\ 0 & & & \end{smallmatrix}}\ar[dr]\ar@{.}[ll]
\\
& {\begin{smallmatrix}0&&&\\0&1&1\\ 0 &0&0 \\ 0 & & &
\end{smallmatrix}}\ar[dr]\ar[ur] & & {\begin{smallmatrix}0&&&\\1&1
&1\\ 0 &0&1 \\ 0 & & & \end{smallmatrix}}\ar[dr]\ar[ur]\ar@{.}[ll]
&&{\begin{smallmatrix}1&&&\\ 0&1 &0\\ 0 &1&1 \\ 0 & & &
\end{smallmatrix}}\ar[dr]\ar[uur]\ar@{.}[ll] &&
{\begin{smallmatrix}0&&&\\ 0&0 &0\\ 1 &1&0 \\ 0 & & &
\end{smallmatrix}}\ar@{.}[ll]
\\
{\begin{smallmatrix}0&&&\\0&0 &1\\ 0 &0&0 \\ 0 & & & \end{smallmatrix}}\ar[ur]\ar[dr] &&{\begin{smallmatrix}0&&&\\0&1 &1\\ 0 &0&1 \\ 0 & & & \end{smallmatrix}}\ar[ur]\ar[ddr]\ar[r]\ar@{.}[ll] & {\begin{smallmatrix}0&&&\\
0&1 &1\\ 0 &1&1 \\ 0 & & & \end{smallmatrix}}\ar[r]& {\begin{smallmatrix}1&&&\\
1&2 &1\\ 0 &1&1 \\ 0 & & & \end{smallmatrix}}\ar[ddr]\ar[ur]&& {\begin{smallmatrix}0&&&\\
0&0 &0\\ 0 &1&1 \\ 0 & & & \end{smallmatrix}}\ar[ddr]
\ar[ur]\ar@{.}[ll]
\\
& {\begin{smallmatrix}0&&&\\0&0 &1\\ 0 &0&1 \\ 0 & & &
\end{smallmatrix}}\ar[ur]
\\ & &&
{\begin{smallmatrix}1&&&\\
0&1 &0\\ 0 &0&0\\ 0 & & & \end{smallmatrix}}\ar[uur]&&
{\begin{smallmatrix}0&&&\\1&1 &1\\ 0 &1&1 \\ 0 & & & \end{smallmatrix}}\ar[uur]\ar[dr]\ar@{.}[ll]&& {\begin{smallmatrix}0&&&\\
0&0 &0\\ 0 &1&1\\ 1 & & & \end{smallmatrix}}\ar@{.}[ll]
\\
&&&&&& {\begin{smallmatrix}0&&&\\
1&1 &1\\ 0 &1&1 \\ 1 & & & \end{smallmatrix}\ar[ur]}}$$ with
$17$ indecomposable objects in ${\rm smon}(Q, I, A)$, but only $12$
indecomposable objects in $A\mbox{-}{\rm mod}\otimes
(kQ/I)\mbox{-}{\rm mod}$.

\vskip5pt

A $\m$-module $Y\in {\rm sepi}(Q, I, A)$ if and only if
$Y_\alpha$ and $Y_\gamma$ are
epic and ${\rm Im} {Y_\beta} = {\rm Ker} {Y_\alpha}$.  By the Auslander-Reiten quiver of $\m$ we get the
Auslander-Reiten quiver of ${\rm sepi}(Q, I, A)$:
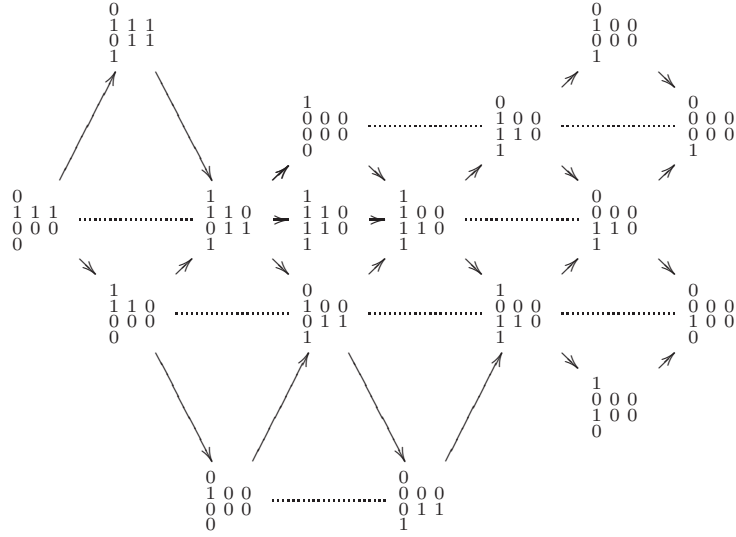
$$\xymatrix@R=0.2cm@C=0.2cm{&{\begin{smallmatrix}0&&&\\1&1 &1\\ 0 &1&1 \\ 1& & & \end{smallmatrix}}\ar[ddr]&&&&&
{\begin{smallmatrix}0&&&\\ 1&0&0\\ 0 &0&0 \\ 1& & &
\end{smallmatrix}} \ar[dr]
\\
&&& {\begin{smallmatrix}1&&&\\ 0&0&0\\ 0&0&0 \\ 0& & &
\end{smallmatrix}}\ar[dr]& & {\begin{smallmatrix}0&&&\\ 1&0&0\\ 1
&1&0 \\ 1 & & & \end{smallmatrix}}\ar[ur]\ar[dr]\ar@{.}[ll]& &
{\begin{smallmatrix}0&&&\\ 0&0&0\\ 0&0&0 \\ 1 & & &
\end{smallmatrix}}\ar@{.}[ll]
\\
{\begin{smallmatrix}0&&&\\ 1&1 &1\\ 0 &0&0 \\ 0 & & &
\end{smallmatrix}}\ar[dr]\ar[uur] & & {\begin{smallmatrix}1&&&\\ 1&1
&0\\ 0 &1&1 \\ 1 & & &
\end{smallmatrix}}\ar[r]\ar[dr]\ar[ur]\ar[ur]\ar@{.}[ll]
&{\begin{smallmatrix}1&&&\\ 1&1&0\\ 1 &1&0 \\ 1 & & &
\end{smallmatrix}}\ar[r]& {\begin{smallmatrix}1&&&\\ 1&0 &0\\ 1 &1&0
\\ 1 & & & \end{smallmatrix}}\ar[dr]\ar[ur] &&
{\begin{smallmatrix}0&&&\\ 0&0 &0\\ 1&1&0 \\ 1 & & &
\end{smallmatrix}}\ar[dr]\ar[ur]\ar@{.}[ll]
\\
&{\begin{smallmatrix}1&&&\\ 1&1 &0\\ 0 &0&0 \\ 0 & & &
\end{smallmatrix}}\ar[ur]\ar[ddr] && {\begin{smallmatrix}0&&&\\ 1&0
&0\\ 0 &1&1 \\ 1& & & \end{smallmatrix}}\ar[ur]\ar[ddr]\ar@{.}[ll]&
& {\begin{smallmatrix}1&&&\\ 0&0 &0\\ 1 &1&0 \\ 1 & & &
\end{smallmatrix}}\ar[dr]\ar[ur]\ar@{.}[ll] &&
{\begin{smallmatrix}0&&&\\ 0&0 &0\\ 1 &0&0 \\ 0 & & &
\end{smallmatrix}}\ar@{.}[ll]
\\
& & & && & {\begin{smallmatrix}1&&&\\ 0&0 &0\\ 1&0&0 \\ 0& & &
\end{smallmatrix}}\ar[ur]
\\
& & {\begin{smallmatrix}0&&&\\ 1&0 &0\\ 0 &0&0\\ 0 & & &
\end{smallmatrix}}\ar[uur]&& {\begin{smallmatrix}0&&&\\ 0&0 &0\\
0&1&1 \\ 1 & & & \end{smallmatrix}}\ar[uur]\ar@{.}[ll]}$$ Form the
Auslander-Reiten quivers of ${\rm smon}(Q, I, A)$ and ${\rm sepi}(Q,
I, A)$, we already see a {\bf unique} {\rm RSS} equivalence $F: {\rm
smon}(Q, I, A)\cong {\rm sepi}(Q, I, A)$, given
by
\begin{align*} & \begin{smallmatrix}1\\
0\end{smallmatrix}\otimes P(1) = \begin{smallmatrix}0&&&\\ 0&0 &1\\
0 &0&0\\ 0 & & &
\end{smallmatrix}
\mapsto \begin{smallmatrix}0&&&\\ 1&1 &1\\ 0 &0&0\\ 0 & & &
\end{smallmatrix} = \begin{smallmatrix}1\\
0\end{smallmatrix}\otimes I(1), \ \ \ \ \
\begin{smallmatrix}1\\
0\end{smallmatrix}\otimes P(2) = \begin{smallmatrix}0&&&\\ 0&1 &1\\
0 &0&0\\ 0 & & &
\end{smallmatrix}
\mapsto \begin{smallmatrix}1&&&\\ 1&1 &0\\ 0 &0&0\\ 0 & & &
\end{smallmatrix} = \begin{smallmatrix}1\\
0\end{smallmatrix}\otimes I(2)
\\ &
\begin{smallmatrix}1\\
0\end{smallmatrix}\otimes P(4) = \begin{smallmatrix}0&&&\\ 1&1 &1\\
0 &0&0\\ 0 & & &
\end{smallmatrix}
\mapsto \begin{smallmatrix}0&&&\\ 1&0 &0\\ 0 &0&0\\ 0 & & &
\end{smallmatrix} = \begin{smallmatrix}1\\
0\end{smallmatrix}\otimes I(4)
\\ \\&
\begin{smallmatrix}1\\
1\end{smallmatrix}\otimes P(1) = \begin{smallmatrix}0&&&\\ 0&0 &1\\
0 &0&1\\ 0 & & &
\end{smallmatrix}
\mapsto \begin{smallmatrix}0&&&\\ 1&1 &1\\ 0 &1&1\\ 1 & & &
\end{smallmatrix}= \begin{smallmatrix}1\\
1\end{smallmatrix}\otimes I(1), \ \ \ \ \ \begin{smallmatrix}0&&&\\
0&1 &1\\ 0 &0&1\\ 0 & & &
\end{smallmatrix}
\mapsto \begin{smallmatrix}1&&&\\ 1&1 &0\\ 0 &1&1\\ 1 & & &
\end{smallmatrix}, \ \ \ \ \
\begin{smallmatrix}0&&&\\ 1&1 &1\\ 0 &0&1\\ 0 & & &
\end{smallmatrix}
\mapsto \begin{smallmatrix}0&&&\\ 1&0 &0\\ 0 &1&1\\ 1 & & &
\end{smallmatrix}
\\&
\begin{smallmatrix}0\\
1\end{smallmatrix}\otimes P(1) = \begin{smallmatrix}0&&&\\ 0&0 &0\\
0 &0&1\\ 0 & & &
\end{smallmatrix}
\mapsto \begin{smallmatrix}0&&&\\ 0&0 &0\\ 0 &1&1\\ 1 & & &
\end{smallmatrix}= \begin{smallmatrix}0\\
1\end{smallmatrix}\otimes I(1)
\\ \\ &\begin{smallmatrix}1\\
1\end{smallmatrix}\otimes P(2) = \begin{smallmatrix}0&&&\\ 0&1 &1\\
0 &1&1\\ 0 & & &
\end{smallmatrix}
\mapsto \begin{smallmatrix}1&&&\\ 1&1 &0\\ 1 &1&0\\ 1 & & &
\end{smallmatrix}= \begin{smallmatrix}1\\
1\end{smallmatrix}\otimes I(2)
\\ \\ &
\begin{smallmatrix}1\\
0\end{smallmatrix}\otimes P(3) = \begin{smallmatrix}1&&&\\ 0&1 &0\\
0 &0&0\\ 0 & & &
\end{smallmatrix}
\mapsto \begin{smallmatrix}1&&&\\ 0&0 &0\\ 0 &0&0\\ 0 & & &
\end{smallmatrix} =
\begin{smallmatrix}1\\
0\end{smallmatrix}\otimes I(3), \ \ \ \ \
\begin{smallmatrix}1&&&\\ 1&2 &1\\ 0 &1&1\\ 0 & & &
\end{smallmatrix}
\mapsto \begin{smallmatrix}1&&&\\ 1&0 &0\\ 1 &1&0\\ 1 & & &
\end{smallmatrix}, \ \ \ \ \
\begin{smallmatrix}1&&&\\ 0&1 &0\\ 0 &1&1\\ 0 & & &
\end{smallmatrix}
\mapsto \begin{smallmatrix}1&&&\\ 0&0 &0\\ 1 &1&0\\ 1 & & &
\end{smallmatrix}
\\& \begin{smallmatrix}1\\
1\end{smallmatrix}\otimes P(3) = \begin{smallmatrix}1&&&\\ 0&1 &0\\ 1 &1&0\\
0 & & &
\end{smallmatrix}
\mapsto \begin{smallmatrix}1&&&\\ 0&0 &0\\ 1&0&0\\ 0 & & &
\end{smallmatrix}=\begin{smallmatrix}1\\
1\end{smallmatrix}\otimes I(3)\\ \\ &
\begin{smallmatrix}0&&&\\ 1&1 &1\\ 0 &1&1\\ 0 & & &
\end{smallmatrix}
\mapsto \begin{smallmatrix}0&&&\\ 1&0 &0\\ 1 &1&0\\ 1 & & &
\end{smallmatrix}, \ \ \ \ \
\begin{smallmatrix}0\\
1\end{smallmatrix}\otimes P(2)= \begin{smallmatrix}0&&&\\ 0&0 &0\\ 0
&1&1\\ 0 & & &
\end{smallmatrix}
\mapsto \begin{smallmatrix}0&&&\\ 0&0 &0\\ 1 &1&0\\ 1 & & &
\end{smallmatrix} = \begin{smallmatrix}0\\
1\end{smallmatrix}\otimes I(2), \ \ \ \ \ \begin{smallmatrix}0\\
1\end{smallmatrix}\otimes P(3) = \begin{smallmatrix}0&&&\\
0&0 &0\\ 1&1&0\\ 0 & & &
\end{smallmatrix}
\mapsto \begin{smallmatrix}0&&&\\ 0&0 &0\\ 1 &0&0\\ 0 & & &
\end{smallmatrix}= \begin{smallmatrix}0\\
1\end{smallmatrix}\otimes I(3)\end{align*}
\begin{align*} &
\begin{smallmatrix}1\\
1\end{smallmatrix}\otimes P(4) = \begin{smallmatrix}0&&&\\ 1&1 &1\\
0 &1&1\\ 1 & & &
\end{smallmatrix}
\mapsto \begin{smallmatrix}0&&&\\ 1&0 &0\\ 0 &0&0\\ 1 & & &
\end{smallmatrix} = \begin{smallmatrix} 1\\
1\end{smallmatrix}\otimes I(4), \ \ \ \ \ \begin{smallmatrix}0\\
1\end{smallmatrix}\otimes P(4) = \begin{smallmatrix}0&&&\\
0&0 &0\\ 0 &1&1\\ 1 & & &
\end{smallmatrix}
\mapsto \begin{smallmatrix}0&&&\\ 0&0 &0\\ 0 &0&0\\ 1 & & &
\end{smallmatrix} =\begin{smallmatrix}0\\
1\end{smallmatrix}\otimes I(4).
\end{align*}
\noindent  Since $F$  sends the
projective $\m$-module $\begin{smallmatrix}0&&&\\ 0&0 &1\\ 0 &0&0\\
0 & & &
\end{smallmatrix}$
to $\begin{smallmatrix}0&&&\\ 1&1 &1\\ 0 &0&0\\ 0 & & &
\end{smallmatrix}$ which is not
an injective $\m$-module, $F$ is {\bf not} the Nakayama functor
$\mathcal N_\m$ (this also follows from
Corollary \ref{NakayamaandRSS}).

\vskip5pt

Now, let $A$ be an arbitrary finite-dimensional algebra,  and $\x$ an arbitrary additive subcategory of $A$-{\rm mod} (not necessarily extension-closed). For  $X\in {\rm smon}(Q, I, A)$,
note that $X\in {\rm smon}(Q, I, \x)$ if and only if ${\rm Coker} X_\alpha\in\x, \ X_2/({\rm Im}X_\beta\oplus {\rm Im}X_\gamma)\in\x, \  X_3\in \x$ and $X_4\in \x$.
For  $Y\in {\rm sepi}(Q, I, A)$,
$Y\in {\rm sepi}(Q, I, \x)$ if and only if $Y_1\in\x, \ {\rm Ker} Y_\alpha\in\x,  \  {\rm Ker} Y_\beta\in\x$ and ${\rm Ker} Y_\gamma\in\x$.

\vskip5pt

We have a functor $F: {\rm smon}(Q, I, \x)\longrightarrow{\rm sepi}(Q, I, \x)$, which
sends $X\in {\rm smon}(Q, I, \x)$ to
$$\xymatrix@R=0.2cm@C=0.2cm{X_2/{\rm Im}X_\gamma\ar[drr]^-{Y_\beta} &&&& \\
& &X_1/X_\alpha({\rm Im}X_\gamma)\ar[rr]^-{Y_\alpha}& & X_1/{\rm Im}X_\alpha\\
X_1\ar[urr]^-{Y_\gamma}&&&&}$$ where $Y_\alpha$ is the canonical
epimorphism with kernel ${\rm Im}X_\alpha/X_\alpha({\rm
Im}X_\gamma),$ $Y_\beta$ is induced by $X_\alpha$, and $Y_\gamma$ is
the canonical epimorphism. Note that ${\rm Ker} Y_\alpha\cong X_2/({\rm Im}X_\beta\oplus {\rm Im}X_\gamma)\in\x$,
${\rm Ker} Y_\beta\cong {\rm Im}X_\beta \cong X_3\in\x$, and ${\rm Ker} Y_\gamma \cong X_4\in\x$.

\vskip5pt

Consider a functor  $G: {\rm sepi}(Q, I,
\x)\longrightarrow {\rm smon}(Q, I, \x)$, which sends $Y\in {\rm
sepi}(Q, I, \x)$ to
$$\xymatrix@R=0.2cm@C=0.2cm{{\rm Ker}Y_\beta\oplus 0\ar[drr]^-{X_\beta} &&&& \\
& &{\rm Ker}(Y_\beta, Y_\gamma)\ar[rr]^-{X_\alpha}& & Y_4\\
0\oplus {\rm Ker}Y_\gamma\ar[urr]^-{X_\gamma}&&&&}$$ where ${\rm
Ker}(Y_\beta, Y_\gamma)= {\rm Ker}(Y_3\oplus Y_4\stackrel{(Y_\beta,
Y_\gamma)}\longrightarrow Y_2)$ is in fact the pullback of $Y_\beta$
and $Y_\gamma$, $X_\alpha$ is the projection to $Y_4$, $X_\beta$
and $X_\gamma$ are embeddings. Note that ${\rm Coker} X_\alpha = Y_4/{\rm Ker}(Y_\alpha Y_\gamma)\cong Y_1\in\x$ and $X_2/({\rm Im}X_\beta\oplus {\rm Im}X_\gamma) = {\rm Ker}(Y_\beta, Y_\gamma) / ({\rm Ker} Y_\beta \oplus {\rm Ker} Y_\gamma) \cong {\rm Ker} Y_\alpha \in\x$.

\vskip5pt

Then $G$ is a quasi-inverse
of $F$: we omit the details, but note that the restriction of
$X_\alpha$ to ${\rm Im}X_\gamma$ is a monomorphism, since ${\rm
Im}X_\gamma \cap {\rm Ker}X_\alpha = {\rm Im}X_\gamma \cap {\rm
Im}X_\beta = 0$, and that $X_\alpha({\rm Ker}(Y_\beta, Y_\gamma))
= {\rm Ker}(Y_\alpha Y_\gamma)$ for $Y\in {\rm sepi}(Q, I, \x)$.
For $M\in \x$ we have functorial isomorphisms of left
$\m$-modules:
\begin{align*}& F(M\otimes P(1)) = F(\begin{smallmatrix}0\searrow&&&\\ && 0& \rightarrow &M\\ 0\nearrow &&&
\end{smallmatrix}
)  = (\begin{smallmatrix}0\searrow&&&\\ && M& \rightarrow &M\\
M\nearrow &&&
\end{smallmatrix})\cong M\otimes
I(1)\\ & F(M\otimes P(2)) \cong F(\begin{smallmatrix}0\searrow&&&\\
&& M& \rightarrow &M\\ 0\nearrow &&&
\end{smallmatrix}) =
(\begin{smallmatrix}M\searrow&&&\\ && M& \rightarrow &0\\ M\nearrow
&&&\end{smallmatrix}) \cong M\otimes I(2)
\\ & F(M\otimes P(3))\cong F(\begin{smallmatrix}M\searrow &&&\\ && M& \rightarrow &0\\ 0\nearrow &&&
\end{smallmatrix})  =
(\begin{smallmatrix}M\searrow&&&\\ && 0& \rightarrow &0\\ 0\nearrow
&&&\end{smallmatrix})\cong M\otimes I(3)\\ & F(M\otimes P(4))\cong
F(\begin{smallmatrix}0\searrow&&&\\ && M& \rightarrow &M\\ M\nearrow
&&&\end{smallmatrix}) = (\begin{smallmatrix}0\searrow&&&\\ && 0&
\rightarrow &0\\ M\nearrow &&&
\end{smallmatrix})
\cong M\otimes I(4).
\end{align*}
So $F: {\rm smon}(Q, I, \x)\cong {\rm sepi}(Q, I, \x)$ is
an {\rm RSS} equivalence.

\vskip 15pt

{\bf Acknowledgements.} We are grateful to Claus Michael Ringel for
the helpful discussions,  especially for suggesting the terminology
``separated monic module".

\newpage
$$\ \ \ \ \ \ \ \ \ \ \ \ \ \ \  \xymatrix@R=0.2cm@C=0.3cm{&&&&& {\begin{smallmatrix}\circ&\circ&\circ&\circ\\
&\circ&\circ&\\ &\circ&\bullet& \end{smallmatrix}}\ar[dl]\ar[dr]
\\
&&&&{\begin{smallmatrix}\circ&\circ&\circ&\circ\\
&\circ&\circ&\\ &\bullet&\bullet&
\end{smallmatrix}}\ar[dr]&&
{\begin{smallmatrix}\circ&\circ&\circ&\circ\\
&\circ&\bullet&\\ &\circ&\bullet&\end{smallmatrix}}\ar[dl]\ar[dr]
\\
&&&&& {\begin{smallmatrix}\circ&\circ&\circ&\circ\\
&\circ&\bullet&\\ &\bullet&\bullet&
\end{smallmatrix}}\ar[dl]\ar[dr]\ar[d]\ar@{.}[uu] && {\begin{smallmatrix}\circ&\circ&\bullet&\circ\\
&\circ&\bullet&\\ &\circ&\bullet&
\end{smallmatrix}}\ar[dl]
&\\
&&&&{\begin{smallmatrix}\circ&\circ&\circ&\circ\\
&\circ&\bullet&\\ &\circ&\circ&
\end{smallmatrix}}\ar[dl]\ar[dr]\ar@{.}[uu] & {\begin{smallmatrix}\circ&\circ&\circ&\circ\\
&\bullet&\bullet&\\ &\bullet&\bullet&\end{smallmatrix}}\ar[d]& {\begin{smallmatrix}\circ&\circ&\bullet&\circ\\
&\circ&\bullet&\\
&\bullet&\bullet&\end{smallmatrix}}\ar[dl]\ar[dr]\ar@{.}[uu]
\\
&&&{\begin{smallmatrix}\circ&\circ&\circ&\bullet\\
&\circ&\bullet&\\ &\circ&\circ&
\end{smallmatrix}}\ar[dr]&& {\begin{smallmatrix}\circ&\circ&\bullet&\circ\\
&\bullet&\Box&\\
&\bullet&\bullet&\end{smallmatrix}}\ar[dl]\ar[d]\ar[dr]&&
{\begin{smallmatrix}\circ&\circ&\circ&\circ\\
&\circ&\circ&\\ &\bullet&\circ&\end{smallmatrix}}\ar[dl]\ar@{.}[uu]
\\
&&&&{\begin{smallmatrix}\circ&\circ&\bullet&\bullet\\
&\bullet&\Box&\\ &\bullet&\bullet&\end{smallmatrix}}
\ar[dl] \ar[dr]\ar@{.}[uu]& {\begin{smallmatrix}\circ&\circ&\bullet&\circ\\
&\circ&\bullet&\\ &\circ&\circ&
\end{smallmatrix}}\ar[d] & {\begin{smallmatrix}\circ&\circ&\circ&\circ\\
&\bullet&\bullet&\\ &\bullet&\circ&
\end{smallmatrix}}\ar[dl]\ar[dr]\ar@{.}[uu]
\\
&&&{\begin{smallmatrix}\circ&\circ&\bullet&\circ\\
&\bullet&\bullet&\\ &\bullet&\bullet&
\end{smallmatrix}}\ar[dr]\ar[dl]\ar@{.}[uu] && {\begin{smallmatrix}\circ&\circ&\bullet&\bullet\\
&\bullet&\Box&\\ &\bullet&\circ&\end{smallmatrix}}\ar[dl]\ar[dr]\ar[d]& &{\begin{smallmatrix}\circ&\circ&\circ&\circ\\
&\bullet&\bullet&\\
&\circ&\circ&\end{smallmatrix}}\ar[dl]\ar@{.}[uu]
\\
&&{\begin{smallmatrix}\bullet&\circ&\bullet&\circ\\
&\bullet&\bullet &\\ &\bullet&\bullet&
\end{smallmatrix}}\ar[dr]&&{\begin{smallmatrix}\circ&\circ&\bullet&\circ\\
&\bullet&\bullet &\\ &\bullet&\circ&
\end{smallmatrix}}\ar[dr]\ar[dl]\ar@{.}[uu] &{\begin{smallmatrix}\circ&\circ&\circ&\bullet\\
&\bullet&\bullet &\\ &\bullet&\circ&
\end{smallmatrix}}\ar[d]& {\begin{smallmatrix}\circ&\circ&\bullet&\bullet\\
&\bullet&\Box&\\ &\circ&\circ&
\end{smallmatrix}}\ar[dl]\ar[dr]\ar@{.}[uu]
\\ &&&{\begin{smallmatrix}\bullet&\circ&\bullet&\circ\\
&\bullet&\bullet&\\ &\bullet&\circ&
\end{smallmatrix}}\ar[dr]\ar@{.}[uu]& &{\begin{smallmatrix}\circ&\circ&\bullet&\bullet\\
&\Box&\Box&\\
&\bullet&\circ&\end{smallmatrix}}\ar[d]\ar[dr]\ar[dl]&&
{\begin{smallmatrix}\circ&\circ&\bullet&\bullet\\
&\circ&\bullet&\\ &\circ&\circ&\end{smallmatrix}} \ar[dl]\ar@{.}[uu]
\\
&&&&{\begin{smallmatrix}\bullet&\circ&\bullet&\bullet\\
&\Box&\Box&\\ &\bullet&\circ&
\end{smallmatrix}}\ar@{.}[uu]\ar[dr]\ar[dl]&{\begin{smallmatrix}\circ&\circ&\bullet&\circ\\
&\bullet&\bullet&\\ &\circ&\circ&
\end{smallmatrix}}\ar[d]& {\begin{smallmatrix}\circ&\circ&\bullet&\bullet\\
&\bullet&\bullet&\\ &\bullet&\circ&
\end{smallmatrix}}\ar[dr]\ar[dl]\ar@{.}[uu]
\\
&&&{\begin{smallmatrix}\circ&\circ&\circ&\bullet\\
&\bullet&\bullet&\\ &\circ&\circ&
\end{smallmatrix}}\ar[dr]\ar[dl]\ar@{.}[uu]& &{\begin{smallmatrix}\bullet&\circ&\Box&\bullet\\
&\Box&\Box&\\ &\bullet&\circ& \end{smallmatrix}}\ar[d]\ar[dl]\ar[dr]&& {\begin{smallmatrix}\circ&\circ&\circ&\circ\\
&\bullet&\circ&\\
&\bullet&\circ&\end{smallmatrix}}\ar@{.}[uu]\ar[dl]
\\
&&{\begin{smallmatrix}\circ&\bullet&\circ&\bullet\\
&\bullet&\bullet&\\ &\circ&\circ&
\end{smallmatrix}}\ar[dr]& &{\begin{smallmatrix}\circ&\circ&\bullet&\bullet\\
&\bullet&\bullet&\\ &\circ&\circ& \end{smallmatrix}}\ar[dr]\ar@{.}[uu]\ar[dl]&{\begin{smallmatrix}\bullet&\circ&\bullet&\bullet\\
&\bullet&\bullet&\\ &\bullet&\circ&
\end{smallmatrix}}\ar[d]& {\begin{smallmatrix}\bullet&\circ&\bullet&\circ\\
&\Box&\bullet&\\
&\bullet&\circ&\end{smallmatrix}}\ar[dr]\ar[dl]\ar@{.}[uu]
\\
&&&{\begin{smallmatrix}\circ&\bullet&\bullet&\bullet\\
&\bullet&\bullet&\\ &\circ&\circ&
\end{smallmatrix}}\ar[dr]\ar[dl]\ar@{.}[uu]& &{\begin{smallmatrix}\bullet&\circ&\bullet&\bullet\\
&\Box&\bullet&\\ &\bullet&\circ& \end{smallmatrix}}\ar[d]\ar[dr]\ar[dl]&&{\begin{smallmatrix}\bullet&\circ&\bullet&\circ\\
&\bullet&\bullet&\\ &\circ&\circ&
\end{smallmatrix}}\ar[dl]\ar@{.}[uu]
\\
&&{\begin{smallmatrix}\circ&\circ&\bullet&\circ\\
&\circ&\circ&\\ &\circ&\circ&
\end{smallmatrix}}\ar[dr]\ar@{.}[uu]& &{\begin{smallmatrix}\bullet&\bullet&\bullet&\bullet\\
&\Box&\bullet&\\ &\bullet&\circ& \end{smallmatrix}}\ar[dr]\ar[dl]\ar@{.}[uu]&{\begin{smallmatrix}\circ&\circ&\circ&\circ\\
&\bullet&\circ&\\ &\circ&\circ&
\end{smallmatrix}}\ar[d]& {\begin{smallmatrix}\bullet&\circ&\bullet&\bullet\\
&\bullet&\bullet&\\
&\circ&\circ&\end{smallmatrix}}\ar[dr]\ar[dl]\ar@{.}[uu]
\\
&&&{\begin{smallmatrix}\bullet&\circ&\bullet&\circ\\
&\bullet&\circ&\\ &\bullet&\circ&
\end{smallmatrix}}\ar[dl]\ar[dr]\ar@{.}[uu]& &{\begin{smallmatrix}\bullet&\bullet&\bullet&\bullet\\
&\Box&\bullet&\\ &\circ&\circ& \end{smallmatrix}}\ar[d]\ar[dr]\ar[dl]&&{\begin{smallmatrix}\circ&\circ&\circ&\bullet\\
&\circ&\circ&\\ &\circ&\circ&
\end{smallmatrix}}\ar[dl]\ar@{.}[uu]
\\
&&{\begin{smallmatrix}\bullet&\circ&\circ&\circ\\
&\bullet&\circ&\\ &\bullet&\circ&
\end{smallmatrix}}\ar[dr]\ar@{.}[uu]& &{\begin{smallmatrix}\bullet&\circ&\bullet&\circ\\
&\bullet&\circ&\\ &\circ&\circ& \end{smallmatrix}}\ar[dr]\ar[dl]\ar@{.}[uu]&{\begin{smallmatrix}\bullet&\bullet&\bullet&\bullet\\
&\bullet&\bullet&\\ &\circ&\circ&
\end{smallmatrix}}\ar[d]& {\begin{smallmatrix}\circ&\bullet&\circ&\bullet\\
&\bullet&\circ&\\
&\circ&\circ&\end{smallmatrix}}\ar[dr]\ar[dl]\ar@{.}[uu]
\\
&&&{\begin{smallmatrix}\bullet&\circ&\circ&\circ\\
&\bullet&\circ&\\ &\circ&\circ&
\end{smallmatrix}}\ar[dr]\ar@{.}[uu]& &{\begin{smallmatrix}\bullet&\bullet&\bullet&\bullet\\
&\bullet&\circ&\\ &\circ&\circ& \end{smallmatrix}}\ar[dr]\ar[dl]&&{\begin{smallmatrix}\circ&\bullet&\circ&\circ\\
&\bullet&\circ&\\ &\circ&\circ&
\end{smallmatrix}}\ar@{.}[uu]\ar[dl]
\\
&&&&{\begin{smallmatrix}\bullet&\bullet&\circ&\bullet\\
&\bullet&\circ&\\ &\circ&\circ&
\end{smallmatrix}}\ar[dl]\ar[dr]\ar@{.}[uu]& &{\begin{smallmatrix}\bullet&\bullet&\bullet&\circ\\
&\bullet&\circ&\\ &\circ&\circ&
\end{smallmatrix}}\ar[dr]\ar[dl]\ar@{.}[uu]
\\
&&&{\begin{smallmatrix}\circ&\bullet&\circ&\bullet\\
&\circ&\circ&\\ &\circ&\circ&
\end{smallmatrix}}\ar[dr]\ar@{.}[uu]& &{\begin{smallmatrix}\bullet&\bullet&\circ&\circ\\
&\bullet&\circ&\\ &\circ&\circ& \end{smallmatrix}}\ar[dl]\ar[dr]\ar@{.}[uu]&&{\begin{smallmatrix}\bullet&\circ&\bullet&\circ\\
&\circ&\circ&\\ &\circ&\circ&
\end{smallmatrix}}\ar[dl]\ar@{.}[uu]
\\
&&&&{\begin{smallmatrix}\circ&\bullet&\circ&\circ\\
&\circ&\circ&\\ &\circ&\circ&
\end{smallmatrix}}\ar@{.}[uu]& &{\begin{smallmatrix}\bullet&\circ&\circ&\circ\\
&\circ&\circ&\\ &\circ&\circ& \end{smallmatrix}}\ar@{.}[uu]}
$$

\vskip10pt

Pu Zhang

School of Mathematics,
\ Shanghai Jiao Tong University,  \ Shanghai 200240, \ China

\vskip5pt

Bao-Lin Xiong

Department of Mathematics,
\ Beijing University of Chemical Technology, Beijing 100029, \ China
\end{document}